%% file: fib.tex
\documentclass[preprint,11pt]{elsarticle}
\usepackage[bookmarks=true,colorlinks=true,linkcolor=blue]{hyperref}

\pdfoutput=1

\biboptions{sort&compress}

\usepackage{color}
\usepackage{amssymb,amsmath}
\usepackage{amsthm}

\usepackage{calc}
\usepackage[margin=1.in]{geometry}

\usepackage{algorithm} 
\usepackage{algpseudocode} 

\input texFiles/defs

\usepackage{tikz}

\input texFiles/trimFig.tex

\begin{document}

\small

\begin{frontmatter}
\title{A stable partitioned FSI algorithm for incompressible flow and elastic solids}

\author[llnl]{J. W. Banks\fnref{llnlThanks}}
\ead{banks20@llnl.gov}

\author[llnl]{W. D. Henshaw\corref{cor1}\fnref{llnlThanks}}
\ead{henshaw1@llnl.gov}

\author[rpi]{D. W. Schwendeman\fnref{donThanks}}
\ead{schwed@rpi.edu}

\address[llnl]{Centre for Applied Scientific Computing, Lawrence Livermore National Laboratory, Livermore, CA 94551, USA}

\address[rpi]{Department of Mathematical Sciences, Rensselaer Polytechnic Institute, Troy, NY 12180, USA}

\cortext[cor1]{Corresponding author. Mailing address: Centre for Applied Scientific Computing, L-422, Lawrence Livermore National Laboratory, Livermore, CA 94551, USA. Phone: 925-423-2697. Fax: 925-424-2477. }

\fntext[llnlThanks]{This work was performed under the auspices of the U.S. Department of Energy (DOE) by
  Lawrence Livermore National Laboratory under Contract DE-AC52-07NA27344 and by 
  DOE contracts from the ASCR Applied Math Program.}

\fntext[donThanks]{This research was supported by Lawrence Livermore National Laboratory under
subcontract B548468, and by the National Science Foundation under grant DMS-1016188.}

\begin{abstract}

A stable partitioned algorithm for coupling incompressible flows with
compressible elastic solids is described.  This added-mass partitioned (AMP)
scheme requires no sub-iterations, can be made fully second- or higher-order
accurate, and remains stable even in the presence of strong added-mass effects.
The approach extends the scheme of Banks et al.~\cite{fsi2012,sjogreenBanks2012}
for compressible flow, and uses Robin (mixed) boundary conditions with the fluid
and solid solvers at the interface.  The AMP Robin conditions are derived from a
local characteristic decomposition in the solid at the interface.  Two forms of
the Robin conditions are derived depending on whether the fluid equations are
advanced with a fractional-step method or not.  A normal mode analysis for a
discretization of an FSI model problem is performed to show that the new AMP
algorithm is stable for any ratio of the solid and fluid densities, including
the case of very light solids when the added-mass effects are large.  In
contrast, it is shown that a traditional partitioned algorithm involving a
Dirichlet-Neumann coupling for the same FSI problem is formally unconditionally
unstable for any ratio of densities.  Exact traveling wave solutions are derived
for three FSI model problems of increasing complexity, and these solutions are
used to verify the stability and accuracy of the corresponding numerical results
obtained from the AMP algorithm for the cases of light, medium and heavy solids.

\end{abstract}


\end{frontmatter}

\clearpage
\tableofcontents

\clearpage
\input texFiles/intro.tex

\input texFiles/governingEquations.tex


\input texFiles/partitionedAlgorithms

\input texFiles/modelProblems
\input texFiles/ampAlgorithm1d


\input texFiles/modelProblemAnalysis

\input texFiles/algorithm.tex


\input texFiles/results

\input texFiles/conclusions

\appendix

\input texFiles/exactSolutions

\input texFiles/procedures
\bibliographystyle{elsart-num}
\bibliography{journal-ISI,jwb,henshawPapers,henshaw,fsi}

\end{document}

%% file: texFiles/defs.tex
\renewcommand{\url}[1]{}
\newcommand{\citeCount}[1]{}
\newtheorem{theorem}{Theorem}

\newtheorem{proposition}{Proposition}

\newtheorem*{algorithmAMP}{AMP Algorithm (added-mass partitioned) }
\newtheorem*{algorithmTP}{TP Algorithm (traditional partitioned)}

\newtheorem*{ModelProblemIA}{Model Problem MP-IA}
\newtheorem*{ModelProblemVA}{Model Problem MP-VA}
\newtheorem*{ModelProblemVE}{Model Problem MP-VE}

\newcommand{\dx}{\Delta x}
\newcommand{\dy}{\Delta y}
\newcommand{\dt}{\Delta t}

\newcommand{\bogus}[1]{{}}

\newcommand{\ev}{\mathbf{ e}}

\newcommand{\iv}{\mathbf{ i}}

\newcommand{\nv}{\mathbf{ n}}

\newcommand{\qv}{\mathbf{ q}}

\newcommand{\tv}{\mathbf{ t}}
\newcommand{\uv}{\mathbf{ u}}
\newcommand{\vv}{\mathbf{ v}}

\newcommand{\xv}{\mathbf{ x}}

\newcommand{\Iv}{\mathbf{ I}}

\newcommand{\Lv}{\mathbf{ L}}

\newcommand{\half}{{1\over2}}

\newcommand{\Real}{{\mathbb R}}
\newcommand{\Complex}{{\mathbb C}}

\newcommand{\Ac}{{\mathcal A}}
\newcommand{\Bc}{{\mathcal B}}
\newcommand{\Dc}{{\mathcal D}}

\newcommand{\Oc}{{\mathcal O}}

\newcommand{\tauv}{\boldsymbol{\tau}}

\newcommand{\sigmav}{\boldsymbol{\sigma}}

\newcommand{\grad}{\nabla}

\newcommand{\tableFont}{\scriptsize}
\newcommand{\tableFontSize}{\small}
\newcommand{\num}[2]{#1e#2} 

\newcommand{\rateLabel}{rate}

\clearpage

\newcommand{\rhos}{\bar{\rho}}
\newcommand{\cp}{\bar{c}_p}
\newcommand{\cs}{\bar{c}_s}

\newcommand{\zs}{{\bar {z}}}  
\newcommand{\zp}{z_p}  
\newcommand{\us}{\bar{u}}

\newcommand{\vs}{\bar{v}}
\newcommand{\sigmas}{\bar{\sigma}}

\newcommand{\zf}{z_f}

\newcommand{\usv}{\bar{\uv}}
\newcommand{\vsv}{\bar{\vv}}

\newcommand{\ush}{\hat{\us}}
\newcommand{\usvh}{\hat{\usv}}
\newcommand{\vvh}{\hat{\vv}}
\newcommand{\vh}{\hat{v}}
\newcommand{\ph}{\hat{p}}

\newcommand{\vvs}{\bar{\vv}}
\newcommand{\uvs}{\bar{\uv}}

\newcommand{\sigmavs}{\bar{\sigmav}}
\newcommand{\sigmasv}{\bar{\sigmav}}

\newcommand{\qsv}{\bar{\qv}}

\newcommand{\vvI}{\vv^I}

\newcommand{\vvIf}{\vv^f}

\newcommand{\amp}{A}

\newcommand{\bfss}{\sffamily\bfseries}

\newcommand{\normalss}{\sffamily}

\newcommand{\Hs}{\bar{H}}

\newcommand{\aHat}{\hat{a}}
\newcommand{\bHat}{\hat{b}}

\newcommand{\aTilde}{\tilde{a}}
\newcommand{\bTilde}{\tilde{b}}

\newcommand{\dHat}{\hat{d}}

\newcommand{\qHat}{\hat{q}}
\newcommand{\vHat}{\hat{v}}

\newcommand{\sigmaHat}{\hat{\sigma}}

\newcommand{\lambdas}{\bar{\lambda}}
\newcommand{\mus}{\bar{\mu}}

\newcommand{\OmegaF}{\Omega^F}
\newcommand{\OmegaS}{\Omega^S}

\newcommand{\OmegaFh}{\Omega^F_h}
\newcommand{\OmegaSh}{\Omega^S_h}

\newcommand{\GammaI}{\Gamma}
\newcommand{\GammaIh}{\Gamma_h}

\renewcommand{\zp}{\bar{z}_p}
\renewcommand{\zs}{\bar{z}_s}
\newcommand{\zb}{\bar{z}}

\newcommand{\MR}{M_r}
\newcommand{\Mr}{\mathcal{M}}

\newcommand{\phis}{\phi_s}
\newcommand{\phib}{\phi_b}

\newcommand{\phim}{\phi_-}
\newcommand{\phip}{\phi_+}

\newcommand{\kk}{k_x} 
\newcommand{\kx}{k_x} 

\newcommand{\thetaj}{\theta} 

\newcommand{\Achar}{\Ac}
\newcommand{\Bchar}{\Bc}

\newcommand{\Af}{A_f}
\newcommand{\Bf}{B_f}
\newcommand{\Df}{\Dc}
\newcommand{\As}{A_s}
\newcommand{\Bs}{B_s}
\newcommand{\Ds}{\bar{\Dc}}

\newcommand{\Ck}{C_k}
\newcommand{\Sk}{S_k}
\newcommand{\Ca}{C_\alpha}
\newcommand{\Sa}{S_\alpha}

\newcommand{\Sas}{S_a}
\newcommand{\Cas}{C_a}
\newcommand{\Sbs}{S_b}
\newcommand{\Cbs}{C_b}

\newcommand{\dsf}{\delta}

\newcommand{\expkw}{e^{i(k x -\omega t)}}
\newcommand{\ampe}{{\bar u_{{\rm max}}}}

\newcommand{\ly}{\lambda_y}
\newcommand{\lx}{\lambda_x}

\newcommand{\tanv}{\ev}
\newcommand{\blue}{\color{blue}}

\newcommand{\strutt}{\rule{0pt}{10pt}}

%% file: texFiles/trimFig.tex
\newlength{\tfwidth}
\newlength{\tfheight}
\newlength{\tfxa}
\newlength{\tfxb}
\newlength{\tfya}
\newlength{\tfyb}
\newcommand{\trimFigNoBox}[6]{%
\setlength\fboxsep{1pt}
\setlength\fboxrule{0.0pt}
\fbox{\includegraphics[width=#2, clip, trim=#3 #4 #5 #6]{#1}}%
}

\newcommand{\trimFig}[6]{%
\setlength{\tfwidth}{(#2+#2*\real{#3})+#2*\real{#4}}
\setlength{\tfheight}{(#2+#2*\real{#5})+#2*\real{#6}}%
\setlength{\tfxa}{\tfwidth*\real{#3}}%
\setlength{\tfxb}{\tfwidth*\real{#4}}%
\setlength{\tfya}{\tfheight*\real{#5}}%
\setlength{\tfyb}{\tfheight*\real{#6}}%
\trimFigNoBox{#1}{#2}{\tfxa}{\tfya}{\tfxb}{\tfyb}%
}

%% file: texFiles/intro.tex
\section{Introduction}

Fluid-structure interaction problems, which involve the coupled evolution of fluids and structures,
are of great interest due to their importance in a wide range of applications such as those found 
in structural engineering and biomedicine. Numerical algorithms used to solve
FSI problems can be classified into two main categories. In monolithic schemes, the
fluid and structural equations are combined into a single large system of equations.
In partitioned schemes (also known as modular or sequential schemes) there are separate fluid and structural solvers
which are coupled at the interface. 
Strongly coupled partitioned schemes perform multiple iterations per time-step to solve the
coupled equations while loosely coupled schemes use only one or a few iterations.
Partitioned schemes are often preferred since they can make use of existing solvers and can be more
efficient than monolithic schemes. Partitioned schemes, however, can have stability issues in certain regimes as
will be discussed further below.

In this manuscript we consider an important class of FSI problems where partitioned schemes often 
have difficulty. This class involves the
coupling of incompressible fluids and relatively {\em light} structures where the {\em added-mass} effects are large.
The added-mass effect arises from the fact that the force required to accelerate a solid structure immersed in a fluid
must also account for the acceleration of the surrounding fluid.
Partitioned schemes for this class of problem generally suffer from poor convergence, requiring many sub-iterations per time step,
or are unstable, when the added-mass effects are large. 
The instabilities in partitioned schemes can be traced to the application of the interface conditions which must also
be partitioned between the fluid and solid solvers.
The traditional partitioned algorithm uses the velocity of the solid as a boundary condition on the fluid
and the traction from the fluid as a forcing on the solid. This is sometimes referred to as Dirichlet-Neumann
coupling with the Dirichlet condition being the assignment of the velocity and the Neumann condition the assignment
of the traction. More generally practitioners have considered Robin-Robin coupling (i.e. {\em mixed}
boundary conditions for the fluid and solid) to impose the interface conditions~\cite{BadiaNobileVergara2008}.

In this work we develop a stable partitioned algorithm for coupling
incompressible flows with compressible elastic {\em bulk} solids
Our new added-mass
partitioned (AMP) algorithm is stable, even when added-mass effects are large, and
requires no sub-iterations to converge.  The coefficients in the AMP Robin
interface conditions are derived from a local characteristic decomposition in
the solid, and there is no need to estimate any acceleration parameters as with
some other approaches.  Two forms of the Robin interface conditions are provided
depending on whether the fluid equations are advanced with a fractional-step
method or not. 
We develop and evaluate the scheme
for a linearized FSI problem that couples Stokes flow in the fluid to a linear
elastic solid on fixed fluid and solid reference domains. A normal mode analysis
of a two-dimensional model problem shows that the AMP algorithm is stable for any
ratio of the fluid to solid densities. The analysis also shows that a scheme
based on the traditional Dirichlet-Neumann coupling is {\em formally} unconditionally unstable
for any density ratio in the sense that the scheme may be stable on a coarse grid but
becomes unstable for a sufficiently fine mesh\footnote{Formally, 
a scheme is said to be stable, if for a given time interval, the numerical solutions can be bounded, independent of the
mesh spacing, for 
all mesh spacings sufficiently small~\cite{GustafssonKreissOliger95}. We will sometimes informally refer to a scheme 
being stable for some particular values of the mesh spacing and time-step.}.
The incompressible flow
equations in our implementation of the AMP algorithm are solved with a
fractional-step method based on a velocity-pressure
formulation~\cite{ICNS,splitStep2003} that can be made fully second-order (or
higher-order) accurate.
A key ingredient to this scheme is the proper specification of boundary conditions for the pressure.  
The elastic wave
equation is solved with a second-order accurate upwind scheme for the equations
written as a first-order system following the approach developed for
overlapping grids in Appel\"o et~al.~\cite{smog2012}.  The AMP approximation can
be made fully second-order (or high-order) accurate.  We develop exact traveling
wave solutions to three FSI model problems.  Numerical results using these exact
solutions confirm that our formally second-order accurate implementation of the
AMP algorithm is stable and second-order accurate in the maximum norm. 

In previous work by Banks et~al.~\cite{fsi2012},
a stable partitioned algorithm for {\em compressible} flows and compressible elastic solids was developed that overcomes the
added-mass effect and uses deforming composite grids
to treat large solid motions.
Compressible flows are somewhat easier to deal with than incompressible flows since
the added-mass effects are more localized
due to the finite speeds of wave propagation in the fluid~\cite{vanBrummelen2009}. 
The scheme in~\cite{fsi2012}, which was based on the analysis
of Banks and Sj\"ogreen~\cite{sjogreenBanks2012}, 
used a local characteristic analysis
of a fluid-structure Riemann problem to define  an impedance weighted averaging of the fluid and
solid interface values (i.e.~a Robin-Robin coupling). 
This algorithm was extended to the case of rigid bodies
in~\cite{lrb2013} where it was shown that the scheme remains stable even for rigid bodies of zero mass.

In recent work~\cite{fis2013r}, we have also developed a stable partitioned FSI algorithm for
incompressible flow coupled to structural shells (or beams) that remains
stable even for light structures when added mass effects are large. The case of shells
or beams is different from the situation of bulk solids since the fluid traction 
enters as a body forcing in the structural equations, instead of as a boundary condition. 
The AMP algorithm for shells is thus different than the one for bulk solids
presented in this article, although both algorithms impose generalized Robin boundary conditions
on the fluid domain.

FSI is a very active field of
research with many publications, see for example~\cite{DoneaGiulianiHalleux1982,CebralLohner1997,Lohner99,LeTallecMouro2001,PipernoFarhat2001,SchaferTeschauer2001,GlowinskiPanHelsaJosephPeriaux2001,Guruswamy2002,MillerColella2002,KuhlHulshoffBorst2003,ArientiHungMoranoShepherd2003,Farhat2006,HronTurekMonolithic2006,Tezduyar2006,SchaferHeckYigit2006,AhnKallinderis2006,CirakDeiterdingMauch2007,vanLoon2007,TezduyarSathe2007,BorazjaniGeSotiropolous2008,Wilcox2010,BartonObadiaDrikakis2011,DegrooteVierendeels2011,TezduyarTakizawaBrummerChen2011,CrosettoDeparisFouresteyQuarteroni2011,Degroote2011,GretarssonKwatraFedkiw2011,HouWangLayton2012}. There has also
been much work concerning partitioned algorithms and the added-mass effect, and we outline some of this work now.
For the case of incompressible fluids coupled to thin structural {\em shells} there has been some success in 
treating the added-mass instability. 
The  {\em kinematically coupled} 
scheme of Guidoboni et~al.~\cite{GuidoboniGlowinskiCavalliniCanic2009}, later extended to the
$\beta$-scheme by \v{C}ani\'c, Muha and Buka\v{c}~\cite{CanicMuhaBukac2012}, is a stable partitioned scheme
that uses operator splitting. Nobile and Vergara~\cite{NobileVergara2008}, among others, have developed 
stable semi-monolithic schemes for the case of thin shells. 
Even fully monolithic schemes for thin-shells~\cite{FigueroaVignonClementelJansenHughesTaylor2006} are
not overly expensive since the number of degrees of freedom for the structural shell is usually small compared
to those in the fluid.
For bulk solids, however, apparently all previous
partitioned approaches for coupling incompressible flows and light solids are
either unstable or require multiple sub-iterations per time-step to converge (tens or
hundreds of iterations are typical) with the number of iterations generally
increasing as the solid get lighter. Indeed, it is
common to resort to fully monolithic schemes~\cite{BatheZhang2004,Heil2004,HronTurekMonolithic2006,GeeKuttlerWall2011} or semi-monolithic
schemes~\cite{FernandezGerbeauGrandmont2007,AstorinoChoulyFernandez2009,GretarssonKwatraFedkiw2011}.

A number of authors have analyzed the added-mass effect and the stability of FSI
algorithms~\cite{Conca1997,LeTallecMani2000,MokWallRamm2001,ForsterWallRamm2007,AstorinoChoulyFernandez2009,vanBrummelen2009,Gerardo_GiordaNobileVergara2010}.
Causin, Gerbeau and Nobile~\cite{CausinGerbeauNobile2005}, for example, analyze a model
problem of a structural shell coupled to an incompressible fluid and show that
the traditional scheme can be unconditionally unstable over a certain range of
parameters.
There are numerous approaches that have been developed to reduce the number of iterations required per time-step
in partitioned schemes~\cite{MokWallRamm2001}.
Robin-Robin interface conditions are used to stabilize partitioned schemes by, for example,
Nobile, Vergara and co-workers~\cite{NobileVergara2008,BadiaNobileVergara2008,BadiaNobileVergara2009,Gerardo_GiordaNobileVergara2010,NobileVergara2012} and Astorino, Chouly and Fernandez~\cite{AstorinoChoulyFernandez2009}.
Yu, Baek and Karniadakis~\cite{YuBaekKarniadakis2013}, 
and Baek and Karniadakis~\cite{BaekKarniadakis2012} have developed {\em fictitious pressure} and
{\em fictitious mass} algorithms which incorporate additional terms into the governing equations to account
for added-mass effects. The optimal values for parameters are estimated from analysis and computations and the number of sub-iterations
required using Aitken acceleration (10's of iterations) is similar to optimal Robin-Robin approaches.
To stabilize the FSI problem and reduce the number of sub-iterations,
Riemslagh, Virendeels and Dick~\cite{RiemslaghVierendeelsDick2000} and 
Degroote et~al.~\cite{DegrooteSwillensBruggemanHaeltermanSegersVierendeels2010}, 
and Raback, Ruokolainen, and Lyly~\cite{RabackRuokolainenLyly2001}
have developed an {\em interface artificial compressibility method}
that adds a source term to the fluid continuity equation near the interface; 
the effect of the source term goes away as the sub-iterations converge.
Idelsohn et~al.~\cite{IdelsohnDelPinRossiOnate2012},  and Badia, Quaini and Quarteroni~\cite{BadiaQuainiQuarteroni2008}, 
form approximate factorizations of the fully
monolithic scheme to construct partitioned schemes but these still may require many iterations to converge. 
Degroote et~al.~\cite{DegrooteBruggemanHaeltermanVierendeels2008,DegrooteBatheVierendeels2009}
use reduced order models and Aitken acceleration methods to reduce the number of
iterations in partitioned schemes.

The remainder of the manuscript is organized as follows. In Section~\ref{sec:governingEquations} we define the governing equations and 
present the essential form of the AMP interface conditions. 
In Section~\ref{sec:partitionedAlgorithms} there is a discussion of partitioned algorithms and how the AMP
interface conditions are incorporated into the fluid and solid solution schemes. 
The three model problems we consider are defined in Section~\ref{sec:modelProblems}.
In Section~\ref{sec:AMP1D} the form of the AMP algorithm for a one-dimensional model problem is discussed.
A normal mode stability analysis of the AMP scheme as well as the traditional (Dirichlet-Neumann) and anti-traditional (Neumann-Dirichlet)
schemes are presented in Section~\ref{sec:analysis}.
The FSI time-stepping algorithm for the AMP scheme is outlined in Section~\ref{sec:timeSteppingAlgorithm}.
Section~\ref{sec:numericalResults} provides numerical results confirming the stability and accuracy
of the scheme. Conclusions are given in Section~\ref{sec:conclusions}.

%% file: texFiles/governingEquations.tex
\section{Governing equations and interface conditions} \label{sec:governingEquations}

We consider a fluid-structure interaction problem in which an incompressible
fluid in a domain $\OmegaF$ is coupled to a compressible elastic solid in
a domain $\OmegaS$.  The interface where the two domains meet is assumed to be
smooth and is denoted by $\Gamma$.  For the purposes of this paper, we consider
an FSI problem consisting of small perturbations about an equilibrium state.
Thus, we assume that the nonlinear convection terms in the fluid are negligible and
the solid is linearly elastic, and we consider the problem to be linearized
about a fixed interface position.  We also neglect the effects of gravity or other body forces.
Under these assumptions, the flow in the
(fixed) fluid domain is governed by the linear Stokes equations
\begin{alignat}{3}
  &  \rho\frac{\partial\vv}{\partial t} 
                 =  \grad\cdot\sigmav  , \quad&& \xv\in\OmegaF ,  \label{eq:stokes3d} \\
  & \grad\cdot\vv =0,  \quad&& \xv\in\OmegaF ,  \label{eq:fluidDiv3d}
\end{alignat}
where $\xv$ is position, $t$ is time, $\rho$ is the (constant) fluid density, and $\vv=\vv(\xv,t)$ is the fluid velocity.
The fluid stress tensor, $\sigmav=\sigmav(\xv,t)$, is given by 
\begin{equation}
  \sigmav = -p \Iv + \tauv,\qquad \tauv = \mu \left[ \grad\vv + (\grad\vv)^T \right],
  \label{eq:fluidStress}
\end{equation}
where $p$ is the pressure, $\Iv$ is the identity tensor, $\tauv$ is the viscous stress tensor, and $\mu$ is the (constant) fluid viscosity. 
For future reference, 
the components of a vector, such as $\vv$ will be denoted by $v_m$, $m=1,2,3$, (i.e. $\vv=[v_1, v_2, v_3]^T$), while components
of a tensor such as $\sigmav$, will be denoted by $\sigma_{mn}$, $m,n=1,2,3$.
In the (fixed) solid domain, the solid displacement $\usv=\usv(\xv,t)$ and velocity $\vsv=\vsv(\xv,t)$ are governed by 
\begin{alignat}{3}
  \frac{\partial\usv}{\partial t} &= \vsv,  \quad&& \xv\in\OmegaS , \label{eq:usFromvs} \\
  \rhos \frac{\partial\vsv}{\partial t} &= \grad \cdot\sigmasv  , \quad&& \xv\in\OmegaS , \label{eq:elasticWave}
\end{alignat}
where $\rhos$ is the (constant) solid density and $\sigmasv=\sigmasv(\xv,t)$ is the 
Cauchy stress tensor\footnote{
Over-bars on symbols are used throughout the paper to denote quantities belonging to the solid.}.
The stress tensor is given by
\begin{align*}
  \sigmasv &= \lambdas(\grad\cdot\usv)\Iv + \mus\left[ \grad\usv + \grad\usv^T \right],
\end{align*}
where $\lambdas$ and $\mus$ are (constant) Lam\'e parameters.  The matching conditions at the fluid-solid interface $\Gamma$ are
\begin{alignat}{3}
  \vv(\xv,t) &= \vvs(\xv,t), \quad && \xv\in\Gamma,                    \label{eq:interfaceV} \\
  \sigmav(\xv,t)\nv &= \sigmavs(\xv,t)\nv, \quad && \xv\in\Gamma, \label{eq:interfaceStress}
\end{alignat}
where $\nv$ is the outward unit normal to the fluid domain.  The problem is closed by specifying initial conditions for $\vv$ in $\OmegaF$ and for $\usv$ and $\vsv$ in $\OmegaS$, and by specifying suitable conditions on the boundaries of the fluid and solid domains not included in $\Gamma$.

The elastic wave equations in~\eqref{eq:usFromvs} and~\eqref{eq:elasticWave} form a hyperbolic system with a characteristic structure and finite wave speeds.  These wave speeds, given by
\[
\bar c_p=\sqrt{\frac{\lambdas+2\mus}{\rhos}},\qquad \bar c_s=\sqrt{\frac{\mus}{\rhos}},
\]
are the propagation speeds of p-waves and s-waves, respectively, in the elastic
solid.  At the interface the equations can be locally transformed into normal
and tangential coordinates. 
By considering the 
components of the equations in the normal direction, one can define characteristics
and characteristic variables~\cite{Whitham74} that indicate the propagation of information in the direction normal to the interface.
The {\em incoming} and {\em outgoing} characteristic variables
normal to the interface are given by
\begin{alignat}{3}
 \Achar(\sigmasv,\vsv)  &  = \nv^T\sigmasv\nv  - \zp \nv^T\vsv,  \qquad&& \xv\in\GammaI, \quad\text{(incoming)},  \label{eq:charA} \\
 \Achar_m(\sigmasv,\vsv)&  = \ev_m^T\sigmasv\nv - \zs \ev_m ^T\vsv, \quad m=1,2,  \qquad&& \xv\in\GammaI, \quad\text{(incoming)},\\
 \Bchar(\sigmasv,\vsv)  &  = \nv^T\sigmasv\nv  + \zp \nv^T\vsv,  \qquad&& \xv\in\GammaI,  \quad\text{(outgoing)}, \label{eq:charB} \\
 \Bchar_m(\sigmasv,\vsv)&  = \ev_m^T\sigmasv\nv   + \zs \ev_m ^T\vsv, \quad m=1,2, \qquad&& \xv\in\GammaI,\quad\text{(outgoing)}.       \label{eq:charBm}
\end{alignat}
Here, $\zp=\rhos\bar c_p$ and $\zs=\rhos\bar c_s$ are the solid impedances for
p-waves and s-waves, respectively, and $\ev_m$, $m=1,2$ denote mutually
orthogonal unit vectors tangent to the interface.  By {\em outgoing} we mean the
characteristics that leave the solid domain through the interface, and vice versa for the {\em incoming} characteristics.

Partitioned algorithms impose the interface conditions~\eqref{eq:interfaceV} and~\eqref{eq:interfaceStress}
by defining segregated conditions for the fluid and solid domains. 
In the traditional partitioned algorithm, the interface condition
in~\eqref{eq:interfaceV} is taken to be the (Dirichlet) boundary condition for the fluid
(velocity from solid), while~\eqref{eq:interfaceStress} is taken to be the (Neumann)
boundary condition for the solid (traction from fluid).  
Alternatively, the approach we develop in this article is
based on using the incoming and outgoing solid characteristic variables
normal to the interface given in~\eqref{eq:charA}--\eqref{eq:charBm}. 
In this approach a Robin interface condition for the fluid
is defined in terms of the outgoing solid characteristic variables, namely
\begin{equation}
\Bchar(\sigmav,\vv) =\Bchar(\sigmasv,\vsv), \qquad \Bchar_m(\sigmav,\vv) =\Bchar_m(\sigmasv,\vsv), \quad m=1,2,\qquad \xv\in\Gamma,
\label{eq:AMPoutgoing}
\end{equation}
while a Robin interface condition for the solid is defined in terms of the fluid variables using the incoming characteristic variables
\begin{equation}
\Achar(\sigmasv,\vsv) = \Achar(\sigmav,\vv), \qquad \Achar_m(\sigmasv,\vsv) = \Achar_m(\sigmav,\vv), \quad m=1,2,\qquad \xv\in\Gamma.
\label{eq:AMPincoming}
\end{equation}
This, in a nutshell, is the key ingredient of the added-mass partitioned (AMP)
algorithm.  These conditions are linearly independent combinations of the
interface conditions in~\eqref{eq:interfaceV} and~\eqref{eq:interfaceStress}, and
are thus equivalent to the original conditions.

In practice, numerical computations, including those presented in Section~\ref{sec:numericalResults}, are often based on alternative forms of the governing equations.
For the case of the fluid, the equations may be written in the velocity-pressure form
\begin{alignat}{3}
  &  \rho\frac{\partial\vv}{\partial t} + \nabla p
                 =  \mu\Delta\vv  , \quad&& \xv\in\OmegaF ,  \label{eq:NS3dv} \\
  & \Delta p = 0,  \quad&& \xv\in\OmegaF . \label{eq:pressurePoisson}
\end{alignat}
We use a numerical approximation of this form of the equations to advance the fluid variables following the fractional-step method described in~\cite{ICNS}.  
As discussed in~\cite{ICNS},
an additional boundary condition is required for this form of the equations, and a suitable choice is given by
\begin{align}
  & \grad\cdot\vv =0,  \qquad \xv\in\partial\OmegaF.
\label{eq:NS3dDiv}
\end{align}
For the solid domain, we consider the equations of linear elasticity written as the first-order system
\begin{alignat}{3}
  \frac{\partial\usv}{\partial t} &= \vsv,  \quad&& \xv\in\OmegaS , \label{eq:solidDisplacement} \\
  \rhos \frac{\partial\vsv}{\partial t} &= \grad \cdot\sigmasv  , && \xv\in\OmegaS ,\label{eq:solidMomemtum}\\
  \frac{\partial\sigmasv}{\partial t} &= \lambdas(\grad\cdot\vsv)\Iv + \mus\left[ \grad\vsv + \grad\vsv^T \right],\quad && \xv\in\OmegaS , \label{eq:solidStress}
\end{alignat}
and use a numerical approximation of this form of the equations to advance the solid variables following the second-order upwind (Godunov) method discussed in~\cite{smog2012}.

%% file: texFiles/partitionedAlgorithms.tex
\section{Partitioned algorithms for FSI} \label{sec:partitionedAlgorithms}

In this section we describe a traditional partitioned (TP) algorithm and compare
it to our new {\em added-mass} partitioned (AMP) algorithm.  For both
algorithms, we assume that the discrete solution of the FSI problem is known at
time $t^{n-1}$ with fluid state
$\qv^{n-1}=(\vv_\iv^{n-1},p_\iv^{n-1},\sigmav_\iv^{n-1})$ and solid state
$\qsv^{n-1}=(\uvs_\iv^{n-1},\vvs_\iv^{n-1},\sigmasv_\iv^{n-1})$ on a grid with
mesh points $\xv_\iv$.  The algorithms then describe how to obtain the solutions
$\qv^{n}$ and $\qsv^{n}$ at a time $t^n=t^{n-1}+\dt$, where $\dt$ is a chosen
time step.  For accuracy or 
stability reasons, multiple sub-iterations per time-step may be needed
in the TP algorithm~\cite{vanBrummelen2009}, 
and we let $\qv^{(k)}$
and $\qsv^{(k)}$ denote the $k\sp{{\rm th}}$ iterates which are approximations
of $\qv^{n}$ and $\qsv^{n}$, respectively.  
The AMP algorithm, in contrast, requires no sub-iterations, 
and remains stable even for problems where added-mass effects are significant.
Although not required, in practice we generally use a predictor-corrector algorithm, with a single
correction step, to advance the fluid since it has
a larger stability region that includes the imaginary axis~\cite{ICNS}. 
We let $\qv^{(p)}$ and $\qsv^{(p)}$ denote solutions
at the predictor step. 

We begin with a discussion of a TP algorithm that is typical of the traditional
algorithm used in practice.  In this traditional FSI algorithm, the interface
conditions in~\eqref{eq:interfaceV} and~\eqref{eq:interfaceStress} are assigned
separately in the two main steps, either the step to advance the solid or the
step to advance the fluid.  In particular, the velocity interface condition
in~\eqref{eq:interfaceV} is associated with the fluid domain so that, from the
point of view of the fluid, the interface is thought of as a no-slip moving
wall.  The traction condition in~\eqref{eq:interfaceStress}, on the other hand,
is associated with the solid domain so that, from the solid's point of view, the
interface is thought of as a no-slip traction wall.
\begin{algorithmTP}~~
\begin{enumerate}
  \item Choose some initial guess, e.g., $\qv^{(0)}=\qv^{n-1}$, $\qsv^{(0)}=\qsv^{n-1}$, and set $k=1$.
  \item Solve for $\qsv^{(k)}$ in the solid domain using the traction boundary condition 
   $(\sigmasv\nv)^{(k)}=(\sigmav\nv)^{(k-1)}$ on the interface.
  \item Solve for $\qv^{(k)}$ in the fluid domain using the velocity boundary condition $\vv^{(k)}=\vsv^{(k)}$ on the interface. 
  \item Set $k\leftarrow k+1$ and iterate steps 2--4, as needed.
\end{enumerate}
\label{alg:TP}
\end{algorithmTP}
\noindent The TP algorithm given above may be unstable or require many
under-relaxed sub-iterations per time-step to converge when the added-mass
effects are large (the {\em light} solid case). In fact, we show later that this
scheme can be unconditionally unstable no matter how heavy the solid is in comparison
to the fluid.

The AMP algorithm, which requires no sub-iterations, uses a different
construction of the interface conditions obtained from the characteristic-based
Robin (mixed) boundary conditions. 
The AMP algorithm outlined below contains the essential ingredients of the algorithm we use in practice.
It uses a single stage predictor; one can optionally include a correction step if desired.
More details are provided in subsequent discussions and in Section~\ref{sec:timeSteppingAlgorithm}.

\begin{algorithmAMP}~~

\begin{enumerate}
  \item Advance the solution in the solid domain to obtain $\qsv^{(p)}$, and compute predicted values of the outgoing characteristic variables, $\Bchar(\sigmasv^{(p)},\vsv^{(p)})$ and $\Bchar_m(\sigmasv^{(p)},\vsv^{(p)})$, $m=1,2$.
  \item Advance the solution in the fluid domain to obtain $\qv^{(p)}$ using the boundary conditions $\Bchar(\sigmav^{(p)},\vv^{(p)}) =\Bchar(\sigmasv^{(p)},\vsv^{(p)})$ and $\Bchar_m(\sigmav^{(p)},\vv^{(p)}) =\Bchar_m(\sigmasv^{(p)},\vsv^{(p)})$, $m=1,2$, obtained from~\eqref{eq:AMPoutgoing}.
  \item Define the interface traction~$(\sigmav\nv)^I$ to be the traction at the interface of the fluid domain, and compute an interface velocity $\vv^{I}$ using an impedance weighted average of the velocities at the interface of the fluid and solid domain (as is defined later).
  \item Apply solid interface conditions on $\qsv^{(p)}$ using $\vv^{I}$ and $(\sigmav\nv)^I$.
  \item Set $\qv^{n}=\qv^{(p)}$ and $\qsv^{n}$=$\qsv^{(p)}$. 
\end{enumerate}
\label{alg:AMP}
\end{algorithmAMP}

Having outlined the basic steps of the two FSI algorithms, we now provide some
further details of the new AMP algorithm.  In Step~1 of the AMP algorithm, the
discrete solution in the solid domain is advanced a time step $\dt$.  This can
be done, for example, using an approximation of the first-order system
in~\eqref{eq:solidDisplacement}--\eqref{eq:solidStress} as described
in~\cite{smog2012}.  Advancing the discrete solution in the fluid domain, as
given in Step~2, may be done, for example, using the velocity-divergence form of
the equations in~\eqref{eq:stokes3d} and~\eqref{eq:fluidDiv3d} along with
conditions on the interface given by
\begin{alignat}{3}
  & -p + \nv^T\tauv\nv  + \zp \nv^T\vv =\Bchar(\sigmasv^{(p)},\vsv^{(p)}), \qquad&& \xv\in\GammaI,  \label{eq:vdBCa} \\
  & \ev_m^T\tauv\nv   + \zs \ev_m ^T\vv =\Bchar_m(\sigmasv^{(p)},\vsv^{(p)}),\qquad m=1,2, \qquad&& \xv\in\GammaI, \label{eq:vdBCb}
\end{alignat}
(and some boundary conditions for $\xv\in\partial\OmegaF\backslash\GammaI$ which we do not discuss).  We have used~\eqref{eq:fluidStress} to eliminate $\sigmav$ in~\eqref{eq:vdBCa} and~\eqref{eq:vdBCb} to reveal the dependence on the pressure $p$ and the viscous stress $\tauv=\mu \left[ \grad\vv + (\grad\vv)^T \right]$.
This shows that~\eqref{eq:vdBCa} involves both $p$ and $\vv$ while~\eqref{eq:vdBCb} only depends on $\vv$. 

We prefer to advance the discrete solution in the fluid domain (for Step 2)
using a fractional-step method based on the velocity-pressure form of the
equations given in~\eqref{eq:NS3dv} and~\eqref{eq:pressurePoisson}.  In this
approach, the velocity and pressure are advanced in separate stages.  The two
conditions at the interface given by~\eqref{eq:vdBCb}, with $m=1,2$, along with
the divergence condition in~\eqref{eq:NS3dDiv} provide suitable conditions for
the numerical integration of~\eqref{eq:NS3dv} to advance the velocity.  For the
numerical solution of the Poisson problem for pressure
in~\eqref{eq:pressurePoisson} we use a boundary condition based on an alternate
version of the interface condition in~\eqref{eq:vdBCa}.  To motivate this
condition, we first consider the linear Taylor approximation
\[
\vv(\xv,t-\dt)\approx\vv(\xv,t)-\dt\frac{\partial\vv}{\partial t}(\xv,t)
\]
and a similar approximation for the solid velocity.  These Taylor approximations are used in~\eqref{eq:vdBCa} to obtain
\begin{alignat}{3}
 \nv^T\sigmav\nv + \zp\dt\,\nv^T\frac{\partial\vv}{\partial t}
      = \nv^T\sigmasv\nv + \zp\dt\,\nv^T\frac{\partial\vsv}{\partial t} ,\qquad\xv\in\GammaI,
 \label{eq:GeneralizedCharPtemp}
\end{alignat}
assuming that $\vv(\xv,t-\dt)=\vsv(\xv,t-\dt)$ for $\xv\in\GammaI$ according to~\eqref{eq:interfaceV}.  
Note that even though a first-order accurate Taylor approximation was used in deriving~\eqref{eq:GeneralizedCharPtemp}, 
the condition is actually identically true since it is a linear combination of~\eqref{eq:interfaceStress} 
and the time derivative of~\eqref{eq:interfaceV}.
Eliminating $\partial\vv/\partial t$ in~\eqref{eq:GeneralizedCharPtemp} 
using the fluid momentum equation~\eqref{eq:NS3dv} gives
\begin{alignat}{3}
  -p  - \frac{\zp\dt}{\rho}\,\frac{\partial p}{\partial n}   + \nv^T\tauv\nv 
  - \frac{\mu\zp\dt}{\rho}\,\nv^T(\grad\times\grad\times\vv) 
      = \nv^T\sigmasv\nv + \zp\dt\,\nv^T\frac{\partial\vsv}{\partial t},\qquad\xv\in\GammaI  .
 \label{eq:GeneralizedCharP}
\end{alignat}
The condition in~\eqref{eq:GeneralizedCharP} forms a suitable Robin condition
for pressure and is the key ingredient of the AMP algorithm for fractional step fluid solvers. 
 Note that the Robin condition only depends on $\vv$ through
the viscous traction and diffusion operator, which generally
makes~\eqref{eq:GeneralizedCharP} a better boundary condition to use when solving the
pressure equation separately.  We have also made the substitution $\Delta\vv=
-\grad\times\grad\times\vv$, valid when $\grad\cdot\vv=0$, since this improves
the stability of implicit time-stepping schemes~\cite{splitStep2003}.

The discussion of our fractional-step scheme to advance the solution in the fluid domain for Step 2 may be summarized 
by defining the following two fluid sub-problems:
\paragraph{Velocity sub-problem} 
Assuming a known solution in the fluid domain at time $t-\dt$ 
and a predicted solution in the solid domain at time $t$, 
a discrete solution to the velocity at time $t$ is determined by solving an appropriate discretization of
\begin{alignat}{3}
  &  \rho \frac{\partial \vv}{\partial t} 
               - \mu\Delta\vv = -\grad p,  \qquad&& \xv\in \OmegaF , \label{eq:fluidMomentumRiemann}  
\intertext{with boundary conditions}
  &\ev_m^T\vv +\frac{1}{\zs} \ev_m^T\tauv\nv\ = \ev_m^T\vsv +\frac{1}{\zs}\,\ev_m^T\sigmasv\nv,\qquad m=1,2, 
               \qquad&& \xv\in\GammaI,  \label{eq:solidVtangentialVsubproblem} \\
& \grad\cdot\vv =0,  \qquad&& \xv\in\GammaI.  \nonumber 
\end{alignat}
along with suitable boundary conditions for $\xv\in\partial\OmegaF\backslash\GammaI$.
\paragraph{Pressure sub-problem}
 Given the discrete solution for $\vv$ at time $t$ from the velocity sub-problem, 
the pressure is found using a discrete form of
\begin{alignat*}{3}
 &  \Delta p = 0 , \quad&& \xv\in\OmegaF, 
\end{alignat*}
with the Robin boundary condition 
\begin{alignat}{3}
 &     -p  - \frac{\zp\dt}{\rho} \frac{\partial p}{\partial n}
   = - \nv^T\tauv\nv + \frac{\mu\zp\dt}{\rho} \nv^T(\grad\times\grad\times\vv) 
     + \nv^T\sigmasv\nv +\zp\dt \, \nv^T\frac{\partial \vsv}{\partial t} 
            ,\qquad&& \xv\in\GammaI. \label{eq:AddedMassPressureBC}   
\end{alignat}
and suitable boundary conditions for $\xv\in\partial\OmegaF\backslash\GammaI$.

\bigskip
Moving on to Step 3 in the AMP algorithm, 
the interface traction, denoted by $(\sigmav\nv)^I$, is defined to be that from the fluid at time $t$,
\begin{align}
  (\sigmav\nv)^I = -p\nv + \tauv\nv,\qquad \xv\in\GammaI,
\end{align}
since the fluid velocity and pressure have already incorporated the primary AMP interface conditions.
The interface velocity, denoted by $\vv^I$, can be defined either from the
fluid velocity or can be computed from the characteristic relations
in~\eqref{eq:AMPoutgoing}, which may be written in the form
\begin{alignat}{3}
    \nv^T\vv &= \nv^T\vsv + \frac{1}{\zp}\nv^T( \sigmasv\nv-\sigmav\nv) ,\qquad&&\xv\in\GammaI,\label{eq:AMPa2} \\
    \ev_m^T\vv &= \ev_m^T\vsv\ + \frac{1}{\zs}\ev_m^T( \sigmasv\nv-\sigmav\nv)  , \qquad m=1,2, \qquad&&\xv\in\GammaI. \label{eq:AMPb2}
\end{alignat}
The choice given by~\eqref{eq:AMPa2} and~\eqref{eq:AMPb2} 
is the better
conditioned approximation for the case of {\em heavy} solids, while defining $\vv^I$
from the fluid velocity is better for {\em light} solids.  To smoothly accommodate
both limits, we define the interface velocity as an impedance-weighted
average of the two choices, namely
\begin{alignat}{3}
&  \nv^T\vvI = \frac{\zf}{\zf+\zp} \nv^T\vv + \frac{\zp}{\zf+\zp} \nv^T\vsv  
              + \frac{1}{\zf+\zp}\nv^T\big( \sigmasv\nv -(\sigmav\nv)^I \big),\qquad &&\xv\in\GammaI,  \label{eq:ImpedenceWeightedVn} \\
& \ev_m^T\vvI = \frac{\zf}{\zf+\zs} \ev_m^T\vvIf + \frac{\zs}{\zf+\zs} \tv_m^T\vsv 
               + \frac{1}{\zf+\zs}\tv_m^T\big( \sigmasv\nv -(\sigmav\nv)^I \big),         \quad m=1,2, \qquad &&\xv\in\GammaI.
\label{eq:ImpedenceWeightedVt}
\end{alignat}
Here, $\zf=\rho v_f$ is the fluid impedance, where $v_f$ is some suitable measure of the velocity in the 
fluid, as discussed further in Section~\ref{sec:analysis}.
Note that in practice, the algorithm is found to be very insensitive to the particular choice of $v_f$. 
This insensitivity is also confirmed in the theoretical results in Section~\ref{sec:analysis}.
Also note the form of the impedance weighted averaged in~\eqref{eq:ImpedenceWeightedVn} and~\eqref{eq:ImpedenceWeightedVt}
are the same form as those appearing in the added mass algorithm for compressible fluids~\cite{fsi2012}.
Finally, the interface values $\vv^I$ and $(\sigmav\nv)^I$ are used to assign boundary conditions on the solid in Step 4 of the AMP algorithm.

%% file: texFiles/modelProblems.tex
\section{FSI Model problems} \label{sec:modelProblems}

Three FSI model problems, of increasing complexity, are now defined. 
Model problem MP-IA, for an inviscid incompressible fluid and {\em acoustic} solid (defined below), 
is used in the two-dimensional analysis of partitioned schemes
in Section~\ref{sec:analysis}, as well as being the basis for the one-dimensional model problem discussed in Section~\ref{sec:AMP1D}.
The second model problem, MP-VA, includes the effects of viscosity in the fluid but retains the acoustic solid.
Model problem MP-VE includes viscosity in the fluid and treats a linearly elastic solid. 
Exact traveling wave solutions to these model problems are
given in~\ref{sec:travelingWave}, while numerical simulations are
given in Section~\ref{sec:numericalResults}.
In all cases the fluid domain is the rectangular region  $\OmegaF=(0,L)\times(-H,0)$, 
the solid domain is $\OmegaS=(0,L)\times(0,\Hs)$ and the interface is 
$\GammaI=\{ (x,y) \,\vert\, x\in(0,L),  y=0 \}$,
see Figure~\ref{fig:elasticShellCartoon}. 
We consider solutions that are periodic in the $x$-direction with period $L$.

\input texFiles/modelProblemGeometryFig

\begin{ModelProblemIA}\label{MP-IA}
Model problem MP-IA defines an inviscid incompressible fluid 
and an ``acoustic'' solid, that only supports vertical motion, 
\begin{equation}
\begin{aligned}
& \text{Fluid:}\;  \left\{ 
   \begin{alignedat}{3}
  &  \rho\frac{\partial\vv}{\partial t} 
                 + \grad p =0, \quad&& \xv\in\OmegaF ,  \\
  & \grad\cdot\vv =0,  \quad&& \xv\in\OmegaF ,   \\
  & v_2(x,-H,t)=0  ~~\text{or}~~ p(x,-H,t)=0 ,
   \end{alignedat}  \right. 
\qquad
 \text{Solid:} \; \left\{ 
   \begin{alignedat}{3}
  &  \rhos\frac{\partial^2\us_2}{\partial t^2}  = \rhos\cp^2\Delta\us_2, \quad&& \xv\in\OmegaS, \\
  &    u_2(x,\Hs,t)=0, 
   \end{alignedat}  \right.  \\
&\text{Interface:}\quad v_2=\frac{\partial\us_2}{\partial t},
       \quad -p=\rhos\cp^2\frac{\partial\us_2}{\partial y}, \qquad \xv\in\GammaI.
\end{aligned}
  \label{eq:MP-IA}
\end{equation}
\end{ModelProblemIA}

\vskip\baselineskip
\begin{ModelProblemVA}\label{MP-VA}
Model problem MP-VA defines a viscous incompressible fluid 
and an ``acoustic'' solid, that only supports vertical motion, 
\begin{equation*}
\begin{aligned}
& \text{Fluid:} \; \left\{ 
   \begin{alignedat}{3}
  &  \rho\frac{\partial\vv}{\partial t} 
                 + \grad p =\mu\Delta\vv, \quad&& \xv\in\OmegaF ,  \\
  & \grad\cdot\vv =0,  \quad&& \xv\in\OmegaF ,  \\
  &   \vv(x,-H,t)=0 ,  
   \end{alignedat}  \right. 
\qquad
\text{Solid:} \; \left\{ 
   \begin{alignedat}{3}
  &  \rhos\frac{\partial^2\us_2}{\partial t^2}  = \rhos\cp^2\Delta\us_2, \quad&& \xv\in\OmegaS, \\
  &    u_2(x,\Hs,t)=0, 
   \end{alignedat}  \right.  \\
&\text{Interface:}\quad v_1=0, \quad v_2=\frac{\partial\us_2}{\partial t}, 
       \quad -p + 2\mu \frac{\partial v_2}{\partial y} =\rhos\cp^2\frac{\partial\us_2}{\partial y}, \qquad \xv\in\GammaI.
\end{aligned}
\end{equation*}
\end{ModelProblemVA}


\vskip\baselineskip
\begin{ModelProblemVE}\label{MP-VE}
Model problem MP-VE defines a viscous incompressible fluid 
and an elastic solid, 
\begin{equation*}
\begin{aligned}
& \text{Fluid:} \; \left\{ 
   \begin{alignedat}{3}
  &  \rho\frac{\partial\vv}{\partial t} 
                 + \grad p =\mu\Delta\vv, \quad&& \xv\in\OmegaF ,  \\
  & \grad\cdot\vv =0,  \quad&& \xv\in\OmegaF ,  \\
  &   \vv(x,-H,t)=0 ,  
   \end{alignedat}  \right. 
\qquad
\text{Solid:} \; \left\{ 
   \begin{alignedat}{3}
  &  \rhos\frac{\partial^2\usv}{\partial t^2}  = (\lambdas+\mus)\grad(\grad\cdot\usv)+\mus\Delta\usv, \quad&& \xv\in\OmegaS, \\
  &    \usv(x,\Hs,t)=0, 
   \end{alignedat}  \right.  \\
&\text{Interface:}\quad \vv=\frac{\partial\usv}{\partial t}, \quad
  \mu(  \frac{\partial v_1}{\partial y} + \frac{\partial v_2}{\partial x}) =
               \mus(  \frac{\partial \us_1}{\partial y} + \frac{\partial \us_2}{\partial x}), 
  \quad -p + 2\mu \frac{\partial v_2}{\partial y} = 
        \lambdas\grad\cdot\usv +2\mus\frac{\partial\us_2}{\partial y} ,
\quad \xv\in\GammaI .
\end{aligned}
\end{equation*}
\end{ModelProblemVE}

%% file: texFiles/modelProblemGeometryFig.tex
%
%
{
\newcommand{\lbfont}{\small}
\newcommand{\ysa}{3.05}
\newcommand{\ysb}{5.}
\newcommand{\yb}{3}
\newcommand{\xL}{8}
\begin{figure}[hbt]
\newcommand{\textFont}{\normalss}
\begin{center}
\begin{tikzpicture}[scale=.75]
\useasboundingbox (0,0) rectangle (\xL,5.25);  
\draw[thick,red,fill=red,opacity=0.25] (0,\ysa) rectangle (\xL,\ysb);
\draw[thick,black] (4,4.) node {solid: $\OmegaS$};
\draw[thick,blue,fill=blue,opacity=0.25] (0,0) rectangle (\xL,\yb);
\draw[thick,blue] (0,0) rectangle (\xL,\yb);
\draw[thick,black] (4,1.25) node {fluid: $\OmegaF$};
\draw[thick,black] (4,\yb) node[anchor=north] {interface: $\GammaI$};

\draw[-,thick,red] (0  ,\ysa) -- (\xL,\ysa);
\draw[-,thick,red] (0  ,\ysb) -- (\xL,\ysb);
\draw[-,thick,red] (0  ,\ysa) -- (0  ,\ysb) node[anchor=east,black] {\lbfont$y=\Hs$};
\draw[-,thick,red] (\xL,\ysa) -- (\xL,\ysb);
\draw[-,thick,blue] (0,\yb) -- (\xL,\yb);
\draw[-,thick,blue] (\xL,0) -- (\xL,\yb);
\draw[-,thick,blue] (0,0) node[anchor=north,black] {\lbfont$x=0$} -- (\xL,0) node[anchor=north,black] {\lbfont$x=L$};
\draw[-,thick,blue] (0,0) node[anchor=east,black] {\lbfont$y=-H$} -- (0,\yb) node[anchor=east,black] {\lbfont$y=0$};
%
%
\end{tikzpicture}
\end{center}
\caption{The geometry for the 2D FSI model problems.}
\label{fig:elasticShellCartoon}
\end{figure}
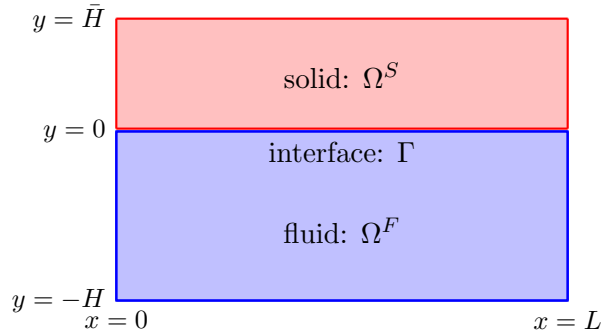
}

%% file: texFiles/ampAlgorithm1d.tex
\section{An illustration of the AMP algorithm}\label{sec:AMP1D}

As an illustration of the AMP algorithm, it is helpful to consider, as a
representative FSI problem, a one-dimensional version of model problem MP-IA given in~\eqref{eq:MP-IA}.
While this one-dimensional problem is
relatively simple, it does provide useful insight into the properties of the algorithm.  
To be consistent with the geometry of the two-dimensional problem, 
the one-dimensional problem varies in the $y$-direction, with $x$ fixed.
Consider then an inviscid fluid in the domain $\OmegaF=\{ y\in (-H,0)\}$ adjacent to a solid in
the domain $\OmegaS=\{ y\in (0,\infty)\}$ with interface at $y=0$ (and where $\Hs\rightarrow\infty$).  For this
one-dimensional problem, the velocity-pressure form of the equations
in~\eqref{eq:NS3dv} and~\eqref{eq:pressurePoisson} for the fluid 
reduces to 
\begin{equation}
\rho \frac{\partial v}{\partial t} + \frac{\partial p}{\partial y} = 0,\qquad \frac{\partial^2 p}{\partial y^2} = 0,\qquad y\in\OmegaF,
\label{eq:onedimensionalfluid}
\end{equation}
where $v=v_2$
and at $y=-H$. We take the traction boundary condition $p(-H,t)=0$. 
The equations in~\eqref{eq:solidMomemtum} and~\eqref{eq:solidStress} for the solid become
\begin{equation}
\rhos \frac{\partial \vs}{\partial t} = \frac{\partial \sigmas}{\partial y}, \qquad \frac{\partial\sigmas}{\partial t}  = \rhos\bar c_p \frac{\partial \vs}{\partial y}, \qquad y\in\OmegaS,
\label{eq:onedimensionalsolid}
\end{equation}
where $\vs=\vs_2$ and $\sigmas=\rhos\cp^2\partial_y\us_2$.
The system of equations for the solid can be written in the characteristic form
\[
\frac{d}{dt} ( \sigmas \pm \zp \vs )=0, \quad \text{along} \quad\frac{d y}{dt}=\mp\bar c_p,
\]
and the definitions for the incoming and outgoing characteristics at the interface in~\eqref{eq:charA} and~\eqref{eq:charB}, respectively, become
\[
\Achar(\sigmas,\vs) = \sigmas - \zp\vs,\qquad \Bchar(\sigmas,\vs) = \sigmas + \zp\vs,\qquad y=0.
\]
The standard interface conditions are
\begin{equation}
v(0,t)=\vs(0,t),\qquad \sigma(0,t)=\sigmas(0,t),
\label{eq:onedimesionalinterface}
\end{equation}
where $\sigma=-p$. The equivalent AMP interface conditions corresponding to~\eqref{eq:AMPoutgoing} and~\eqref{eq:AMPincoming} are
\begin{equation}
\sigma + \zp v=\Bchar(\sigmas,\vs),\qquad
\Achar(\sigmas,\vs) = \sigma - \zp v, \qquad y=0,  \label{eq:amp1d}
\end{equation}
where the first and second conditions in~\eqref{eq:amp1d} are thought of as the interface conditions for the fluid and solid,
  respectively.
The alternate AMP conditions corresponding~\eqref{eq:GeneralizedCharPtemp} and~\eqref{eq:AMPincoming}, 
suitable when solving the pressure equation, are
\begin{equation}
\sigma + \zp\dt\, \frac{\partial v}{\partial t} = \sigmas+ \zp\dt\,\frac{\partial\vs}{\partial t},\qquad
\Achar(\sigmas,\vs) = \sigma - \zp v, \qquad y=0.
\label{eq:onedimesionalinterface-alt}
\end{equation}
By combining the first two conditions in~\eqref{eq:amp1d} and~\eqref{eq:onedimesionalinterface-alt}, it can be shown that
\[
\frac{\partial}{\partial t}(v-\vs)=-\frac{1}{\dt}(v-\vs),\qquad y=0,
\]
and thus if $v(0,t)=\vs(0,t)$ at $t=0$, then the velocities on the interface are equal for all time.  Also, if $v(0,t)-\vs(0,t)$ is equal to some small nonzero value (from numerical error), then the difference in the velocity rapidly approaches zero.  Thus, the alternate version of the interface conditions in~\eqref{eq:onedimesionalinterface-alt} is essentially equivalent to those in~\eqref{eq:onedimesionalinterface}.

Given the solution at time $t-\dt$ the AMP algorithm proceeds first by taking a
time step of the solution in the solid domain to obtain predicted values,
$\sigmas^{(p)}$ and $\vs^{(p)}$, which then may be used to obtain outgoing
characteristic data on the interface.  Note that in one dimension, this data,
if obtained using the method of characteristics, is completely determined from
the solution in the solid domain at $t-\dt$, and is independent of the solution
in the fluid domain.

We now advance the solution in the fluid domain using the outgoing interface
data obtained from the solid domain, which is a key step in the AMP algorithm.
First consider the pressure sub-problem,
\[
\frac{\partial^2 \sigma}{\partial y^2} =0 , \qquad y\in(-H,0),
\]
with boundary conditions
\[
\sigma=0,\quad y=-H, \qquad \sigma + \frac{\zp\dt}{\rho}\, \frac{\partial\sigma}{\partial y} = \sigmas^{(p)}+ \zp\dt\,\frac{\partial\vs^{(p)}}{\partial t},\quad y=0.
\]
The boundary condition at $y=0$ is obtained from the outgoing condition in~\eqref{eq:onedimesionalinterface-alt} and the momentum equation for the fluid in~\eqref{eq:onedimensionalfluid}. 
The solution of the pressure sub-problem is
\begin{equation}
\sigma(y,t)=\sigma^I(t)\left(1+\frac{y}{H}\right)
\label{eq:pressureoned}
\end{equation}
where
\begin{equation}
\sigma^I(t) = \frac{\MR}{1+\MR}\left( \sigmas^{(p)}(0,t) + \zp\dt\frac{\partial\vs^{(p)}}{\partial t}(0,t) \right).
\label{eq:tractiononed}
\end{equation}
Observe that the traction at the interface, $\sigma^I$, is determined by the (known) outgoing characteristic data from the solid domain and involves the {\em added-mass} ratio $\MR$ given by
\[
\MR=\frac{\rho H}{\zp\dt}=\frac{\rho H}{\rhos\cp\dt} = \frac{\zf}{\zp},
\]
where $\zf=\rho H/\dt$ (or $\zf=\rho v_f$ with $v_f=H/\dt$) is a measure of the fluid impedance.  
The quantity $\MR$ may be interpreted as the ratio of the mass of the entire fluid domain, $\rho H$, to 
the mass of the solid displaced by its characteristic velocity over a time $\dt$, given by $\rhos\cp\dt$. 
It may also be interpreted as a ratio of the fluid to solid impedances. When $\MR$ is small, added-mass
effects are in some sense small. When $\MR$ is large, added-mass effects are large. It may seem odd that
$\MR$ becomes large (i.e. added-mass effects become large) as $\dy=\cp\dt$ becomes small (with $\rho H/\rhos$ fixed),
but this is consistent with the analysis of the traditional-partitioned scheme in Section~\ref{sec:TPanalysis}.  

In the context of the one-dimensional model problem, the velocity sub-problem is
\[
\rho\frac{\partial v}{\partial t}=\frac{\partial\sigma}{\partial y} , \qquad y\in(-H,0),
\]
with the boundary condition
\begin{equation}
\frac{\partial v}{\partial y}=0,\qquad y=0.
\label{eq:velocitybconed}
\end{equation}
Integrating the momentum equation and using the solution in~\eqref{eq:pressureoned} of the pressure sub-problem, we find
\begin{equation}
v(y,t)=v(y,t-\dt)+\frac{1}{\rho H}\int_{t-\dt}^{t}\sigma^I(\tau)\,d\tau.
\label{eq:velocityoned}
\end{equation}
The solution for velocity in~\eqref{eq:velocityoned} is spatially uniform so that the boundary condition in~\eqref{eq:velocitybconed} is satisfied identically.  Setting $y=0$ in~\eqref{eq:velocityoned} and using~\eqref{eq:tractiononed}, we obtain
\[
v(0,t)=v(0,t-\dt)+
    \frac{\MR}{\rho H(1+\MR)}
    \left[\int_{t-\dt}^{t}\sigmas^{(p)}(0,\tau)\,d\tau+\zp\dt\left(\vs^{(p)}(0,t)-\vs^{(p)}(0,t-\dt)\right)\right].
\]
Assuming the fluid and solid velocities on the interface are equal at $t-\dt$, and using a simple approximation of the integral gives
\begin{equation}
 v(0,t)=\frac{\MR}{1+\MR}v(0,t-\dt)+\frac{1}{1+\MR}\left[\vs^{(p)}(0,t)+\frac{1}{\zp}\sigmas^{(p)}(0,t)\right].
\label{eq:interfacevelocityoned}
\end{equation}
We observe that the velocity on the interface at time $t$ given
in~\eqref{eq:interfacevelocityoned} from the solution of the equations in the
fluid domain is an added-mass-weighted average of the velocity at $t-\dt$ and a
velocity determined by the outgoing characteristic data from the solution in the
solid domain.

The interface velocity $v^I$ may be taken as the velocity in~\eqref{eq:velocityoned} (or the approximation in~\eqref{eq:interfacevelocityoned}) determined by the solution in the fluid domain, or it may be taken as the velocity,
\[
\vs^{(p)}(0,t)+\frac{1}{\zp}\left(\sigmas^{(p)}(0,t)-\sigma^I(t)\right),
\]
determined by the outgoing characteristic condition at the interface (c.f.~\eqref{eq:AMPa2}).  
In fact, by eliminating $v(0,t-\dt)$ in~\eqref{eq:interfacevelocityoned} in terms of $\sigma^I$ using the approximation
$v(0,t)=v(0,t-\dt)+ \sigma^I\dt/(\rho H)$ to~\eqref{eq:velocityoned}, equation~\eqref{eq:interfacevelocityoned}
can be written as an impedance-weighted average of these two values 
\[
v^I(t)=v(0,t) = \frac{\zf}{\zf+ \zp} v(0,t) + \frac{\zp}{\zf+ \zp} \vs^{(p)}(0,t) + \frac{1}{\zf + \zp} \left(\sigmas^{(p)}(0,t) -\sigma^I(t)\right),
\]
which is the one-dimensional version of the formula in~\eqref{eq:ImpedenceWeightedVn}. 
In this one-dimensional problem, the definition of $v_f$ as $v_f=H/\dt$ naturally appears.
Given the velocity and traction at the interface, $v^I$ and $\sigma^I$, respectively, the solution for the solid can now be fully determined by solving the governing equations in~\eqref{eq:onedimensionalsolid} together with the boundary data
\begin{align*}
    \sigmas(0,t) - \zp \vs(0,t) &= \sigma^I(t) - \zp v^I(t) ,  
\end{align*}
from the incoming characteristic.

The above description follows the AMP algorithm for a simple FSI problem in one
dimension. The key step in the algorithm, as mentioned previously, is the
incorporation of outgoing characteristic data from the solid domain into the
fluid domain.  The application of the AMP algorithm to this simple FSI problem
also illustrates the contribution of the added-mass ratio $\MR$.  In addition,
it can be confirmed that for $\rho H \ll \zp\dt$ ($\MR\ll
1$) the algorithm approaches the standard TP algorithm (defined previously in
Algorithm~\ref{alg:TP}) with the interface velocity being primarily transmitted
from the solid, while the traction comes primarily from the fluid.  On the other
hand, for $\rho H \gg \zp\dt$ ($\MR\gg 1$) the AMP algorithm approaches an {\em
anti-traditional algorithm} with the roles of fluid and solid reversed in the
application of the boundary conditions at the interface for each domain. 
Note that 
the AMP scheme always tends to the anti-traditional scheme for $\dt$ sufficiently small; this is
consistent with the analysis in Sections~\ref{sec:TPanalysis} and~\ref{sec:ATPanalysis} which indicates 
that the anti-traditional scheme is stable when $\dy$ and $\dt$ become sufficiently small, while the TP algorithm is not.

%


%% file: texFiles/modelProblemAnalysis.tex

\section{Analysis of a two-dimensional acoustic solid and an inviscid incompressible fluid} \label{sec:analysis}

In this section, we perform a stability analysis of the AMP algorithm applied to
the two-dimensional FSI model problem MP-IA.  In this problem, the fluid is
taken as inviscid and the solid is treated as an acoustic solid that only
supports motion in the vertical direction. The governing equations for the model
are given by~\eqref{eq:MP-IA}.
We consider a semi-infinite solid domain ($\Hs=\infty$) and look for solutions in the solid that have finite $L_2$-norm (which implies the solutions decay to zero as $y\rightarrow\infty$).  The bottom boundary condition on the fluid, at $y=-H$, is chosen as $\sigma=0$ (i.e. $p=0$).

\subsection{Discretization}

We discretize the fluid and solid variables in the $x$-direction using a uniform
grid with spacing $\dx=L/(N_x+1)$. The grid points are given by $x_\ell=\ell\dx$,
$\ell=0,1,\ldots,N_x$.  The choice of discrete approximations to the
$x$-derivatives in the equations plays little role in the analysis.  Therefore, it is
convenient to use a pseudo-spectral approximation by expanding
each component $q$ of the solution $(\qv,\qsv)$ in a discrete Fourier series of
the form
\begin{align}
   q_\ell(y,t) = \sum_{k=-N_x/2}^{N_x/2} e^{2\pi i k x_\ell/L} \qHat_{k}(y,t), \qquad \ell=0,1,2,\ldots,N_x, \label{eq:discreteFourier}
\end{align}
where $q_\ell(y,t)\approx q(x_\ell,y,t)$ and $\qHat_{k}(y,t)$ are Fourier coefficients, 
and where $N_x$ is assumed to be even for convenience.
Transforming the governing equations to Fourier space leads to the following equations for the Fourier coefficients of the solution variables:
\begin{alignat}{3}
&   \rhos \partial_t \vs = i\kx \sigmas_{21} + \partial_y \sigmas_{22} , \quad&& y\in(0,\infty), \label{eq:vEqn} \\
&   \partial_t \sigmas_{22} = \rhos \cp^2 \partial_y\vs , \quad&&y\in(0,\infty), \label{eq:s22Eqn} \\
&   \partial_t \sigmas_{21} = i\kx \rhos \cp^2 \vs , \quad&&y\in(0,\infty) , \label{eq:s21Eqn} \\
& \rho \partial_t v_1 = i\kx      \sigma, \quad &&y\in(-H,0), \label{eq:v1HatEqn} \\
& \rho \partial_t v_2 = \partial_y\sigma, \quad && y\in(-H,0), \label{eq:v2HatEqn} \\
&  i\kx v_1 + \partial_y v_2 =0 ,                  \quad&&y\in(-H,0) , \label{eq:divHatEqn}
\end{alignat}
where $\kx = 2\pi k/L$ is a normalized wave number in the $x$-direction and the
equations hold for each $k=-N_x/2,-N_x/2+1,\ldots,N_x/2$. The hats on the
Fourier coefficients, along with the $k$ subscripts, have been dropped for
notational convenience.  Note that the horizontal component of the fluid
velocity, $v_1$, is decoupled from the other equations and boundary conditions, 
and can be determined once $v_2$ is known. Equations~\eqref{eq:v1HatEqn}--\eqref{eq:divHatEqn} are used to obtain
\[
\partial_y^2 \sigma - \kx^2 \sigma =0,      \qquad y\in(-H,0) ,
\]
which is a Laplace equation for pressure (in Fourier space).  The equations for the Fourier coefficients of the solid variables will be discretized using an upwind scheme.  In preparation for this, it is convenient to define the characteristic variables
\[
a=\sigmas_{22}-\rhos\cp\vs,\qquad b=\sigmas_{22}+\rhos\cp\vs,\qquad d=\sigmas_{21},
\]
which, from~\eqref{eq:vEqn}--\eqref{eq:s21Eqn}, satisfy
\begin{alignat}{3}
&   \partial_t a + \cp \partial_y a =  i\kx\cp\, d , \qquad&& y\in(0,\infty),  \label{eq:aHatEqnII}\\
&   \partial_t b - \cp \partial_y b = -i\kx\cp\, d , \qquad&& y\in(0,\infty),  \\
&   \partial_t d = \frac{i\kx\cp}{2}(b-a) , \qquad&& y\in(0,\infty).  \label{eq:dHatEqnII} 
\end{alignat}

The solid variables are discretized in the $y$-direction using a uniform grid with
spacing $\dy$.  The grid points for the solid domain are defined as $y_j=
(j-\half)\dy$, $j=0,1,2,\ldots$. 
Following the analysis in~\cite{sjogreenBanks2012,lrb2013}, the grid is staggered with respect to the
interface at $y=0$ and has a ghost point at $j=0$; use of a non-staggered grid leads to similar results.
Let  $\dt$ be the time step and let $t^n=n\dt$.  The
equations~\eqref{eq:aHatEqnII}--\eqref{eq:dHatEqnII} are discretized using the
first-order upwind scheme,
\begin{alignat}{3}
&   \frac{ a_j^{n+1}- a_j^n}{\dt} + \cp\frac{( a_j^n - a_{j-1}^n)}{\dy} =   i\kx\cp\, d_j^{n}, \label{eq:aEqnIII} \\
&   \frac{ b_j^{n+1}- b_j^n}{\dt} - \cp\frac{( b_{j+1}^n - b_j^n)}{\dy} = - i\kx\cp\, d_j^{n}, \label{eq:bEqnIII} \\
&   \frac{ d_j^{n+1}- d_j^n}{\dt} = \frac{i\kx\cp}{2}(b_j^{n+1} - a_j^{n+1}), \label{eq:dEqnIII}
\end{alignat}
where $a_j^{n}\approx a(y_j,t^n)$, $b_j^{n}\approx b(y_j,t^n)$ and
$d_j^{n}\approx d(y_j,t^n)$. 
Since we have used a non-dissipative pseudo-spectral approximation, rather than a two-dimensional upwind scheme, 
the right-hand side terms in~\eqref{eq:dEqnIII} are taken at time $t^{n+1}$ in order to stabilize the approximation.
A von Neumann stability analysis of equations~\eqref{eq:aEqnIII} and~\eqref{eq:bEqnIII} (for a
periodic problem in $y$) shows that this time-stepping scheme is stable under
reasonable conditions on the constants $\ly=\cp\dt/\dy$ and $\lx=\cp\kx\dt$.
The fluid variables are kept continuous in $y$ in order to simplify the
presentation; a discrete version can be introduced but this makes only
insignificant changes to the fundamental results.  We are thus led to the
following semi-discrete approximation for the solid and fluid equations,
\begin{alignat}{3}
&   a_j^{n+1} = a_j^n -\ly( a_j^n - a_{j-1}^n) - i \lx d_j^{n} , \quad&& j=1,2,3,\ldots\;, \label{eq:aEqn} \\
&   b_j^{n+1} = b_j^n +\ly( b_{j+1}^n - b_j^n) + i \lx d_j^{n},  \quad&& j=0,1,2,\ldots\;, \label{eq:bEqn}  \\
&   d_j^{n+1} = d_j^n +  \frac{i\lx}{2} ( b_j^{n+1} - a_j^{n+1}),    \quad&& j=0,1,2,\ldots\;,  \\
&   v^{n+1} = v^n + \frac{\dt}{\rho} \partial_y \sigma^{n+1},             \quad&& y\in(-H,0), \\
&  0=\partial_y^2 \sigma^{n+1}- \kk^2 \sigma^{n+1},                             \quad&& y\in(-H,0), \label{eq:sigmaEqn}
\end{alignat}
where $v^n(y)=v(y,t^n)$ and $\sigma^n(y)=\sigma(y,t^n)$.  For reference, the discrete values for the (Fourier coefficients of the) velocity and components of stress in the solid are related to the characteristic variables by
\[
\vs_j^n=\frac{1}{2\zp}(b_j^n-a_j^n),\qquad\sigmas_{22,j}^n=\frac{1}{2}(b_j^n+a_j^n),\qquad \sigmas_{21,j}^n=d_j^n.
\]
The conditions at the top and bottom boundaries are 
\begin{align}
&   \vert a_j^{n}\vert^2 + \vert b_j^{n}\vert^2 + \vert d_j^{n}\vert^2 \rightarrow 0, 
         \quad\text{as $j\rightarrow \infty$} ,  \label{eq:BCtop}\\
&   \sigma^n(-H)=0. \label{eq:BCbottom}
\end{align}
Initial conditions are required to define $a_j^{0}$, $b_j^{0}$ and $d_j^{0}$ for $j=0,1,2,3,\ldots$, as well as $v^0(y)$ for  $y\in(-H,0)$, but these conditions are of no significance in the subsequent stability analysis.

Approximations to the interface conditions are needed to complete the discrete
formulation of the FSI model problem, and various choices are possible depending
on the algorithm used.  
For the AMP algorithm, the pressure Robin
condition~\eqref{eq:AddedMassPressureBC} is imposed on the fluid; this is derived from the
outgoing characteristic variables~\eqref{eq:AMPoutgoing}. In addition, the 
incoming characteristic variable~\eqref{eq:AMPincoming} is used to specify a condition
on the solid.  The latter condition uses an interface velocity defined
in~\eqref{eq:ImpedenceWeightedVn} based on an impedance-weighted average.  
First-order accurate approximations of these conditions are
\begin{alignat}{3}
&  \sigma^{n+1}(0) + \frac{\zp\dt}{\rho} \sigma_y^{n+1}(0) = \sigmas_{22,1}^{n+1} + \zp(\vs_1^{n+1} - v^{I,n})
                                           = b_1^{n+1} -\zp v^{I,n}, && \label{eq:1dSigmaBC}\\
&  a_0^{n+1} = \sigma^{I,n+1} - \zp v^{I,n+1}, &&\label{eq:1dSolidBC}
\end{alignat}
where $\sigma^{I,n}\equiv\sigma^n(0)$ and $v^{I,n}$ are the interface traction and velocity, respectively.  The latter is given by the impedance-weighted average
\begin{alignat}{3}
&   v^{I,n+1} = \frac{\zf}{\zf+\zp} v^{n+1}(0) + \frac{\zp}{\zf+\zp}\vs_1^{n+1} + \frac{1}{\zf+\zp}(\sigma_{22,1}^{n+1} - \sigma^{I,n+1} ) && \nonumber\\
& \hphantom{v^{I,n+1}} = \gamma v^{n+1}(0) + \frac{1-\gamma}{\zp}( b_1^{n+1} - \sigma^{I,n+1} ), && \label{eq:1dVelocityBC}\end{alignat}
where
\[
\gamma=\frac{\zf}{\zf+\zp},
\]
and where $\zf$ is a measure of the fluid impedance as introduced in Section~\ref{sec:partitionedAlgorithms}.

The traditional partitioned (TP) algorithm uses the velocity from the solid as a boundary condition for the fluid and the traction from the fluid as a boundary condition for the solid.  First-order accurate approximations of these conditions are
\begin{alignat}{3}
 &   v^{n+1}(0)=\vs_1^{n+1}  
        , &&   \label{eq:1dTradv} \\
 &   \sigmas^{n+1}_{22,0} = \sigma^{n+1}(0). &&\label{eq:1dTradsigma}
\end{alignat}
It is also instructive to consider an anti-traditional algorithm in which
the roles of the fluid and solid are reversed in the application of the
partitioned interface conditions.  Approximations of these conditions are
\begin{alignat}{3}
&   \sigma^{n+1}(0) = \sigmas^{n+1}_{22,1}\,,  &&  \label{eq:1dAntiTradsigma} \\
&   a_0^{n+1} = \sigma^{n+1}(0) - \zp v^{n+1}(0) .&& \label{eq:1dAntiTradv} 
\end{alignat}
The condition in~\eqref{eq:1dAntiTradv} is equivalent to setting the velocity of solid at the interface equal to the velocity in the fluid since the stresses have already been equilibrated in~\eqref{eq:1dAntiTradsigma}. 

\subsection{Stability analysis}

\renewcommand{\aHat}{{\tilde a}}
\renewcommand{\bHat}{{\tilde b}}
\renewcommand{\dHat}{{\tilde d}}
\renewcommand{\vHat}{{\tilde v}}
\renewcommand{\sigmaHat}{{\tilde\sigma}}
\renewcommand{\aTilde}{\alpha}
\renewcommand{\bTilde}{\beta}

To analyze the stability of the approximation
in~\eqref{eq:aEqn}--\eqref{eq:sigmaEqn}, we look for solutions of the form \[
a_j^n=\amp^n\aHat_j, \qquad b_j^n=\amp^n\bHat_j, \qquad d_j^n=\amp^n\dHat_j,
\qquad v^n(y)=\amp^n\vHat(y), \qquad \sigma^n(y)=\amp^n\sigmaHat(y),
\] 
where $\amp\in\Complex$ is the amplification factor,
$(\aHat_j,\bHat_j,\dHat_j)$ are grid functions, and $(\vHat(y),\sigmaHat(y))$
are functions for $y\in[-H,0]$. 
 Substituting these forms
into~\eqref{eq:aEqn}--\eqref{eq:sigmaEqn} results in the following system of
equations
\begin{alignat}{3}
&  \amp \aHat_j = \aHat_j -\ly( \aHat_j - \aHat_{j-1}) - i \lx \dHat_j , \qquad&& j=1,2,3,\ldots\,, \label{eq:aHatEqn} \\
&  \amp  \bHat_j = \bHat_j +\ly( \bHat_{j+1} - \bHat_j) + i \lx \dHat_j,  \qquad&& j=0,1,2,\ldots\,, \label{eq:bHatEqn}  \\
&  \amp  \dHat_j = \dHat_j + \frac{i\lx}{2}\amp  ( \bHat_j - \aHat_j),    \qquad&& j=0,1,2,\ldots\,,  \label{eq:dHatEqn} \\
&  \amp  \vHat(y) = \vHat(y) + \amp\dt\,\sigmaHat^I\,\frac{\kk\cosh(\kk (y+H))}{\rho\sinh(\kk H)},  \qquad&& y\in[-H,0], \label{eq:vHatEqn} \\
&  \sigmaHat(y)= \sigmaHat^I\, \frac{\sinh(\kk (y+H))}{\sinh(\kk H)},          \qquad&& y\in[-H,0].\label{eq:sigmaHatEqn}
\end{alignat}
To obtain~\eqref{eq:vHatEqn} and~\eqref{eq:sigmaHatEqn},  the
ODE in~\eqref{eq:sigmaEqn} was integrated using the boundary conditions $\sigmaHat(-H)=0$
and $\sigmaHat(0)=\sigmaHat^I$, where the fluid stress at the interface, $\sigmaHat^I$, is as
yet unspecified.  Note that the solution in the one-dimensional case,
$\kk=0$, is a special case that can be found by taking the limit $\kk\rightarrow
0$ in~\eqref{eq:vHatEqn} and~\eqref{eq:sigmaHatEqn}. 
Together with the interface
conditions, equations~\eqref{eq:aHatEqn}--\eqref{eq:sigmaHatEqn} define
a homogeneous set of difference equations. These equations have
non-trivial solutions (i.e. eigenfunctions) only for particular values of $\amp$ (i.e. eigenvalues).
The scheme is said to be weakly stable, or stable in the sense of
Godunov-Ryabenkii, for given values of the parameters $\dt$, $\dy$, $\kx$, $H$, etc., 
if there are no solutions (i.e.~eigenfunctions) to this
eigenvalue problem with $\vert \amp\vert>1$. 
The {\em region of stability} of the scheme is the region in parameter space where there are no roots with $\vert \amp\vert>1$. 
In order to delineate the stability region 
it is generally easier to look for regions where the scheme is {\em not} stable. Therefore, we 
assume that $\vert \amp\vert>1$ and look for solutions to the eigenvalue problem.  

Equations~\eqref{eq:aHatEqn}--\eqref{eq:dHatEqn} are a system of
constant-coefficient difference equations for the solid characteristic
variables.  After solving for $\dHat_j$ from~\eqref{eq:dHatEqn} and substituting
into~\eqref{eq:aHatEqn} and~\eqref{eq:bHatEqn}, $\aHat_j$ and
$\bHat_j$ are found to satisfy the recursion
\begin{equation}
 \begin{bmatrix} \aHat_j \\ \bHat_j \end{bmatrix} = K \begin{bmatrix} \aHat_{j-1} \\ \bHat_{j-1} \end{bmatrix}, \qquad K =\begin{bmatrix}  ( 1 - \theta^2)/R & -\theta \\ \theta & R \end{bmatrix} ,
\label{eq:recursion}
\end{equation}
where
\begin{equation}
R = r-\theta, \qquad r = \frac{\amp-1+\ly}{\ly}, \qquad \theta = -\frac{\amp\lx^2}{2\ly(\amp-1)} . 
\label{eq:RandTheta}
\end{equation}
The eigenvalues, $\phi_-$ and $\phi_+$, of
the $2\times2$ matrix $K$ in~\eqref{eq:recursion}
are given by
\begin{equation}
\phim = B - \sqrt{ B^2 - 1}, \qquad \phip=1/\phim, \qquad B = \frac{1}{2}\left(R + \frac{1 - \thetaj^2}{R}\right).
\label{eq:phipm}
\end{equation}
Note that the product of the eigenvalues is one.  
If both eigenvalues 
had magnitude equal to one, then there would be no solution of the recursion
in~\eqref{eq:recursion} satisfying the boundary condition~\eqref{eq:BCtop}.
Therefore, for a valid solution, there is one eigenvalue (either $\phi_-$ or $\phi_+$),
denoted by $\phis$ (for {\em small} $\phi$), that has magnitude strictly less
than one, i.e.~$\vert\phis\vert<1$.  Let $\phib=1/\phis$ (for {\em big} $\phi$)
be the other eigenvalue with $\vert\phib\vert>1$.  Note that when $\theta=0$,
$0<\ly\le 1$, and $\vert\amp\vert>1$, then $\phis=\phim=1/r$ (but in general
$\phis$ is not always equal to $\phim$).  Given that the eigenvalues are
distinct, $K$ can be diagonalized as 
\[
K = S^{-1} \Phi S , \qquad
\Phi =\begin{bmatrix}  \phis &  0 \\
                           0 & \phib \end{bmatrix}, \qquad
S =\begin{bmatrix}  R-\phis & \thetaj \\
                    R-\phib & \theta \end{bmatrix}.
\]
By setting 
\begin{align}
&   \begin{bmatrix} \aTilde_j \\ \bTilde_j \end{bmatrix}
   =  S \begin{bmatrix} \aHat_j \\ \bHat_j \end{bmatrix} = 
   \begin{bmatrix} (R-\phis)\aHat_j + \theta \bHat_j \\
                   (R-\phib)\aHat_j + \theta \bHat_j
   \end{bmatrix}, \label{eq:abTildeEqn}
\end{align}
the general solution of the recursion in~\eqref{eq:recursion} is given by 
\begin{equation}
\begin{array}{l}
   \aTilde_j = \phis^j \aTilde_0, \smallskip\\
   \bTilde_j = \phib^j \bTilde_0.
\end{array}
\label{eq:abBarEqn}
\end{equation}
Since the solution must decay to zero as $j\rightarrow \infty$ and since $\vert\phib\vert>1$, it follows that $\bTilde_0=0$ and thus $\bTilde_j=0$ for all $j=0,1,2,\ldots$.  Since $\bTilde_j=0$, it follows from~\eqref{eq:abTildeEqn} that
\begin{equation}
\bHat_j = Q \aHat_j,\qquad \aTilde_j=(\phib-\phis)\aHat_j,\qquad Q = \frac{\phib-R}{\theta},
\label{eq:Qdef}
\end{equation}
and thus from~\eqref{eq:abBarEqn}, 
\[
\aHat_j=\phi_s^j\aHat_0\,.
\]
Whence,
\begin{proposition}
Solutions to the system in~\eqref{eq:aHatEqn}--\eqref{eq:sigmaHatEqn} with $\vert\amp\vert>1$ satisfying the boundary conditions in~\eqref{eq:BCtop} and~\eqref{eq:BCbottom} are given by
\[
\begin{array}{c}
\displaystyle{
\aHat_j=\phi_s^j\aHat_0,\qquad \bHat_j=Q\aHat_j,\qquad \dHat_j=\frac{i\lx(Q-1)\amp}{2(\amp-1)}\aHat_j,\qquad j=0,1,2,\ldots\,,
} \bigskip\\
\displaystyle{
\sigmaHat(y)= \sigmaHat^I\, \frac{\sinh(\kk (y+H))}{\sinh(\kk H)},\qquad \vHat(y)=\sigmaHat^I\,\frac{\amp\dt}{\amp-1}\;\frac{\kk\cosh(\kk (y+H))}{\rho\sinh(\kk H)},\qquad y\in[-H,0],
}
\end{array}
\]
where $\aHat_0$ and $\sigmaHat^I$ are free constants.
\end{proposition}

The remaining constraints needed to complete the eigenvalue problem, which then determines the amplification factor $\amp$, are provided by two interface conditions.  Choices for these interface conditions are discussed in the subsections below.


\subsection{AMP interface conditions} 
Imposing the AMP interface conditions in~\eqref{eq:1dSigmaBC} and~\eqref{eq:1dSolidBC} with the formula for the interface velocity in~\eqref{eq:1dVelocityBC} gives
\begin{align}
&  \amp\sigmaHat^{I} = \eta\left( \amp\bHat_1 - \zp \vHat^{I}\right), \label{eq:AmpSigmaHat} \\
&  \aHat_0 = \sigmaHat^{I} - \zp \vHat^{I} , \label{eq:AmpaHat} \\
&  \vHat^{I}  = \gamma \vHat(0) + \frac{1-\gamma}{\zp}( \bHat_1 - \sigmaHat^{I} ), \label{eq:AmpvHatI}
\end{align}
where
\begin{align*}
 & \eta \equiv \frac{1}{1 + \frac{\zp\kk\dt}{\rho}\coth(\kk H)}  = \frac{1}{1 + \frac{\rhos\lx}{\rho}\coth(\kk H)}. 
\end{align*}
From~\eqref{eq:vHatEqn} and~\eqref{eq:AmpSigmaHat} we also have
\begin{align}
  (\amp-1) \vHat(0) = \frac{1-\eta}{\zb}(\amp \bHat_1 - \zb\vHat^I). \label{eq:AmpvHat} 
\end{align}
Using $\bHat_1=Q\phis \aHat_0$ to eliminate $\bHat_1$ in terms of $\aHat_0$ and using~\eqref{eq:AmpSigmaHat} to
eliminate $\sigmaHat^{I}$, equations~\eqref{eq:AmpaHat}, \eqref{eq:AmpvHat} and~\eqref{eq:AmpvHatI} can be written
in the form 
%
\begin{align}
\begin{bmatrix}
   \amp(1-\eta Q\phis)            & 0      & \amp+\eta \\
   (1-\eta)\amp Q \phis           & 1-\amp    & \eta-1 \\
  -(1-\gamma)(1-\eta)Q\phis  & -\gamma &  1-(1-\gamma)\eta \amp^{-1} 
\end{bmatrix} 
\begin{bmatrix} \aHat_0 \\ \zb \vHat(0) \\ \zb\vHat^I \end{bmatrix} =0 .  \label{eq:AMPmatrix}
\end{align}
The determinant of the matrix in~\eqref{eq:AMPmatrix} must be zero for non-trivial solutions to exist, and this leads to the following equation for $\amp$:
\begin{align*}
   \big(\amp-(1-\gamma)\big) \left(\amp-\eta - \big( (2\eta-1)\amp -\eta \big)Q \phis\right) =0.
\end{align*}
Since $\amp-(1-\gamma)\ne0$ when $\vert\amp\vert>1$ and $0\le \gamma\le 2$, we have the following result.
\begin{theorem}
The AMP interface approximation to the
scheme~\eqref{eq:aEqn}--\eqref{eq:sigmaEqn} with boundary
conditions~\eqref{eq:BCtop} and~\eqref{eq:BCbottom}, interface
conditions~\eqref{eq:1dSigmaBC} and~\eqref{eq:1dSolidBC}, and interface
velocity~\eqref{eq:1dVelocityBC} is weakly stable if and only if $0\le \gamma\le 2$ and there are no roots $\amp$ with $\vert \amp\vert>1$ to the equation
\begin{align}
   f(\amp) \equiv& \amp -\eta - \Big( (2\eta-1)\amp -\eta \Big)Q \phis =0. \label{eq:fAMP}
\end{align}
\end{theorem}
\noindent Note that $f(\amp)$ does not depend on $\gamma$, i.e.~the weighting in the definition of the interface velocity, and thus the stability of the scheme is independent of $\gamma$ provided $0\le \gamma\le 2$.  From a numerical conditioning point of view, however, we recommend the choice $\gamma=\zf/(\zf+\zp)$ as used in~\eqref{eq:1dVelocityBC}.

\paragraph{Reduction to a polynomial equation for $\amp$} To examine the stability of the AMP algorithm we must consider all solutions of $f(\amp)=0$.  The function $f(\amp)$ in~\eqref{eq:fAMP} involves square roots of $\amp$ through the eigenvalues, $\phis$ and $\phib$, of the matrix $K$ in~\eqref{eq:recursion} and the definition of $Q$ in~\eqref{eq:Qdef}.  
To ensure that all possible solutions are found, 
a polynomial $P(\amp)$ can be derived whose roots include
all the roots of $f(\amp)$.  To this end, 
note that from the definition for $Q$ and $\phis\phib=1$ it follows that 
\begin{align*}
  Q\phis = \frac{1-R\phis}{\theta}.
\end{align*}
Thus $f$ is a linear function of $\phis$ since $f$ only depends on the product $Q\phis$. 
Note 
also that $Q\phis$ is a rational function of $\amp$ according the definitions of $R$ and $\theta$ in~\eqref{eq:RandTheta}.  
Thus, $f(\amp)=0$ can be written in the form
\begin{align*}
    \phis = F(\amp) ,
\end{align*}
where $F$ is a rational function of $\amp$. 
Since $\phis$ is equal to one of the two eigenvalues, $\phi_\pm$, defined in~\eqref{eq:phipm}, it follows that 
\begin{equation}
  \pm \sqrt{ B^2 - 1} = F(\amp)-B,
\label{eq:Fratio}
\end{equation}
where the sign in front of the square root is chosen to be consistent with $\vert\phis\vert<1$.
Upon squaring both sides of~\eqref{eq:Fratio}, which introduces new roots, one obtains
%
\begin{equation}
   P(\amp)=0,
\label{eq:PAMP}
\end{equation}
where $P(\amp)$ is a polynomial (as described below).  We have thus shown the following:
\begin{proposition}
  Any solution of $f(\amp)=0$ is also a root of the polynomial $P(\amp)$. 
\end{proposition}
All of the roots of $P(\amp)$ can be found numerically using well-developed
algorithms (e.g.~solving the eigenvalue problem for the companion matrix).  
In
deriving $P(\amp)$, however, new roots have been introduced. Therefore, each
root $\amp^*$ of $P(\amp)$ must be checked to ensure that it is also solution of
$f(\amp)=0$.
\begin{proposition}
  A root, $\amp^*$, of $P(\amp)$ is also a root of $f(\amp)$, if $\phis$ in the definition of $f$ satisfies
\begin{equation*}
    \phis = \begin{cases}
                \phi_-(\amp^*) & \text{if $\vert\phi_-(\amp^*)\vert<1 $}, \\
                \phi_+(\amp^*) & \text{otherwise} ,
            \end{cases}
\end{equation*}
and if $f(\amp^*)=0$.
\end{proposition}
%
\newcommand{\Pamp}{P_{\rm AMP}}
The polynomial in~\eqref{eq:PAMP} is (after removing uninteresting factors of $\amp-1$) a polynomial of degree 5 given by 
\[
\Pamp(\amp)=\alpha_5\amp^5+\alpha_4\amp^4+\alpha_3\amp^3+\alpha_2\amp^2+\alpha_1\amp+\alpha_0,
\]
where the coefficients in the polynomial are
\[
\begin{array}{l}
\alpha_0=-\eta^2+2 \eta^2 \ly, \smallskip\\
\alpha_1=-8 \eta^2 \ly+5 \eta^2-\lx^2 \eta^2, \smallskip\\
\alpha_2=10 \eta^2 \ly-9 \eta^2+1+3 \lx^2 \eta^2+4 \eta \ly-2 \eta-2 \ly, \smallskip\\
\alpha_3=-2 \ly \lx^2-2 \eta \lx^2+4 \ly+\lx^2+4 \eta \ly \lx^2+7 \eta^2-3+6 \eta-4 \eta^2 \ly-8 \eta \ly-2 \eta^2 \ly \lx^2-2 \lx^2 \eta^2, \smallskip\\
\alpha_4=3-2 \ly-2 \eta^2-\lx^2+4 \eta \ly-6 \eta+2 \eta \lx^2, \smallskip\\
\alpha_5=-1+2 \eta.
\end{array}
\]
It turns out that $\Pamp(\amp)$ is the same polynomial generated by the Maple symbolic algebra program when asked to solve $f(\amp)=0$.

To determine if there are any roots of $\Pamp(\amp)$ with $\vert\amp\vert>1$, we
evaluate the roots numerically for an array of parameter values in the region
$0<\lx<1.2$, $0<\ly<1.2$ and $0\le \eta \le 1$ with $N_\amp$ equally spaced
points in each parameter direction.  The stability region of the scheme is the
region in the space of $(\lx,\ly,\eta)$ where the magnitude of the largest
valid root satisfies $\vert\amp\vert\le 1$.  This region is computed on a grid with
$N_\amp=800$ points for each parameter.  In
Figure~\eqref{fig:acousticSolidStokestability}, we plot the curve in the
$(\lx,\ly)$ plane for which $\max_{0\le\eta\le1}\vert\amp\vert=1$.  (The
curve for which $\max_{0\le\eta\le1}\vert\amp\vert=1.1$ is also plotted for
reference.)  Observe that the region of stability contains the region
$0\le\eta\le1$ and $\lx^2 + \ly^2 \le \kappa^2$
(for $\lx\ge0$ and $\ly\ge 0$) where the maximum value for $\kappa$ appears to be $\kappa=1$. 
Further evidence for $\kappa=1$ being the upper bound is provided by the restriction of the analysis to one dimension given below.
The curve $\lx^2 + \ly^2 = (.99)^2$
is shown in the figure as an example.  Since $\lx=\cp\kx\dt$ and $\ly=\cp\dt/\dy$,
the condition $\lx^2+\ly^2\le 1$ implies the following result.

\begin{theorem}
By numerical evaluation of the stability polynomial, $\Pamp(\amp)$,
we have found that the AMP algorithm is weakly stable (i.e.~there are no valid roots of $\Pamp(\amp)$ with $\vert\amp\vert>1$) 
provided $\lx^2+\ly^2\le 1$, i.e.
 \begin{align*}
    \dt \le \frac{1}{\cp}\left[ \frac{1}{\dy^2} + \kx^2 \right]^{-1/2}. 
 \end{align*}
This is a sufficient but not a necessary condition. 
\end{theorem}

{
\newcommand{\figWidth}{8cm}
\begin{figure}
\begin{center}
  \includegraphics[width=\figWidth]{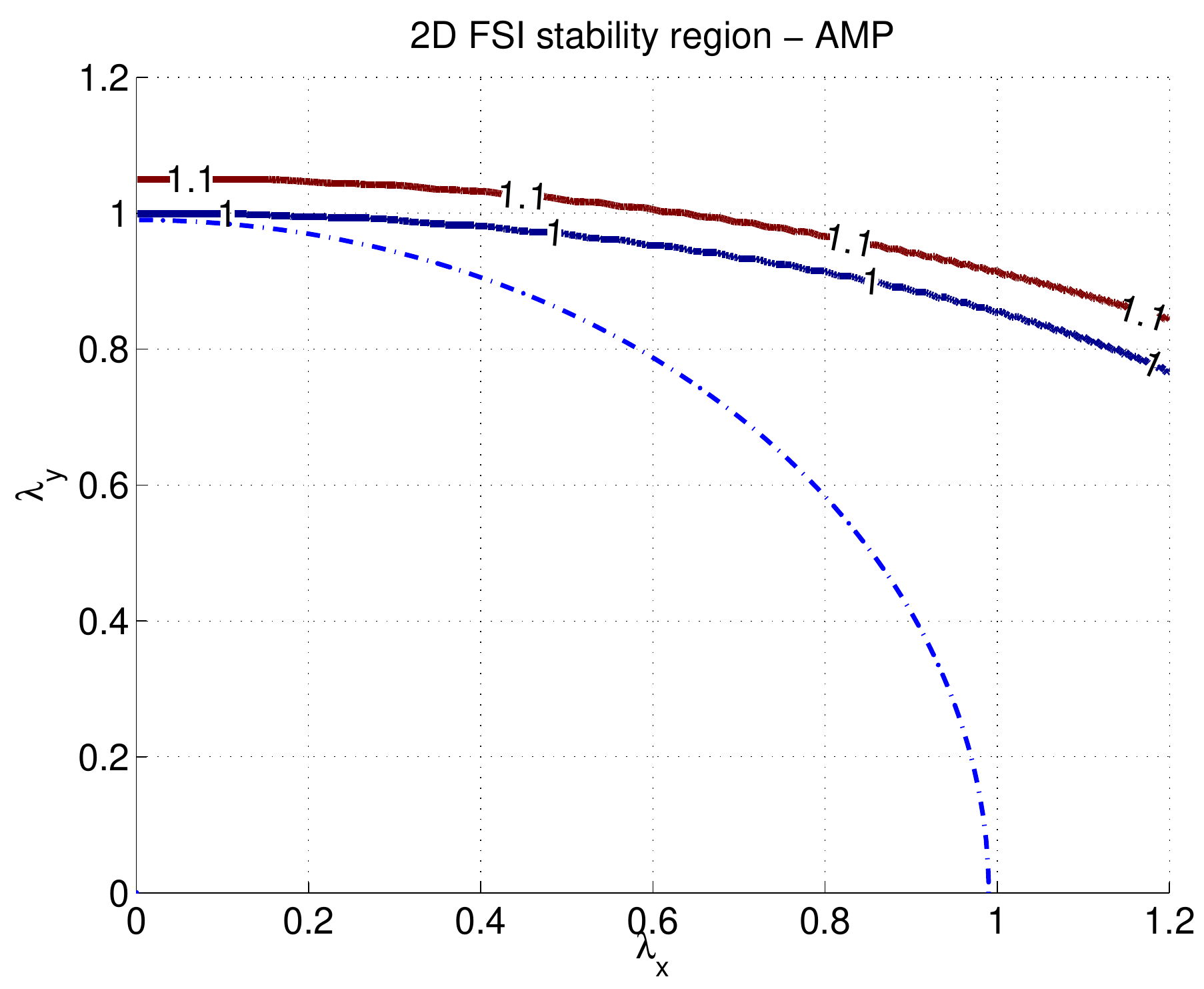}
  \caption{Stability region for the two-dimensional model problem using the AMP algorithm.  The curve marked 1 is the curve for which $\max_{0\le\eta\le1}\vert\amp\vert=1$.  The curve for which $\max_{0\le\eta\le1}\vert\amp\vert=1.1$ is shown for reference.  The dashed curve shows the curve $\lx^2+\ly^2=(.99)^2$.}
  \label{fig:acousticSolidStokestability}
\end{center}
\end{figure}
}



\paragraph{Restriction to one space dimension} 
The stability of the one-dimensional version of the AMP algorithm 
may be examined by considering the special case of the two-dimensional algorithm with $\lx=0$ and $Q=0$.  In this case there are two solutions of $f(\amp)=0$ given by
\begin{align*}
   \amp = 1-2\ly \qquad\text{or}\qquad \amp =\eta. 
\end{align*}
We are interested in whether there are solutions with $\vert\amp\vert>1$.  Since $0\le \eta\le 1$ we have the result.
\begin{theorem}
The one-dimensional AMP algorithm is weakly stable provided  $0\le\gamma\le 2$ and $0 \le \ly \le 1$, i.e.
 \begin{align*}
    \dt \le \frac{\dy}{\cp}. 
 \end{align*}
\end{theorem}


For comparison, we now discuss briefly the corresponding stability results for the scheme given by~\eqref{eq:aEqn}--\eqref{eq:sigmaEqn} with boundary conditions~\eqref{eq:BCtop} and~\eqref{eq:BCbottom}, and with interface conditions given by either the 
traditional (TP) algorithm or the anti-traditional algorithm.

\subsection{TP interface conditions} \label{sec:TPanalysis}

The interface conditions for the TP algorithm are given by~\eqref{eq:1dTradv} and~\eqref{eq:1dTradsigma}.  
Following the previous analysis, $\amp$ is found to satisfy 
\begin{align}
&   f_{{\rm TP}}(\amp) \equiv \amp(Q+1) - \frac{\eta-1}{\eta} (\amp-1)(Q-1)\phis =0 . \label{eq:traditionalSchemeF}
\end{align}


\paragraph{Restriction to one space dimension} 

In one dimension, the condition in~\eqref{eq:traditionalSchemeF} reduces to the polynomial equation
\begin{align*}
P_{{\rm TP}}(\amp) = \amp^2 + (\Mr + \ly-1) \amp - \Mr = 0, \quad\text{where}~~
  \Mr = \frac{\rho H}{\rhos \dy}.
\end{align*}
(Here, all of the roots of this quadratic are valid since no spurious roots were introduced in deriving the equation.)  The product of the roots is $-\Mr$ and thus a necessary condition for stability is that $\vert \Mr\vert\le 1$.
Thus the mass of the solid in a thin strip one grid cell wide, $\rhos\dy L$, 
must be larger than the mass of the entire fluid domain $\rho H L$. 
Using the theory of von Neumann polynomials~\cite{Miller1971,Strikwerda89} one can show the following necessary and
sufficient condition.  
\begin{theorem} \label{lemma:traditional}
The one-dimensional TP algorithm is weakly stable if and only if 
\begin{align*}
     \dt \le \frac{2}{\cp}\Big( \dy - \frac{\rho H}{\rhos} \Big). 
\end{align*}
\end{theorem}

This is a curious result.  For a given {\em large} value of $\dy$ the one-dimensional scheme may be stable, but it will eventually become unstable as the mesh is refined and $\dy$ becomes less than $(\rho H)/\rhos$.  The one-dimensional TP algorithm is thus formally {\em unconditionally unstable}.  
The numerical results in Section~\ref{sec:travelingWaveResults} show that this behaviour also holds for 
the second-order accurate discretization of the two-dimensional model problem MP-VE.

\subsection{Anti-traditional interface conditions} \label{sec:ATPanalysis}

The interface conditions for the anti-traditional algorithm are given by~\eqref{eq:1dAntiTradsigma} and~\eqref{eq:1dAntiTradv}.  The analysis for this choice leads to the constraint
\begin{align}
&   f_{{\rm AT}}(\amp) = 1- \amp + \half\left( \frac{2\eta-1}{\eta} \amp -1 \right) (Q+1)\phis=0 .  \label{eq:antiTraditional2d}
\end{align}

\paragraph{Restriction to one space dimension} 
In one dimension, \eqref{eq:antiTraditional2d} reduces to the polynomial constraint
\begin{align*}
     P_{{\rm AT}}(\amp) = \amp^2 + \left( \frac{\ly}{2}\left(1+\frac{\ly}{\Mr}\right) -2\right)\amp + 1-\frac{\ly}{2} = 0.  
\end{align*}
An analysis of the roots of $P_{{\rm AT}}(\amp)$ leads us to the following result:
\begin{theorem}
  The one-dimensional anti-traditional algorithm is weakly stable if and only if $\ly\le 4$ and 
\begin{align*}
   \dt  \le \frac{\dy}{\cp}\left( \sqrt{\Mr^2+8\Mr}-\Mr \right). 
\end{align*}
\end{theorem}

We observe that the one-dimensional anti-traditional scheme is (weakly) stable for any value of the added-mass ratio $\Mr=(\rho H)/(\rhos \dy)$ provided that $\dt$ is taken sufficiently small. 





%% file: texFiles/algorithm.tex
\newcommand{\Lsv}{\bar{\Lv}}
\newcommand{\bc}[1]{\mbox{\bfss#1}}   
\newcommand{\cc}[1]{\mbox{$//$  #1}}  
\newcommand{\ia}{\quad}        
\newcommand{\ib}{\ia\quad}     
\newcommand{\ic}{\ib\quad}     
\newcommand{\id}{\ic\quad}     
\newcommand{\ie}{\id\quad}     

\newcommand{\FUNC}[1]{{\color{blue}#1}}
\newcommand{\RETURN}{{\color{blue}Return}}
\newcommand{\IF}{{\color{blue}if}}
\newcommand{\ELSE}{{\color{blue}else}}
\newcommand{\ELSEIF}{{\color{blue}else if}}
\newcommand{\FOR}{{\color{blue}for}}
\newcommand{\COM}[1]{{\color{purple}\em #1}}
%
\section{The AMP FSI time-stepping algorithm} \label{sec:timeSteppingAlgorithm}

We now describe our implementation of the AMP algorithm, which is then used to
obtain the numerical results presented in Section~\ref{sec:numericalResults} for
the three model problems defined in Section~\ref{sec:modelProblems}.  The AMP
algorithm was introduced in predictor form in Section~\ref{sec:partitionedAlgorithms}.
It is defined here in a more general form, as a predictor-corrector time-stepping scheme, in terms of
various procedures which we outline below; additional details of the
procedures are given in~\ref{sec:timeSteppingProcedures}.  As mentioned earlier,
the velocity-pressure form of the fluid equations are solved using a
fractional-step scheme following~\cite{ICNS,splitStep2003}.  The viscous term,
$\mu\Delta\vv$, in the fluid momentum equation is advanced explicitly in time in
the present implementation, but this term could also be treated implicitly with
minor modifications.  The fractional-step scheme for the fluid also uses a
predictor-corrector scheme, consisting of a second-order accurate
Adams-Bashforth predictor followed by a second-order accurate Adams-Moulton
corrector (trapezoidal rule).  Other time-stepping schemes for the fluid could
also be used.  We note that the AMP algorithm is stable even with no corrector
step provided the predictor step in the fluid is stable in isolation.  (The use
of the Adams-Bashforth predictor for the fluid by itself requires sufficient
physical or numerical dissipation for stability.)  The solid is advanced
explicitly as well using a second-order accurate upwind scheme for the solid
equations written as a first-order
system~\eqref{eq:solidDisplacement}--\eqref{eq:solidStress}, following the
approach in~\cite{smog2012}.

Let the grid functions $p_\iv^{n-1}\approx p(\xv_\iv,t^{n-1})$ and $\vv_\iv^{n-1}\approx \vv(\xv_\iv,t^{n-1})$ denote the discrete approximations to the fluid pressure and velocity, respectively, at grid points $\xv_\iv$ and at time $t^{n-1}$.  Let $\qsv^{n-1}_\iv =( \uvs^{n-1}_\iv,\, \vvs^{n-1}_\iv,\,\sigmasv^{n-1}_\iv)^T$ denote a vector of grid functions for the solid displacement, velocity and stress, respectively.

\vskip\baselineskip
\noindent {\bf Begin predictor.}
\medskip

\noindent Stage 1: Advance the solid one time-step to give a vector of predicted values $\qsv^{(p)}_\iv =( \uvs^{(p)}_\iv,\, \vvs^{(p)}_\iv,\,\sigmasv^{(p)}_\iv)^T$.  Assign values at interior, boundary and interface points, including values for $\sigmasv_\iv^{(p)}\nv_\iv$ on the interface. 
\begin{flushleft}\tt\small
 \ia $\qsv^{(p)}_\iv$ = \FUNC{advanceSolid}( $\qsv_\iv^{n-1}$ )
\end{flushleft}


\noindent Stage 2(a): Advance the fluid velocity one time-step to give predicted values $\vv_\iv^{(p)}$.  Assign values at interior, boundary and interface points.

\begin{flushleft}\tt\small
 \ia $\vv_\iv^{(p)}$ = \FUNC{advanceFluid}( $\vv_\iv^{n-1}$, $p_\iv^{n-1}$, $\vv_\iv^{n-2}$, $p_\iv^{n-2}$ )
\end{flushleft}

\noindent Stage 2(b): Compute extrapolated values for the fluid pressure and compute the projected interface velocity excluding traction terms in the formula for the interface velocity since these terms are not known at this stage.  Assign boundary conditions on the fluid velocity using the AMP Robin condition~\eqref{eq:solidVtangentialVsubproblem} for the velocity.
 
\begin{flushleft}\tt\small
\ia $p_\iv^{(e)}$ = \FUNC{extrapolateInTime}( $p_\iv^{n-1}$, $p_\iv^{n-2}$ $p_\iv^{n-3}$ ) \\ 
\ia $\vv_\iv^I$ = \FUNC{projectInterfaceVelocity}( $\vv_\iv^{(p)}$, $p_\iv^{(e)}$, $\vsv_\iv^{(p)}$, $\sigmasv_\iv^{(p)}\nv_\iv$, $\beta=0$ ) \\
\ia $(\vv_\iv^I,(\sigmav\nv)_\iv^I)$ = \FUNC{assignFluidVelocityBoundaryConditions}( $\vv_\iv^{(p)}$, $p_\iv^{(e)}$, $\vsv_\iv^{(p)}$, $\sigmasv_\iv^{(p)}\nv_\iv$ ) 
\end{flushleft}
Here $p_\iv^{(e)}$, needed in the boundary conditions, is a second-order accurate approximation for the fluid pressure on the interface at time $t^{n}$ obtained by extrapolation in time using past values $p_\iv^m$, $m=n-1,n-2,n-3$.  The parameter $\beta$ in the procedure that computes the interface velocity is set to zero which specifies that the traction terms are not used.

\medskip\noindent Stage 3(a): Compute the solid acceleration, $\dot{\vvs}_\iv^{(p)}$, and solve for the fluid pressure using
the AMP Robin pressure boundary condition~\eqref{eq:AddedMassPressureBC}.  Assign values at interior, boundary and interface points, including ghost points. 
\begin{flushleft}\tt\small
\ia $\dot{\vvs}_\iv^{(p)}$ = \FUNC{computeSolidAcceleration}( $\vvs_\iv^{(p)}$, $\vvs_\iv^n$  )  \\
\ia $p_\iv^{(p)}$ = \FUNC{solveFluidPressureEquation}( $\vv_\iv^{(p)}$, $\vsv_\iv^{(p)}$, $\sigmasv_\iv^{(p)}\nv_\iv$, $\dot{\vvs}_\iv^{(p)}$ ) \\
\end{flushleft}
Here $\dot{\vvs}_\iv^{(p)}$, used for the pressure boundary condition, is a second-order accurate approximation to the solid acceleration on the interface at $t^{n}$.

\medskip\noindent Stage 3(b) :
Given $p_\iv^{(p)}$, recompute the interface velocity and traction, and then assign solid boundary conditions.
\begin{flushleft}\tt\small
\ia $\vv_\iv^I$ = \FUNC{projectInterfaceVelocity}( $\vv_\iv^{(p)}$, $p_\iv^{(p)}$, $\vsv_\iv^{(p)}$, $\sigmasv_\iv^{(p)}\nv_\iv$, $\beta=1$ ) \\
\ia $(\vv_\iv^I,(\sigmav\nv)_\iv^I)$ = \FUNC{assignFluidVelocityBoundaryConditions}( $\vv_\iv^{(p)}$, $p_\iv^{(p)}$, $\vsv_\iv^{(p)}$, $\sigmasv_\iv^{(p)}\nv_\iv$ ) \\
\ia $(\vsv_\iv^{(p)},\sigmasv^{(p)}_\iv\nv_\iv)$ = \FUNC{assignSolidBoundaryConditions}( $\vsv_\iv^{(p)}$, $\sigmasv_\iv^{(p)}$, $\vv_\iv^I$, $(\sigmav\nv)_\iv^I$ ) 
\end{flushleft}
The parameter $\beta$ in the procedure that computes the interface velocity is now set to one which specifies that the traction terms are used.

\medskip\noindent {\bf End predictor.}

\bigskip
If the predictor step is used alone (without the corrector step), then discrete
values for $p_\iv^{n}$ and $\vv_\iv^{n}$ in the fluid, and for $\qsv^{n}_\iv =(
\uvs^{n}_\iv,\, \vvs^{n}_\iv,\,\sigmasv^{n}_\iv)^T$ in the solid, are taken from
the corresponding predicted values and the time step to $t^n$ is complete.  If,
on the other hand, a corrector step is used, then stages 4--7 are included in
the algorithm as given below.

\vskip\baselineskip
\noindent {\bf Begin corrector}: (optional)


\medskip\noindent Stage 4 : Correct the solid to obtain $\qsv^{n}_\iv =( \uvs^{n}_\iv,\, \vvs^{n}_\iv,\,\sigmasv^{n}_\iv)^T$.  Assign values at interior, boundary and interface points. 
\begin{flushleft}\tt\small
 \ia $\qsv^{n}_\iv$ = \FUNC{correctSolid}( $\qsv^{(p)}_\iv$, $\qsv_\iv^{n-1}$ )
\end{flushleft}

\noindent Stage 5(a) : Correct the fluid velocity to obtain $\vv_\iv^{n}$.  Assign interior, boundary and interface points.
\begin{flushleft}\tt\small
 \ia $\vv_\iv^{n}$ = \FUNC{correctFluid}( $\vv^{(p)}_\iv$, $p_\iv^{(p)}$, $\vv_\iv^{n-1}$, $p_\iv^{n-1}$ ) 
\end{flushleft}


\noindent Stage 5(b) : Assign fluid boundary conditions (assigns boundary/interface and ghost points).
\begin{flushleft}\tt\small
 \ia $\vv_\iv^I$ = \FUNC{projectInterfaceVelocity}( $\vv_\iv^{n}$, $p_\iv^{(p)}$, $\vsv_\iv^{n}$, $\sigmasv_\iv^{n}\nv_\iv$, $\beta=1$  ) \\
\ia $(\vv_\iv^I,(\sigmav\nv)_\iv^I)$ = \FUNC{assignFluidVelocityBoundaryConditions}( $\vv_\iv^{n}$, $p_\iv^{(p)}$, $\vv_\iv^I$, $\vsv_\iv^{n}$, $\sigmasv_\iv^{n}\nv_\iv$ ) 
\end{flushleft}

\noindent Stage 6: Correct the fluid pressure.  Assign values at interior, boundary and interface points, including ghost points.
\begin{flushleft}\tt\small
\ia $\dot{\vvs}_\iv$ = \FUNC{computeSolidAcceleration}( $\vvs_\iv^{n}$, $\vvs_\iv^{n-1}$ )  \\
\ia $p_\iv^{n}$ = \FUNC{solveFluidPressureEquation}( $\vv_\iv^{n}$, $\vsv_\iv^{n}$, $\sigmasv_\iv^{n}\nv_\iv$, $\dot{\vvs}_\iv$ ) \\
\end{flushleft}

\noindent Stage 7 :  Re-compute interface velocity using corrected pressure, and re-assign the fluid boundary conditions. 
Assign solid boundary conditions using latest values for $\vv_\iv^I$ and $\sigmav_\iv^I$. 
\begin{flushleft}\tt\small
 \ia $\vv_\iv^I$ = \FUNC{projectInterfaceVelocity}( $\vv_\iv^{n}$, $p_\iv^{n}$, $\vsv_\iv^{n}$, $\sigmasv_\iv^{n}\nv_\iv$, $\beta=1$ ) \\
\ia $(\vv_\iv^I,(\sigmav\nv)_\iv^I)$ = \FUNC{assignFluidVelocityBoundaryConditions}( $\vv_\iv^{n}$, $p_\iv^{n}$, $\vv_\iv^I$, $\vsv_\iv^{n}$, $\sigmasv_\iv^{n}\nv_\iv$ ) \\
\ia $(\vsv_\iv^{n},\sigmasv^{n}_\iv\nv_\iv)$ = \FUNC{assignSolidBoundaryConditions}( $\vsv_\iv^{n}$, $\sigmasv_\iv^{n}$, $\vv_\iv^I$, $(\sigmav\nv)_\iv^I$ )
\end{flushleft}


\noindent {\bf end corrector}

       
\vskip\baselineskip
The corrector step described above may be repeated, replacing the predicted states with
the latest solution values. This may permit a somewhat larger time-step but also
involves an additional cost. 
For all calculations
presented in the next section we use the predictor step with just one corrector step.

%% file: texFiles/results.tex
\renewcommand{\tableFontSize}{\scriptsize}
\section{Numerical results} \label{sec:numericalResults}

In this section, we present numerical results based on the model problems in Section~\ref{sec:modelProblems} to verify the accuracy and stability of the AMP algorithm.  All of the model problems consist of a rectangular fluid domain, $\OmegaF=(0,L)\times(-H,0)$, and a rectangular solid domain, $\OmegaS=(0,L)\times(0,\Hs)$, connected by an interface, $\GammaI=\{ (x,y) \,\vert\, x\in(0,L),  y=0 \}$.
We use a Cartesian grid for the fluid domain with $N_j+1$ grid points in each direction so that the grid spacings are $\dx_j=L/N_j$ and $\dy_j=H/N_j$.  (The subscript $j$ is used later to indicate the resolution of the grid.)  The grid for the solid domain is also a Cartesian grid and uses the same grid spacing as the fluid grid.  For all calculations we take
\begin{align*}
&    \rho=1,\qquad L=1, \qquad H=1, \qquad \rhos=\lambdas=\mus=\rho\, \dsf,\qquad\Hs=1/2,
\end{align*}
where the {\em density ratio} $\dsf=\rhos/\rho$ is a parameter that is chosen later to consider the cases of light, medium and heavy solids. The value for the fluid viscosity, $\mu$, is chosen based on the different model problems.
The value for the fluid impedance, $\zf$, used in the velocity projection~\eqref{eq:ImpedenceWeightedVn} and~\eqref{eq:ImpedenceWeightedVt}, is taken as $\zf=\rho\dy_j/\dt_j$ following the discussion in Section~\ref{sec:AMP1D}, however the results are insensitive to this choice.

We begin by considering solutions of the model problems constructed using the
method of analytic solutions. These are compared with numerical approximations
computed using the AMP algorithm for different values of density
ratio.  We then consider numerical solutions of the model problems for cases
where exact traveling wave solutions are known.

\input texFiles/tzErrorTables

\subsection{The method of analytic solutions} \label{sec:TZ}

The method of analytic solutions is a useful technique for constructing exact
solutions of initial-boundary-value problems for partial differential equations
for the purpose of checking the behavior and accuracy of the numerical
implementation of a problem.  This method, also known as the {\it method of
manufactured solutions}~\cite{Roache2002} or {\it twilight-zone
forcing}~\cite{CGNS}, adds forcing functions to the governing equations, 
boundary conditions and interface conditions. These forcing functions are specified so that a chosen
function, $\tilde{\qv}(\xv,t)$, becomes the exact solution of the forced equations,
and thus the error in the discrete solution can be computed exactly.

The method of analytic solutions is applied to the FSI problems using
trigonometric functions.   The exact solutions for the components of displacement, velocity and stress in the solid
are taken to be 
\begin{alignat}{4}
\tilde \us_1&=~~.25\, c_1 \tilde{c}_2 c_t,\qquad &
\tilde \us_2&=~.5\, c_1 \tilde{c}_2 c_t,\qquad&
\tilde \vs_1&=\dot{\us}_1,\qquad&
\tilde \vs_2&=\dot{\us}_2,~ \label{eq:tzTrig2d}\\
\tilde \sigmas_{11}&=-.5\, c_1 c_2 c_t,~&
\tilde \sigmas_{12}&=~.4\, s_1 c_2 c_t,~&
\tilde \sigmas_{21}&=~.4\, s_1 c_2 c_t,\qquad&
\tilde \sigmas_{22}&=~.6\, c_1 s_2 c_t , \nonumber
\end{alignat}
where $(c_1,s_1)=(\cos(2\pi x),\sin(2\pi x))$, $(c_2,s_2)=(\cos(2\pi y),\sin(2\pi y))$, $\tilde{c}_2=\cos(2\pi(y+0.375))$
and $c_t=\cos(2\pi t)$.
The exact solutions for the pressure and the components of velocity in the fluid are
\begin{alignat}{4}
\tilde p&= c_1 s_2 c_t,\qquad
\tilde v_1&=.5 c_1 c_2 c_t, \qquad
\tilde v_2&=.5 s_1 s_2 c_t.
\end{alignat}
These trigonometric functions are substituted into the differential equations governing each domain, and into the boundary and interface conditions, to define forcing functions to the equations so that the functions given above become exact solutions of the initial-boundary-value problem.

Numerical solutions of the (forced) model problems are computed using the AMP
algorithm.  The initial conditions for the numerical calculations are taken from
the exact solutions evaluated at $t=0$.  For each case, maximum-norm errors,
$E_j^{(q)}$, for solution component $q$, are computed on grids of increasing
resolution using grid spacings $\dx_j=\dy_j=h_j=1/(20j)$, $j=1,2,\ldots$.  The
convergence rate, $\zeta_q$, is estimated by a least squares fit to the
logarithm of the error equation, $E_j^{(q)} = C_q h_j^{\zeta_q}$, where $C_q$ is
approximately constant for small grid spacings.  For vector variables, such as
$\vv$ or $\sigmasv$, the error denotes the maximum over all components of the
vector or tensor.

The tables in Figure~\ref{tab:MP-VA_AMP_TZ_trig} give the maximum errors and estimated convergence rates
for the results of the AMP algorithm applied to the FSI problem of a viscous incompressible fluid and acoustic solid that only supports vertical motion (model problem MP-VA).  Similar results are obtained for the other model problems. 
Results are presented for the density ratios $\delta=10^3$, $1$ and $10^{-1}$, referred to as {\em heavy}, {\em medium} and {\em light} solids. 
The fluid viscosity is taken as $\mu=0.05$ for each case. 
The values in the columns labeled ``r'' give the ratio of the error on the current grid to that on the previous coarser grid,
a ratio of $4$ is expected for a second-order accurate method.
The values in the tables show that the scheme is stable and close to second-order accurate in the max-norm. 
For the heavy solid, $\delta=10^3$, the stresses in the solid are very large, $\Oc(\rhos\cp^2)$, compared to the fluid stress,
and this affects the convergence rates on coarser grids.  Note that the errors given in the table (and subsequent tables below) are absolute errors. 
Thus, for the case of a heavy solid, the errors in the solid stress appear to be large even though the corresponding relative errors are small.
\input texFiles/insSolidFig

\subsection{Traveling wave exact solutions} \label{sec:travelingWaveResults}

We now consider the three FSI model problems for cases when the exact solutions
are traveling waves.  The traveling waves have the form $q(x,y,t) =
\hat{q}(y)e^{i(kx-\omega t)}$, where $q$ represents any component of the
variables belonging to the fluid or solid domains.  The solutions are periodic
in $x$ with wave number $k\in\Real$ and have frequency $\omega\in\Complex$ in
time.  For inviscid fluids (as in the MP-IA model problem), ${\rm Im}(\omega)=0$
so that the traveling wave is also periodic in time.  For viscous fluids (as in
the model problems MP-VA and MP-VE), ${\rm Im}(\omega)<0$ so that the amplitude
of the traveling wave decays over time.  For the results shown here, we take
$k=2\pi$ and a value for $\omega$ obtained from solving an appropriate dispersion relation
(see~\ref{sec:travelingWave}) such that the wave travels from left to right
(i.e.~${\rm Re}(\omega)>0$).  In all cases, periodic boundary conditions are
used in the $x$-direction.  The amplitude parameter in the definition of the
traveling wave solutions is chosen as $\ampe=1/10$, which defines the maximum
amplitude of the displacement on the interface.  The initial conditions for the
numerical solutions are taken from the exact traveling wave solution evaluated
at $t=0$.


{
\newcommand{\labelFont}{\scriptsize}
\newcommand{\xa}{1.6}
\newcommand{\ysa}{2.5}
\newcommand{\ysb}{5.5}
\begin{figure}[hbt]
\newcommand{\figWidth}{7.cm}
\newcommand{\trimfig}[2]{\trimFig{#1}{#2}{0.0}{0.}{.0}{.0}}
\begin{center}
\begin{tikzpicture}[scale=1]
  \useasboundingbox (0,.75) rectangle (14.5,6.5);  
  \draw (0,0) node[anchor=south west,xshift=-4pt,yshift=+0pt] {\trimfig{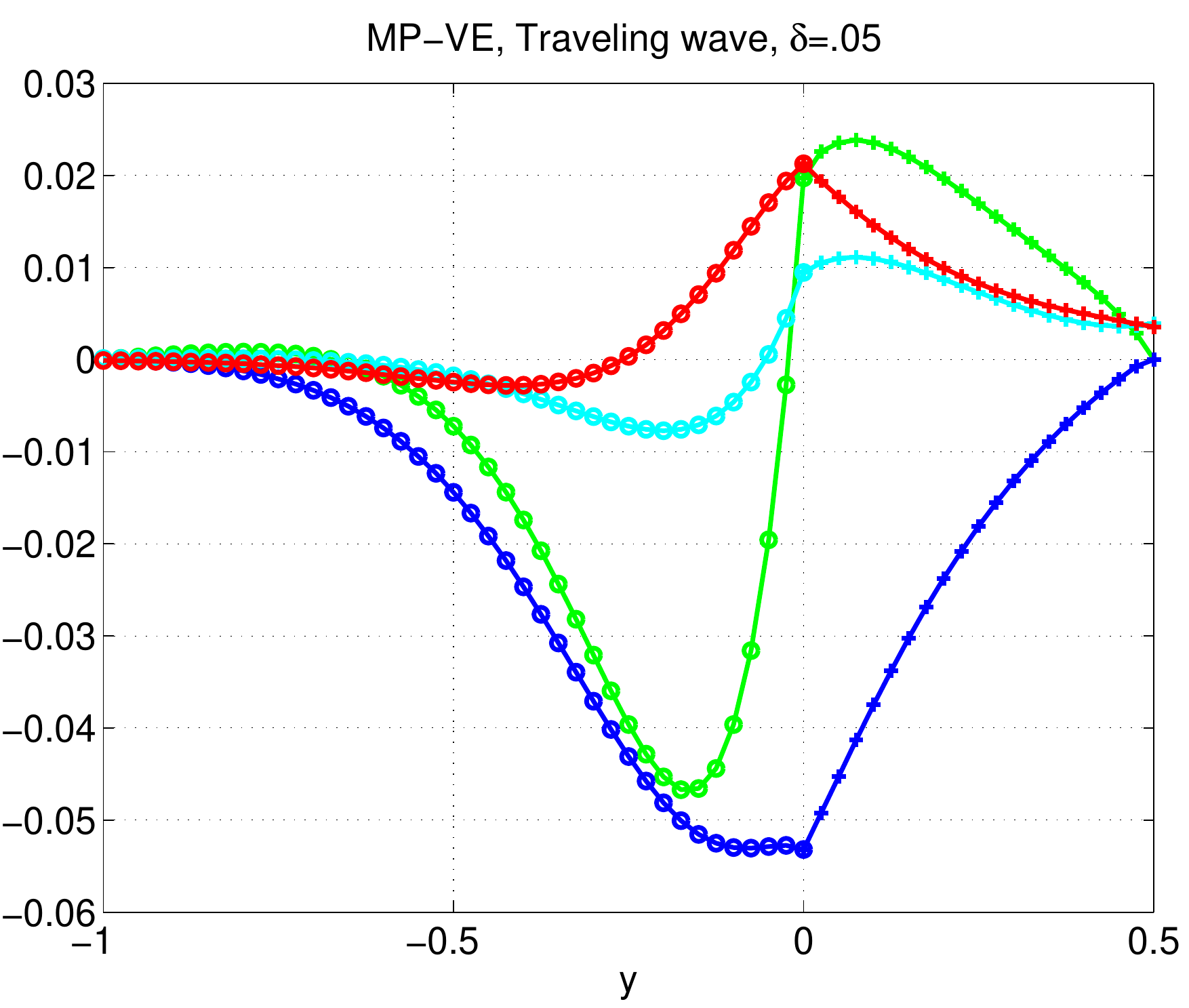}{\figWidth}};
  \draw (7.5, 0) node[anchor=south west,xshift=-4pt,yshift=+0pt] {\trimfig{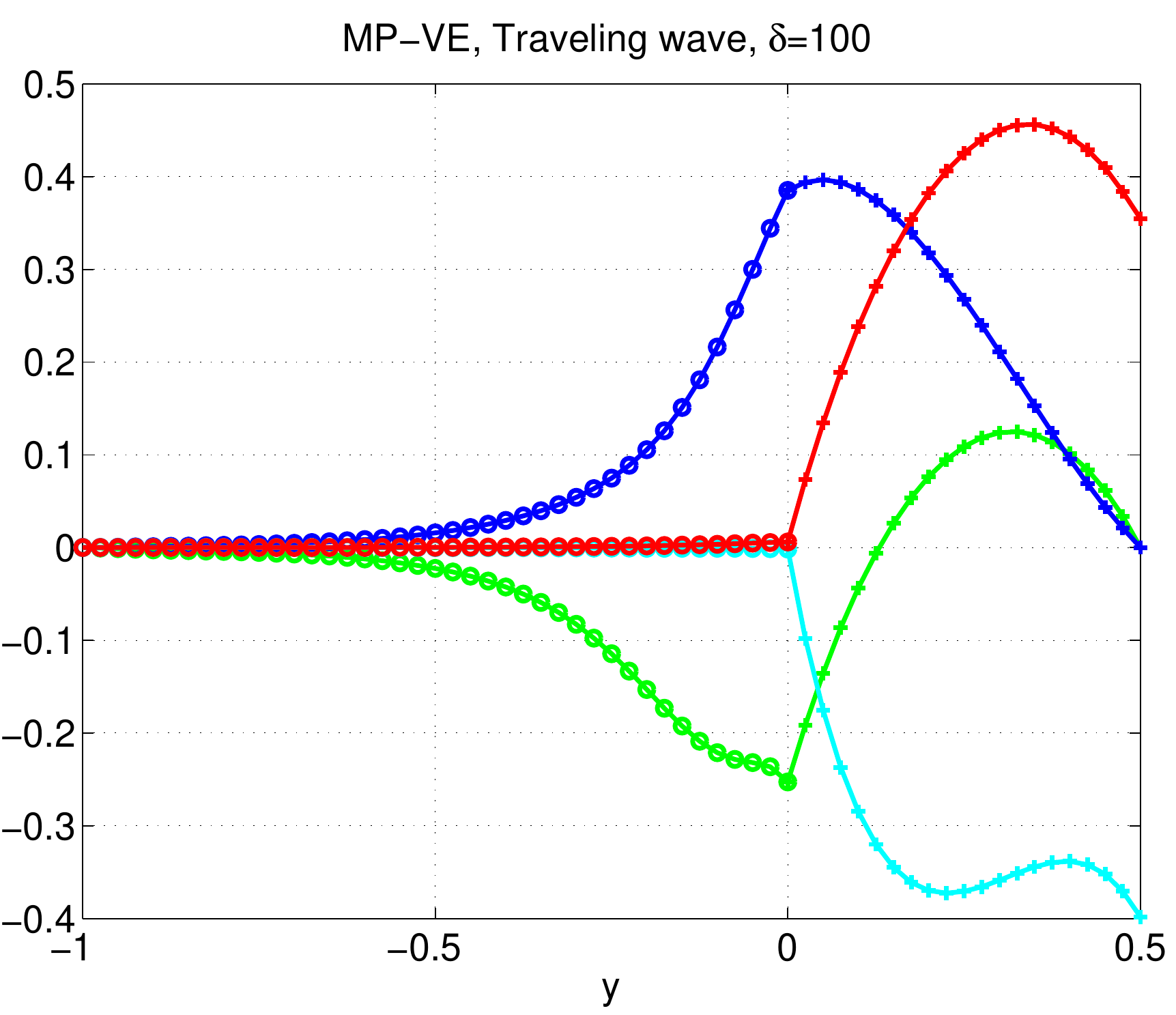}{\figWidth}};
  \draw (4.5,2.0) node[anchor=south east,fill=green!30,xshift=-2pt,yshift= 9pt,inner sep=1pt] {\labelFont$v_1$}; 
  \draw (3.4,1.9) node[anchor=north east,fill=blue!30,xshift=-2pt,yshift= 9pt,inner sep=1pt] {\labelFont$v_2$}; 
  \draw (4.,3.55) node[anchor=north,fill=cyan!30,xshift=-2pt,yshift=-2pt,inner sep=1pt] {\labelFont$\sigma_{21}$}; 
  \draw (4.1,4.4) node[anchor=south east,fill=red!30,xshift=-2pt,yshift=-2pt,inner sep=1pt] {\labelFont$\sigma_{22}$}; 
  \draw (11.,2.1) node[anchor=north east,fill=green!30,xshift=-2pt,yshift= 9pt,inner sep=1pt] {\labelFont$v_1$}; 
  \draw (11.6,3.5) node[anchor=south east,fill=blue!30,xshift=-2pt,yshift= 9pt,inner sep=1pt] {\labelFont$v_2$}; 
  \draw (12.9,1.2) node[anchor=south west,fill=cyan!30,xshift=-2pt,yshift=-2pt,inner sep=1pt] {\labelFont$\sigma_{21}$}; 
  \draw (14.0,5.15) node[anchor=south east,fill=red!30,xshift=-2pt,yshift=-2pt,inner sep=1pt] {\labelFont$\sigma_{22}$}; 
%
%
\end{tikzpicture}
\end{center}
  \caption{Plots of the computed solution along the line $x=2/3$ for the traveling wave solution and model problem MP-VE, for $\delta=0.05$ (left)
   and $\delta=100$ (right). The fluid solution is in the interval $y\in[-1,0]$ while the solid occupies $y\in[0,.5]$. 
   For the right-hand plot, $\sigma_{21}$ and $\sigma_{22}$ are scaled by factors of $50$ and $100$, respectively.
        These plots correspond to the solutions shown in Figure~\ref{fig:TW_EWE}.}

  \label{fig:TW_MPVE_LinePlots}
\end{figure}
}

Figure~\ref{fig:TW_EWE} shows shaded contour plots of the computed solution for
model problem MP-VE with two different density ratios at time $t=0.5$.  The
fluid viscosity is $\mu=.02$ for both cases.  The shaded contours in the upper
solid domain show the behavior of the velocity components, $\vs_1$ and $\vs_2$,
and the stress component, $\sigmas_{22}$, while the shaded contours in the lower
fluid domain show the corresponding behavior of the velocity components, $v_1$
and $v_2$, and the pressure, $p$.  The position of the displaced interface given
by the computed solid displacement, $\uvs_\iv^n$ at $y_\iv=0$ and $t^n=0.5$, is
shown in the plots as a curve superimposed on the shaded contours.  For the
light solid case with $\delta=\rhos/\rho=.05$, the corresponding exact solution
has wave number $k=2\pi$ and (complex) frequency $\omega\approx(1.290,-.5899)$,
while $\omega\approx(6.714,-6.359{\rm e-}3)$ for the heavy solid case with
$\delta=100$.  The behavior of the solution shows an oscillation in the
$x$-direction, as expected, and decay in the $y$-direction away from the
interface at $y=0$.  It can be observed that the contour plots of the velocity
components indicate that the velocity is continuous across the interface at
$y=0$ in agreement with the interface condition in~\eqref{eq:interfaceV}.  A
further verification of the interface conditions in~\eqref{eq:interfaceV}
and~\eqref{eq:interfaceStress} is given in Figure~\ref{fig:TW_MPVE_LinePlots}.
The plots in this figure show the behavior of the components of velocity and
normal stress along the line $x=2/3$ for the density ratios $\delta=.05$ (left)
and $\delta=100$ (right).  Note that each fluid-solid component pair is
continuous at $y=0$ in agreement with the interface conditions.

\begin{figure}[hbt]\tableFontSize
\begin{center}
\begin{tabular}{|c|c|c|c|} \hline
\multicolumn{4}{|c|}{Traveling wave frequencies $\omega$} \\ \hline 
~~$\delta$~~ & MP-IA        & MP-VA              & MP-VE    \\ \hline 
\strutt $10^{-1}$  &  $(15.513,0)$    &  $(2.792,-0.7469)$    &  $(1.905, -0.6524)$          \\ \hline   
\strutt $1$        &  $(16.556,0)$    &  $(8.126,-0.7261)$    &  $(5.082, -0.4619)$            \\ \hline   
\strutt $10^3$     &  $(29.294,0)$    & $(12.163,-9.730\hbox{e-4})$  &  $(6.731, -6.365\hbox{e-4})$         \\ \hline   
\end{tabular}
\end{center}
\caption{Values of the (complex) frequency $\omega=a+ib=(a,b)$ for the exact traveling wave solution used in the numerical simulations of the three model problems for $k=2\pi$, and the parameters $\rho=1$, $H=1$, $\lambdas=\mus=\rhos=\delta$ and $\Hs=1/2$.  For the models with a viscous fluid, we use $\mu=.02$.}
\label{tab:travelingWaveOmega}
\end{figure}

Traveling wave solutions computed using the AMP algorithm may be compared with
the corresponding exact solutions to further verify the stability and accuracy
of the scheme.  We consider the FSI model problems, MP-IA, MP-VA and MP-VE, for
three different density ratios, $\delta=10^3$, 1, $10^{-1}$.  The wave number is
taken to be $k=2\pi$ for each case, and the corresponding values of $\omega$ for
the exact solutions are given in Figure~\ref{tab:travelingWaveOmega}.  The
numerical solutions and corresponding max-norm errors are computed for a
sequence of decreasing grid spacings, as was done previously in
Section~\ref{sec:TZ}, and the convergence rates for the solution components are
estimated.

\input texFiles/twMPIA1ErrorTables

\input texFiles/twMPVA1ErrorTables

\input texFiles/twMPVEErrorTables

The results of this convergence study are given in Figures~\ref{tab:MP-IA1_AMP},
\ref{tab:MP-VA1_AMP} and~\ref{tab:MP-VE_AMP} for the model problems MP-IA, MP-VA
and MP-VE, respectively.  For the case of an inviscid fluid and acoustic solid
supporting vertical motion only (MP-IA), the comparison between the computed
solutions and the exact solutions is made at $t=1.0$.  The results show that the
scheme is stable and close to second-order accurate in the max-norm.  Since the
discretization of the fluid equations uses central finite differences, a
small amount of artificial dissipation is added, proportional to $h_j^2$, to the fluid
momentum equation to smooth boundary layers in the error that otherwise degrade
the max-norm convergence rates slightly (the scheme is stable without
dissipation).  For the cases of a viscous fluid (MP-VA and MP-VE) with
$\mu=.02$, solutions decay in time and the comparison is made at $t=0.3$.  For
both of these viscous cases, one with an acoustic solid supporting vertical
motion only and the other with a compressible elastic solid allowing motion in
both directions, the results show that the AMP algorithm is stable and close to
second-order accurate in the max-norm.

\subsection{Traditional partitioned scheme}
The table in Figure~\ref{tab:MP-VE_Traditional} indicates the stability of the traditional partitioned (TP) algorithm
for the model problem MP-VE for different values of $\delta=\rhos/\rho$.
The stability was determined experimentally by integrating the equations for a large number of time steps, 
and looking for exponential blowup.
For heavy enough solids one might expect the TP algorithm to be stable. 
However, from the theory for the model problem discussed in Section~\ref{sec:analysis} (see Theorem~\ref{lemma:traditional}), 
it was shown that the scheme becomes unstable on a sufficiently fine grid, for any value of $\delta$ no matter how large.
For a given ratio $\delta$, the scheme may be stable on a coarse grid but becomes unstable for all sufficiently fine grids.
The question is whether a similar result holds for the more general model problem MP-VE. 
The results in Figure~\ref{tab:MP-VE_Traditional} strongly suggest that this
behavior also holds for the more complex model problem.  For
$\delta=\rhos/\rho=100$ the scheme is stable for $h_j=1/20$ but unstable for $h_j=1/40$
and smaller. For a larger ratio, $\delta=200$, the scheme is stable for
$h_j=1/20$ and $h_j=1/40$ but unstable for $h_j=1/180$ and smaller. This trend
continues.

\begin{figure}[hbt]\tableFontSize
\begin{center}
\begin{tabular}{|c|c|c|c|c|} \hline
\multicolumn{5}{|c|}{MP-VE, traveling wave, TP algorithm} \\ \hline 
~~$h_j$~~   & $\delta=800$   & $\delta=400$ & $\delta=200$ & $\delta=100$     \\ \hline 
 1/20       & stable        &  stable     & stable      & stable     \\ \hline   
 1/40       & stable        &  stable     & stable      & unstable   \\ \hline   
 1/80       & stable        &  stable     & unstable    & unstable   \\ \hline   
 1/160	    & stable        &  unstable   & unstable    & unstable   \\ \hline   
 1/320	    & unstable      &  unstable   & unstable    & unstable   \\ \hline   
\end{tabular}
\end{center}
\caption{Stability of the traditional partitioned (TP) algorithm for the traveling wave solution for model problem MP-VE for different values of $\delta=\rhos/\rho$ with $\rho=1$.  As predicted by the theory, for any given value of $\delta$, no matter how large, the TP algorithm becomes unstable when the mesh is sufficiently fine.
}
\label{tab:MP-VE_Traditional}
\end{figure}

%% file: texFiles/tzErrorTables.tex
\bogus{
\begin{table}[hbt]\tableFont 
\begin{center}
\begin{tabular}{|l|c|c|c|c|c|c|c|c|c|c|c|} \hline 
grid  & N & $p$ & r & $\vert \mathbf{v} \vert$ & r & $\bar{\mathbf{u}}$ & r & $\bar{\mathbf{v}}$ & r & $\bar{\mathbf{s}}$& r \\ \hline 
                sq20 &     1 & \num{7.0}{-3} &      & \num{8.0}{-3} &      & \num{3.0}{-3} &      & \num{3.4}{-2} &      & \num{8.4}{-3} &       \\ \hline
                sq40 &     2 & \num{1.3}{-3} &  5.4 & \num{1.2}{-3} &  6.5 & \num{8.6}{-4} &  3.5 & \num{7.9}{-3} &  4.2 & \num{1.9}{-3} &  4.5  \\ \hline
                sq80 &     4 & \num{2.7}{-4} &  4.8 & \num{2.0}{-4} &  6.2 & \num{2.3}{-4} &  3.8 & \num{2.1}{-3} &  3.7 & \num{3.8}{-4} &  4.8  \\ \hline
               sq160 &     8 & \num{7.6}{-5} &  3.6 & \num{3.7}{-5} &  5.4 & \num{6.0}{-5} &  3.8 & \num{5.6}{-4} &  3.8 & \num{7.7}{-5} &  5.0  \\ \hline
    rate             &       &  $2.18$       &      &  $2.59$       &      &  $1.88$       &      &  $1.96$       &      &  $2.26$       &       \\ \hline
\end{tabular}
\caption{Inses bulk: tz.trig.am1.rem1, solid=(swefos,fos), AM=1, amImp=0.5, vpf=0.1, ts=pc2, nc=1, normalMotion=1, bcs=periodic, bcsTop=Dirichlet, cdv=0.25, $t=.3$, $\mu=.05$, ad2=0, $\rho=1$, $\bar{\rho}=1e-1$, scaleTZStrain=0, TZ=trig(fx=2,ft=2), cfl=0.5, Thu Jul 18 14:50:09 2013}\label{table:bulk.tz.trig.am1.rem1.pc2.swefos.fos.vector}
\end{center}
\end{table}
}
\newcommand{\tableMPVAIAmpLightTrigTZ}{%
\begin{tabular}{|c|c|c|c|c|c|c|c|c|c|c|} \hline
\multicolumn{11}{|c|}{MP-VA, trigonometric solution, light solid} \\ \hline 
\strutt~~$h_j$~~& $E_j^{(p)}$ & $r$ & $E_j^{(\vv)}$ & $r$ & $E_j^{(\usv)}$ & $r$  & $E_j^{(\vsv)}$  & $r$ & $E_j^{(\sigmasv)}$  & $r$ \\ \hline 
 1/20      & \num{7.0}{-3} &      & \num{8.0}{-3} &      & \num{3.0}{-3} &      & \num{3.4}{-2} &      & \num{8.4}{-3} &       \\ \hline
 1/40      & \num{1.3}{-3} &  5.4 & \num{1.2}{-3} &  6.5 & \num{8.6}{-4} &  3.5 & \num{7.9}{-3} &  4.2 & \num{1.9}{-3} &  4.5  \\ \hline
 1/80      & \num{2.7}{-4} &  4.8 & \num{2.0}{-4} &  6.2 & \num{2.3}{-4} &  3.8 & \num{2.1}{-3} &  3.7 & \num{3.8}{-4} &  4.8  \\ \hline
 1/160	   & \num{7.6}{-5} &  3.6 & \num{3.7}{-5} &  5.4 & \num{6.0}{-5} &  3.8 & \num{5.6}{-4} &  3.8 & \num{7.7}{-5} &  5.0  \\ \hline
\rateLabel &  $2.18$       &      &  $2.59$       &      &  $1.88$       &      &  $1.96$       &      &  $2.26$       &       \\ \hline
\end{tabular}
}

\bogus{
\begin{table}[hbt]\tableFont 
\begin{center}
\begin{tabular}{|l|c|c|c|c|c|c|c|c|c|c|c|} \hline 
grid  & N & $p$ & r & $\vert \mathbf{v} \vert$ & r & $\bar{\mathbf{u}}$ & r & $\bar{\mathbf{v}}$ & r & $\bar{\mathbf{s}}$& r \\ \hline 
                sq20 &     1 & \num{2.1}{-2} &      & \num{2.2}{-2} &      & \num{5.4}{-3} &      & \num{2.5}{-2} &      & \num{4.6}{-2} &       \\ \hline
                sq40 &     2 & \num{6.1}{-3} &  3.4 & \num{3.9}{-3} &  5.6 & \num{1.3}{-3} &  4.0 & \num{6.7}{-3} &  3.8 & \num{9.3}{-3} &  5.0  \\ \hline
                sq80 &     4 & \num{1.8}{-3} &  3.3 & \num{7.4}{-4} &  5.2 & \num{3.2}{-4} &  4.2 & \num{1.9}{-3} &  3.6 & \num{2.2}{-3} &  4.1  \\ \hline
               sq160 &     8 & \num{4.6}{-4} &  4.0 & \num{1.5}{-4} &  4.8 & \num{7.8}{-5} &  4.1 & \num{5.1}{-4} &  3.7 & \num{5.3}{-4} &  4.3  \\ \hline
    rate             &       &  $1.82$       &      &  $2.38$       &      &  $2.04$       &      &  $1.87$       &      &  $2.14$       &       \\ \hline
\end{tabular}
\caption{Inses bulk: tz.trig.am1.rep0, solid=(swefos,fos), AM=1, amImp=0.5, vpf=0.1, ts=pc2, nc=1, normalMotion=1, bcs=periodic, bcsTop=Dirichlet, cdv=0.25, $t=.3$, $\mu=.05$, ad2=0, $\rho=1$, $\bar{\rho}=1.$, scaleTZStrain=0, TZ=trig(fx=2,ft=2), cfl=0.5, Thu Jul 18 14:13:16 2013}\label{table:bulk.tz.trig.am1.rep0.pc2.swefos.fos.vector}
\end{center}
\end{table}
}
\newcommand{\tableMPVAIAmpMediumTrigTZ}{%
\begin{tabular}{|c|c|c|c|c|c|c|c|c|c|c|} \hline
\multicolumn{11}{|c|}{MP-VA, trigonometric solution, medium solid} \\ \hline 
\strutt~~$h_j$~~& $E_j^{(p)}$ & $r$ & $E_j^{(\vv)}$ & $r$ & $E_j^{(\usv)}$ & $r$  & $E_j^{(\vsv)}$  & $r$ & $E_j^{(\sigmasv)}$  & $r$ \\ \hline 
 1/20      & \num{2.1}{-2} &      & \num{2.2}{-2} &      & \num{5.4}{-3} &      & \num{2.5}{-2} &      & \num{4.6}{-2} &       \\ \hline
 1/40      & \num{6.1}{-3} &  3.4 & \num{3.9}{-3} &  5.6 & \num{1.3}{-3} &  4.0 & \num{6.7}{-3} &  3.8 & \num{9.3}{-3} &  5.0  \\ \hline
 1/80      & \num{1.8}{-3} &  3.3 & \num{7.4}{-4} &  5.2 & \num{3.2}{-4} &  4.2 & \num{1.9}{-3} &  3.6 & \num{2.2}{-3} &  4.1  \\ \hline
 1/160	   & \num{4.6}{-4} &  4.0 & \num{1.5}{-4} &  4.8 & \num{7.8}{-5} &  4.1 & \num{5.1}{-4} &  3.7 & \num{5.3}{-4} &  4.3  \\ \hline
\rateLabel &  $1.82$       &      &  $2.38$       &      &  $2.04$       &      &  $1.87$       &      &  $2.14$       &       \\ \hline
\end{tabular}
}

%
\bogus{
\begin{table}[hbt]\tableFont 
\begin{center}
\begin{tabular}{|l|c|c|c|c|c|c|c|c|c|c|c|} \hline 
grid  & N & $p$ & r & $\vert \mathbf{v} \vert$ & r & $\bar{\mathbf{u}}$ & r & $\bar{\mathbf{v}}$ & r & $\bar{\mathbf{s}}$& r \\ \hline 
                sq20 &     1 & \num{1.6}{-2} &      & \num{4.0}{-2} &      & \num{7.6}{-3} &      & \num{4.0}{-2} &      & \num{10.9}{1} &       \\ \hline
                sq40 &     2 & \num{2.1}{-2} &  0.7 & \num{8.0}{-3} &  5.0 & \num{1.5}{-3} &  5.2 & \num{8.1}{-3} &  4.9 & \num{15.5}{0} &  7.0  \\ \hline
                sq80 &     4 & \num{7.6}{-3} &  2.7 & \num{1.9}{-3} &  4.2 & \num{3.4}{-4} &  4.4 & \num{2.2}{-3} &  3.6 & \num{1.6}{0} &  9.5  \\ \hline
               sq160 &     8 & \num{1.5}{-3} &  5.1 & \num{4.6}{-4} &  4.2 & \num{8.0}{-5} &  4.2 & \num{5.8}{-4} &  3.8 & \num{1.9}{-1} &  8.5  \\ \hline
               sq320 &    16 & \num{2.7}{-4} &  5.6 & \num{1.1}{-4} &  4.1 & \num{1.9}{-5} &  4.1 & \num{1.5}{-4} &  3.9 & \num{5.3}{-2} &  3.7  \\ \hline
    rate             &       &  $1.55$       &      &  $2.11$       &      &  $2.14$       &      &  $2.00$       &      &  $2.83$       &       \\ \hline
\end{tabular}
\caption{Inses bulk: tz.trig.am1.rep3, solid=(swefos,fos), AM=1, amImp=0.5, vpf=0.1, ts=pc2, nc=1, normalMotion=1, bcs=periodic, bcsTop=Dirichlet, cdv=0.25, $t=.3$, $\mu=.05$, ad2=0, $\rho=1$, $\bar{\rho}=1e3$, scaleTZStrain=0, TZ=trig(fx=2,ft=2), cfl=0.5, Thu Jul 18 15:01:57 2013}\label{table:bulk.tz.trig.am1.rep3.pc2.swefos.fos.vector}
\end{center}
\end{table}
}
\newcommand{\tableMPVAIAmpHeavyTrigTZ}{%
\begin{tabular}{|c|c|c|c|c|c|c|c|c|c|c|} \hline
\multicolumn{11}{|c|}{MP-VA, trigonometric solution, heavy solid} \\ \hline 
\strutt~~$h_j$~~& $E_j^{(p)}$ & $r$ & $E_j^{(\vv)}$ & $r$ & $E_j^{(\usv)}$ & $r$  & $E_j^{(\vsv)}$  & $r$ & $E_j^{(\sigmasv)}$  & $r$ \\ \hline 
 1/20     & \num{1.6}{-2} &      & \num{4.0}{-2} &      & \num{7.6}{-3} &      & \num{4.0}{-2} &      & \num{1.1}{2} &       \\ \hline 
 1/40     & \num{2.1}{-2} &  0.7 & \num{8.0}{-3} &  5.0 & \num{1.5}{-3} &  5.2 & \num{8.1}{-3} &  4.9 & \num{1.6}{1} &  7.0  \\ \hline 
 1/80     & \num{7.6}{-3} &  2.7 & \num{1.9}{-3} &  4.2 & \num{3.4}{-4} &  4.4 & \num{2.2}{-3} &  3.6 & \num{1.6}{0} &  9.5  \\ \hline  
 1/160	  & \num{1.5}{-3} &  5.1 & \num{4.6}{-4} &  4.2 & \num{8.0}{-5} &  4.2 & \num{5.8}{-4} &  3.8 & \num{1.9}{-1} &  8.5  \\ \hline 
 1/320    & \num{2.7}{-4} &  5.6 & \num{1.1}{-4} &  4.1 & \num{1.9}{-5} &  4.1 & \num{1.5}{-4} &  3.9 & \num{5.3}{-2} &  3.7  \\ \hline
\rateLabel&  $1.55$       &      &  $2.11$       &      &  $2.14$       &      &  $2.00$       &      &  $2.83$       &       \\ \hline 
\end{tabular}
}

\begin{figure}[hbt]\tableFontSize
\begin{center}
\tableMPVAIAmpHeavyTrigTZ
\vskip.5\baselineskip
\tableMPVAIAmpMediumTrigTZ
\vskip.5\baselineskip
\tableMPVAIAmpLightTrigTZ
\caption{Trigonometric exact solution for a viscous incompressible fluid and acoustic solid (model problem MP-VA).
Maximum errors and estimated convergence rates at $t=0.3$ computed using the AMP algorithm for a heavy solid, $\delta=10^3$,
medium solid, $\delta=1$ and light solid, $\delta=10^{-1}$.
}
\label{tab:MP-VA_AMP_TZ_trig}
\end{center}
\end{figure}

%% file: texFiles/insSolidFig.tex
\newcommand{\drawTW}[9]{%
\begin{scope}[#1]
\draw(0,0) node[anchor=south west,xshift=-4pt,yshift=+0pt] {\trimfig{fig/#2}{\figWidth}};
\draw(0,5.7) node[draw,fill=white,anchor=north west,xshift=12pt,yshift=0pt] {\scriptsize #3};
\draw(0,0.0) node[draw,fill=white,anchor=south west,xshift=12pt,yshift=8pt] {\scriptsize #4};
\draw(3.5,.0) node[draw,fill=white,anchor=south east,xshift=-12pt,yshift=8pt] {\scriptsize #5};
\draw (\xcbs,\ycbs) node[anchor=south west,xshift= +0pt,yshift=+0pt] {\trimfigcbs{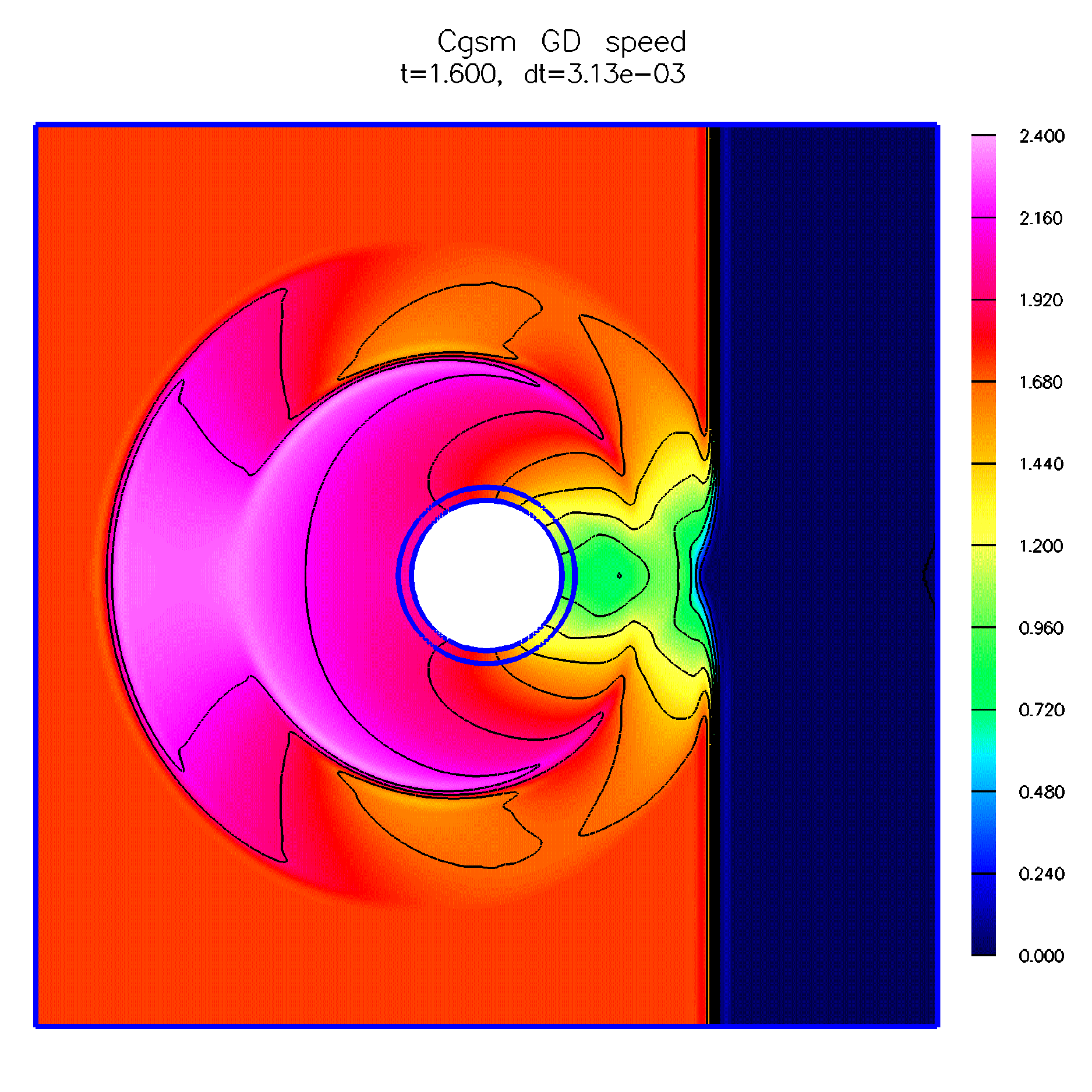}{\cbsWidth}{\cbsHeight}};
\draw (\xcbs,\ycbs) node[anchor=south west,xshift= +8pt,yshift=+1pt] {\scriptsize $#6$};
\draw (\xcbs,\ycbsTop) node[anchor=south west,xshift= +8pt,yshift=-6pt] {\scriptsize $#7$};
\draw (\xcb,\ycb) node[anchor=south west,xshift= +0pt,yshift=+0pt] {\trimfigcb{fig/colourBarLines}{\cbWidth}{\cbHeight}};
\draw (\xcb,\ycb) node[anchor=south west,xshift= +8pt,yshift=+1pt] {\scriptsize $#8$};
\draw (\xcb,\ycbTop) node[anchor=south west,xshift= +8pt,yshift=-6pt] {\scriptsize $#9$};
\end{scope}
}
{
\newcommand{\cbsWidth}{.2cm}
\newcommand{\cbsHeight}{1.9cm}
\newcommand{\xcbs}{4.2cm}
\newcommand{\ycbs}{3.75cm}
\newlength{\ycbsTop}%
\setlength{\ycbsTop}{\ycbs+\cbsHeight}
\newlength{\ycbsMid}%
\setlength{\ycbsMid}{\ycbs+\cbsHeight*\real{.5}}
\newcommand{\xsLabel}{6.5cm}
\newcommand{\ysLabel}{6.5cm}
\newcommand{\trimfigcbs}[3]{\includegraphics[width=#2, height=#3, clip, trim=17cm 2.35cm 1.65cm 2.35cm]{#1}}
\newcommand{\cbWidth}{.2cm}
\newcommand{\cbHeight}{3.6cm}
\newcommand{\xcb}{4.2cm}
\newcommand{\ycb}{0.1cm}
\newlength{\ycbTop}%
\setlength{\ycbTop}{\ycb+\cbHeight}
\newlength{\ycbMid}%
\setlength{\ycbMid}{\ycb+\cbHeight*\real{.5}}
\newcommand{\xLabel}{6.5cm}
\newcommand{\yLabel}{6.5cm}
\newcommand{\trimfigcb}[3]{\includegraphics[width=#2, height=#3, clip, trim=17cm 2.35cm 1.65cm 2.35cm]{#1}}
\newcommand{\figWidth}{4.25cm}
\newcommand{\trimfig}[2]{\trimFig{#1}{#2}{.45}{.45}{.1}{.425}}
\newcommand{\xa}{.08}
\newcommand{\ysa}{6.5}
\newcommand{\ysi}{10.05}
\newcommand{\ysb}{12.05}
\newcommand{\yAxis}{%
    \draw[-,black,thick,xshift=2pt,yshift=-1.7pt] (-.08,\ysa) -- (.08,\ysa);
    \draw[-,black,thick,xshift=2pt,yshift=1.75pt] (-.08,\ysi) -- (.08,\ysi);
    \draw[-,black,thick,xshift=2pt,yshift=-2.4pt] (-.08,\ysb) -- (.08,\ysb);
    \draw[-,black,thick,xshift=2pt,yshift=-2pt] (0,\ysa) -- (0,\ysb);
    \draw (0,\ysa) node[anchor=east,xshift=2pt] {\scriptsize $-1$};
    \draw (0,\ysi) node[anchor=east,xshift=2pt,yshift=1pt] {\scriptsize $0$};
    \draw (0,\ysb) node[anchor=east,xshift=2pt,yshift=-4pt] {\scriptsize $.5$};
    \draw (0,11) node[anchor=east,xshift=2pt,yshift=-4pt] {\scriptsize $y$};
}
\begin{figure}[htb]
\begin{center}
\begin{tikzpicture}[scale=1]
  \useasboundingbox (-.5,.5) rectangle (16.,12.);  
%
%
  \drawTW{xshift= 0.0cm,yshift=6.2cm}{insSolid_mu02_scfp05_EWE_FOS_t0p5_v1}{$\vs_1$}{$v_1$}{$\rhos/\rho=.05$}{-.05}{.05}{-.05}{.05};
  \drawTW{xshift= 5.3cm,yshift=6.2cm}{insSolid_mu02_scfp05_EWE_FOS_t0p5_v2}{$\vs_2$}{$v_2$}{$\rhos/\rho=.05$}{-.10}{.10}{-.10}{.10};
  \drawTW{xshift=10.6cm,yshift=6.2cm}{insSolid_mu02_scfp05_EWE_FOS_t0p5_p_s22}{$\sigmas_{22}$}{$p$}{$\rhos/\rho=.05$}{-.036}{.036}{-.041}{.041};
%
%
  \drawTW{xshift= 0.0cm,yshift= 0.0cm}{insSolid_mu02_scf1e2_EWE_FOS_t0p5_v1}{$\vs_1$}{$v_1$}{$\rhos/\rho=100$}{-.35}{.35}{-.35}{.35};
  \drawTW{xshift= 5.3cm,yshift= 0.0cm}{insSolid_mu02_scf1e2_EWE_FOS_t0p5_v2}{$\vs_2$}{$v_2$}{$\rhos/\rho=100$}{-.55}{.55}{-.55}{.55};
  \drawTW{xshift=10.6cm,yshift= 0.0cm}{insSolid_mu02_scf1e2_EWE_FOS_t0p5_p_s22}{$\sigmas_{22}$}{$p$}{$\rhos/\rho=100$}{-62.}{62.}{-.67}{.67};
%
  \begin{scope}[xshift=\xa cm,yshift=-6.21cm]
    \yAxis
  \end{scope}
   \draw[-,black,thick,xshift=3.5pt,yshift=2pt] (.15,0) -- (3.8,0);
   \draw  (.15,0) node[-,black,anchor=north west, xshift=0pt,yshift=2pt]  {\scriptsize $0$};
   \draw  (3.8,0) node[-,black,anchor=north east, xshift=6pt,yshift=2pt]  {\scriptsize $1$};
   \draw  (2.0,0) node[-,black,anchor=north east, xshift=8pt,yshift=2pt]  {\scriptsize $x$};
   \draw[-,black,thick,xshift=3.5pt,yshift=2.pt]  (.15,-.08) --  (.15,0.08); 
   \draw[-,black,thick,xshift=3.5pt,yshift=2.pt]  (3.8,-.08) --  (3.8,0.08); 
%
%
\end{tikzpicture}
\end{center}
  \caption{Traveling wave solution for a viscous incompressible fluid and elastic solid, model problem MP-VE, at $t=0.5$. Top row: shaded contours of the computed solution for
  $\delta=\rhos/\rho=.05$. Bottom row: shaded contours for $\delta=\rhos/\rho=100$.  The deformed interface is shown
   super-imposed as a curve.  
      }
  \label{fig:TW_EWE}
\end{figure}
}

%% file: texFiles/twMPIA1ErrorTables.tex
\bogus{
\begin{table}[hbt]\tableFont 
\begin{center}
\begin{tabular}{|l|c|c|c|c|c|c|c|c|c|c|c|} \hline 
grid  & N & $p$ & r & $\vert \mathbf{v} \vert$ & r & $\bar{\mathbf{u}}$ & r & $\bar{\mathbf{v}}$ & r & $\bar{\mathbf{s}}$& r \\ \hline 
                sq20 &     1 & \num{8.7}{0} &      & \num{2.2}{0} &      & \num{6.3}{-2} &      & \num{1.8}{0} &      & \num{29.1}{2} &       \\ \hline
                sq40 &     2 & \num{3.0}{0} &  2.9 & \num{7.5}{-1} &  2.9 & \num{2.2}{-2} &  2.9 & \num{6.5}{-1} &  2.8 & \num{97.0}{1} &  3.0  \\ \hline
                sq80 &     4 & \num{7.0}{-1} &  4.3 & \num{1.7}{-1} &  4.4 & \num{5.6}{-3} &  4.0 & \num{1.6}{-1} &  4.0 & \num{24.4}{1} &  4.0  \\ \hline
               sq160 &     8 & \num{1.6}{-1} &  4.4 & \num{3.8}{-2} &  4.5 & \num{1.4}{-3} &  4.1 & \num{3.8}{-2} &  4.2 & \num{60.0}{0} &  4.1  \\ \hline
    rate             &       &  $1.94$       &      &  $1.96$       &      &  $1.86$       &      &  $1.87$       &      &  $1.88$       &       \\ \hline
\end{tabular}
\caption{Inses bulk: fw.mu0.rep3, solid=(swefos,fos), AM=1, amImp=0.5, vpf=0.1, ts=pc2, nc=1, normalMotion=1, bcs=periodic, bcsTop=displacement, cdv=0.25, $t=1.$, $\mu=0.$, ad2=1, $\rho=1$, $\bar{\rho}=1e3$, scaleTZStrain=0, , cfl=0.5, Wed Jul 17 17:48:11 2013}\label{table:bulk.fw.mu0.rep3.pc2.swefos.fos.vector}
\end{center}
\end{table}
}
\newcommand{\tableMPIAIAmpHeavy}{%
\begin{tabular}{|c|c|c|c|c|c|c|c|c|c|c|} \hline
\multicolumn{11}{|c|}{MP-IA, traveling wave, heavy solid} \\ \hline 
\strutt~~$h_j$~~& $E_j^{(p)}$ & $r$ & $E_j^{(\vv)}$ & $r$ & $E_j^{(\usv)}$ & $r$  & $E_j^{(\vsv)}$  & $r$ & $E_j^{(\sigmasv)}$  & $r$ \\ \hline 1/20      & \num{8.7}{0} &      & \num{2.2}{0} &      & \num{6.3}{-2} &      & \num{1.8}{0} &      & \num{2.9}{3} &       \\ \hline	  
 1/40      & \num{3.0}{0} &  2.9 & \num{7.5}{-1} &  2.9 & \num{2.2}{-2} &  2.9 & \num{6.5}{-1} &  2.8 & \num{9.7}{2} &  3.0  \\ \hline 
 1/80      & \num{7.0}{-1} &  4.3 & \num{1.7}{-1} &  4.4 & \num{5.6}{-3} &  4.0 & \num{1.6}{-1} &  4.0 & \num{2.4}{2} &  4.0  \\ \hline
 1/160	   & \num{1.6}{-1} &  4.4 & \num{3.8}{-2} &  4.5 & \num{1.4}{-3} &  4.1 & \num{3.8}{-2} &  4.2 & \num{6.0}{1} &  4.1  \\ \hline
\rateLabel &  $1.94$       &      &  $1.96$       &      &  $1.86$       &      &  $1.87$       &      &  $1.88$       &       \\ \hline
\end{tabular}
}
\bogus{
\begin{table}[hbt]\tableFont 
\begin{center}
\begin{tabular}{|l|c|c|c|c|c|c|c|c|c|c|c|} \hline 
grid  & N & $p$ & r & $\vert \mathbf{v} \vert$ & r & $\bar{\mathbf{u}}$ & r & $\bar{\mathbf{v}}$ & r & $\bar{\mathbf{s}}$& r \\ \hline 
                sq20 &     1 & \num{3.3}{-1} &      & \num{1.6}{-1} &      & \num{1.7}{-2} &      & \num{2.6}{-1} &      & \num{3.5}{-1} &       \\ \hline
                sq40 &     2 & \num{8.9}{-2} &  3.7 & \num{4.5}{-2} &  3.6 & \num{4.9}{-3} &  3.6 & \num{7.1}{-2} &  3.8 & \num{1.0}{-1} &  3.4  \\ \hline
                sq80 &     4 & \num{2.2}{-2} &  4.0 & \num{1.1}{-2} &  3.9 & \num{1.2}{-3} &  3.9 & \num{1.8}{-2} &  3.9 & \num{3.1}{-2} &  3.4  \\ \hline
               sq160 &     8 & \num{5.3}{-3} &  4.2 & \num{2.9}{-3} &  4.0 & \num{3.1}{-4} &  4.0 & \num{4.8}{-3} &  3.8 & \num{8.8}{-3} &  3.5  \\ \hline
    rate             &       &  $1.99$       &      &  $1.94$       &      &  $1.94$       &      &  $1.93$       &      &  $1.77$       &       \\ \hline
\end{tabular}
\caption{Inses bulk: fw.mu0.rep0, solid=(swefos,fos), AM=1, amImp=0.5, vpf=0.1, ts=pc2, nc=1, normalMotion=1, bcs=periodic, bcsTop=displacement, cdv=0.25, $t=1.$, $\mu=0.$, ad2=1, $\rho=1$, $\bar{\rho}=1.$, scaleTZStrain=0, , cfl=0.5, Wed Jul 17 17:36:44 2013}\label{table:bulk.fw.mu0.rep0.pc2.swefos.fos.vector}
\end{center}
\end{table}
}
\newcommand{\tableMPIAIAmpMedium}{%
\begin{tabular}{|c|c|c|c|c|c|c|c|c|c|c|} \hline
\multicolumn{11}{|c|}{MP-IA, traveling wave, medium solid} \\ \hline 
\strutt~~$h_j$~~& $E_j^{(p)}$ & $r$ & $E_j^{(\vv)}$ & $r$ & $E_j^{(\usv)}$ & $r$  & $E_j^{(\vsv)}$  & $r$ & $E_j^{(\sigmasv)}$  & $r$ \\ \hline
 1/20      & \num{3.3}{-1} &      & \num{1.6}{-1} &      & \num{1.7}{-2} &      & \num{2.6}{-1} &      & \num{3.5}{-1} &       \\ \hline
 1/40      & \num{8.9}{-2} &  3.7 & \num{4.5}{-2} &  3.6 & \num{4.9}{-3} &  3.6 & \num{7.1}{-2} &  3.8 & \num{1.0}{-1} &  3.4  \\ \hline
 1/80      & \num{2.2}{-2} &  4.0 & \num{1.1}{-2} &  3.9 & \num{1.2}{-3} &  3.9 & \num{1.8}{-2} &  3.9 & \num{3.1}{-2} &  3.4  \\ \hline
 1/160	   & \num{5.3}{-3} &  4.2 & \num{2.9}{-3} &  4.0 & \num{3.1}{-4} &  4.0 & \num{4.8}{-3} &  3.8 & \num{8.8}{-3} &  3.5  \\ \hline
\rateLabel &  $1.99$       &      &  $1.94$       &      &  $1.94$       &      &  $1.93$       &      &  $1.77$       &       \\ \hline
\end{tabular}
}

\bogus{
\begin{table}[hbt]\tableFont 
\begin{center}
\begin{tabular}{|l|c|c|c|c|c|c|c|c|c|c|c|} \hline 
grid  & N & $p$ & r & $\vert \mathbf{v} \vert$ & r & $\bar{\mathbf{u}}$ & r & $\bar{\mathbf{v}}$ & r & $\bar{\mathbf{s}}$& r \\ \hline 
                sq20 &     1 & \num{2.5}{-1} &      & \num{1.4}{-1} &      & \num{1.2}{-1} &      & \num{1.9}{0} &      & \num{2.5}{-1} &       \\ \hline
                sq40 &     2 & \num{5.9}{-2} &  4.2 & \num{3.3}{-2} &  4.2 & \num{2.9}{-2} &  4.2 & \num{4.5}{-1} &  4.2 & \num{5.9}{-2} &  4.2  \\ \hline
                sq80 &     4 & \num{1.5}{-2} &  4.0 & \num{8.2}{-3} &  4.1 & \num{7.6}{-3} &  3.9 & \num{1.2}{-1} &  3.9 & \num{1.5}{-2} &  4.0  \\ \hline
               sq160 &     8 & \num{3.7}{-3} &  4.0 & \num{2.1}{-3} &  3.9 & \num{1.9}{-3} &  4.0 & \num{3.0}{-2} &  3.9 & \num{3.7}{-3} &  4.0  \\ \hline
    rate             &       &  $2.02$       &      &  $2.01$       &      &  $2.00$       &      &  $1.99$       &      &  $2.02$       &       \\ \hline
\end{tabular}
\caption{Inses bulk: fw.mu0.rem1, solid=(swefos,fos), AM=1, amImp=0.5, vpf=0.1, ts=pc2, nc=1, normalMotion=1, bcs=periodic, bcsTop=displacement, cdv=0.25, $t=1.$, $\mu=0.$, ad2=1, $\rho=1$, $\bar{\rho}=.1$, scaleTZStrain=0, , cfl=0.5, Wed Jul 17 17:32:13 2013}\label{table:bulk.fw.mu0.rem1.pc2.swefos.fos.vector}
\end{center}
\end{table}
}
\newcommand{\tableMPIAIAmpLight}{%
\begin{tabular}{|c|c|c|c|c|c|c|c|c|c|c|} \hline
\multicolumn{11}{|c|}{MP-IA, traveling wave, light solid} \\ \hline 
\strutt~~$h_j$~~& $E_j^{(p)}$ & $r$ & $E_j^{(\vv)}$ & $r$ & $E_j^{(\usv)}$ & $r$  & $E_j^{(\vsv)}$  & $r$ & $E_j^{(\sigmasv)}$  & $r$ \\ \hline
 1/20      & \num{2.5}{-1} &      & \num{1.4}{-1} &      & \num{1.2}{-1} &      & \num{1.9}{0} &      & \num{2.5}{-1} &       \\ \hline     
 1/40      & \num{5.9}{-2} &  4.2 & \num{3.3}{-2} &  4.2 & \num{2.9}{-2} &  4.2 & \num{4.5}{-1} &  4.2 & \num{5.9}{-2} &  4.2  \\ \hline    
 1/80      & \num{1.5}{-2} &  4.0 & \num{8.2}{-3} &  4.1 & \num{7.6}{-3} &  3.9 & \num{1.2}{-1} &  3.9 & \num{1.5}{-2} &  4.0  \\ \hline    
 1/160	   & \num{3.7}{-3} &  4.0 & \num{2.1}{-3} &  3.9 & \num{1.9}{-3} &  4.0 & \num{3.0}{-2} &  3.9 & \num{3.7}{-3} &  4.0  \\ \hline    
\rateLabel &  $2.02$       &      &  $2.01$       &      &  $2.00$       &      &  $1.99$       &      &  $2.02$       &       \\ \hline  
\end{tabular}
}

\begin{figure}[hbt]\tableFontSize
\begin{center}
\tableMPIAIAmpHeavy
\vskip.5\baselineskip
\tableMPIAIAmpMedium
\vskip.5\baselineskip
\tableMPIAIAmpLight
\caption{Traveling wave solution for an inviscid incompressible fluid and acoustic solid (model problem MP-IA).  Maximum errors and estimated convergence rates at $t=1.0$, computed using the AMP scheme for a heavy solid, $\delta=10^3$, medium solid, $\delta=1$ and light solid, $\delta=10^{-1}$.
}
\label{tab:MP-IA1_AMP}
\end{center}
\end{figure}

%% file: texFiles/twMPVA1ErrorTables.tex
\bogus{
\begin{table}[hbt]\tableFont 
\begin{center}
\begin{tabular}{|l|c|c|c|c|c|c|c|c|c|c|c|} \hline 
grid  & N & $p$ & r & $\vert \mathbf{v} \vert$ & r & $\bar{\mathbf{u}}$ & r & $\bar{\mathbf{v}}$ & r & $\bar{\mathbf{s}}$& r \\ \hline 
                sq20 &     1 & \num{2.9}{-1} &      & \num{1.5}{-1} &      & \num{2.1}{-3} &      & \num{2.1}{-2} &      & \num{33.3}{0} &       \\ \hline
                sq40 &     2 & \num{8.5}{-2} &  3.3 & \num{4.1}{-2} &  3.6 & \num{5.0}{-4} &  4.3 & \num{4.8}{-3} &  4.5 & \num{7.5}{0} &  4.4  \\ \hline
                sq80 &     4 & \num{1.5}{-2} &  5.8 & \num{7.7}{-3} &  5.3 & \num{1.2}{-4} &  4.1 & \num{1.1}{-3} &  4.2 & \num{1.8}{0} &  4.2  \\ \hline
               sq160 &     8 & \num{2.0}{-3} &  7.2 & \num{1.4}{-3} &  5.4 & \num{3.9}{-5} &  3.1 & \num{3.7}{-4} &  3.1 & \num{5.7}{-1} &  3.2  \\ \hline
    rate             &       &  $2.40$       &      &  $2.25$       &      &  $1.94$       &      &  $1.96$       &      &  $1.96$       &       \\ \hline
\end{tabular}
\caption{Inses bulk: fw.mu02.rep3, solid=(swefos,fos), AM=1, amImp=0.5, vpf=0.1, ts=pc2, nc=1, normalMotion=1, bcs=periodic, bcsTop=displacement, cdv=0.25, $t=.3$, $\mu=.02$, ad2=0, $\rho=1$, $\bar{\rho}=1e3$, , cfl=0.5, Tue Jul 16  9:30:23 2013}\label{table:bulk.fw.mu02.rep3.pc2.swefos.fos.vector}
\end{center}
\end{table}
}
\newcommand{\tableMPVAIAmpHeavy}{%
\begin{tabular}{|c|c|c|c|c|c|c|c|c|c|c|} \hline
\multicolumn{11}{|c|}{MP-VA, traveling wave, heavy solid} \\ \hline 
\strutt~~$h_j$~~& $E_j^{(p)}$ & $r$ & $E_j^{(\vv)}$ & $r$ & $E_j^{(\usv)}$ & $r$  & $E_j^{(\vsv)}$  & $r$ & $E_j^{(\sigmasv)}$  & $r$ \\ \hline
 1/20      & \num{2.9}{-1} &      & \num{1.5}{-1} &      & \num{2.1}{-3} &      & \num{2.1}{-2} &      & \num{3.3}{1} &       \\ \hline
 1/40      & \num{8.5}{-2} &  3.3 & \num{4.1}{-2} &  3.6 & \num{5.0}{-4} &  4.3 & \num{4.8}{-3} &  4.5 & \num{7.5}{0} &  4.4  \\ \hline 
 1/80      & \num{1.5}{-2} &  5.8 & \num{7.7}{-3} &  5.3 & \num{1.2}{-4} &  4.1 & \num{1.1}{-3} &  4.2 & \num{1.8}{0} &  4.2  \\ \hline 
 1/160	   & \num{2.0}{-3} &  7.2 & \num{1.4}{-3} &  5.4 & \num{3.9}{-5} &  3.1 & \num{3.7}{-4} &  3.1 & \num{5.7}{-1} &  3.2  \\ \hline
\rateLabel &  $2.40$       &      &  $2.25$       &      &  $1.94$       &      &  $1.96$       &      &  $1.96$       &       \\ \hline
\end{tabular}
}
\bogus{
\begin{table}[hbt]\tableFont 
\begin{center}
\begin{tabular}{|l|c|c|c|c|c|c|c|c|c|c|c|} \hline 
grid  & N & $p$ & r & $\vert \mathbf{v} \vert$ & r & $\bar{\mathbf{u}}$ & r & $\bar{\mathbf{v}}$ & r & $\bar{\mathbf{s}}$& r \\ \hline 
                sq20 &     1 & \num{5.4}{-2} &      & \num{4.8}{-2} &      & \num{3.2}{-3} &      & \num{3.6}{-2} &      & \num{7.7}{-2} &       \\ \hline
                sq40 &     2 & \num{1.1}{-2} &  4.8 & \num{1.1}{-2} &  4.4 & \num{8.4}{-4} &  3.9 & \num{9.3}{-3} &  3.9 & \num{1.8}{-2} &  4.4  \\ \hline
                sq80 &     4 & \num{2.3}{-3} &  4.8 & \num{2.3}{-3} &  4.8 & \num{2.0}{-4} &  4.2 & \num{2.2}{-3} &  4.3 & \num{4.0}{-3} &  4.4  \\ \hline
               sq160 &     8 & \num{5.5}{-4} &  4.2 & \num{5.2}{-4} &  4.3 & \num{4.9}{-5} &  4.0 & \num{5.2}{-4} &  4.1 & \num{1.1}{-3} &  3.7  \\ \hline
    rate             &       &  $2.21$       &      &  $2.18$       &      &  $2.02$       &      &  $2.05$       &      &  $2.07$       &       \\ \hline
\end{tabular}
\caption{Inses bulk: fw.mu02.rep0, solid=(swefos,fos), AM=1, amImp=0.5, vpf=0.1, ts=pc2, nc=1, normalMotion=1, bcs=periodic, bcsTop=displacement, cdv=0.25, $t=.3$, $\mu=.02$, ad2=0, $\rho=1$, $\bar{\rho}=1.$, , cfl=0.5, Tue Jul 16  9:24:39 2013}\label{table:bulk.fw.mu02.rep0.pc2.swefos.fos.vector}
\end{center}
\end{table}
}
\newcommand{\tableMPVAIAmpMedium}{%
\begin{tabular}{|c|c|c|c|c|c|c|c|c|c|c|} \hline
\multicolumn{11}{|c|}{MP-VA, traveling wave, medium solid} \\ \hline 
\strutt~~$h_j$~~& $E_j^{(p)}$ & $r$ & $E_j^{(\vv)}$ & $r$ & $E_j^{(\usv)}$ & $r$  & $E_j^{(\vsv)}$  & $r$ & $E_j^{(\sigmasv)}$  & $r$ \\ \hline 1/20      & \num{5.4}{-2} &      & \num{4.8}{-2} &      & \num{3.2}{-3} &      & \num{3.6}{-2} &      & \num{7.7}{-2} &       \\ \hline
 1/40      & \num{1.1}{-2} &  4.8 & \num{1.1}{-2} &  4.4 & \num{8.4}{-4} &  3.9 & \num{9.3}{-3} &  3.9 & \num{1.8}{-2} &  4.4  \\ \hline
 1/80      & \num{2.3}{-3} &  4.8 & \num{2.3}{-3} &  4.8 & \num{2.0}{-4} &  4.2 & \num{2.2}{-3} &  4.3 & \num{4.0}{-3} &  4.4  \\ \hline
 1/160	   & \num{5.5}{-4} &  4.2 & \num{5.2}{-4} &  4.3 & \num{4.9}{-5} &  4.0 & \num{5.2}{-4} &  4.1 & \num{1.1}{-3} &  3.7  \\ \hline
\rateLabel &  $2.21$       &      &  $2.18$       &      &  $2.02$       &      &  $2.05$       &      &  $2.07$       &       \\ \hline
\end{tabular}
}

\bogus{
\begin{table}[hbt]\tableFont 
\begin{center}
\begin{tabular}{|l|c|c|c|c|c|c|c|c|c|c|c|} \hline 
grid  & N & $p$ & r & $\vert \mathbf{v} \vert$ & r & $\bar{\mathbf{u}}$ & r & $\bar{\mathbf{v}}$ & r & $\bar{\mathbf{s}}$& r \\ \hline 
                sq20 &     1 & \num{2.5}{-3} &      & \num{8.3}{-3} &      & \num{2.1}{-3} &      & \num{8.3}{-3} &      & \num{6.2}{-3} &       \\ \hline
                sq40 &     2 & \num{6.2}{-4} &  4.0 & \num{1.9}{-3} &  4.3 & \num{5.3}{-4} &  3.9 & \num{1.9}{-3} &  4.3 & \num{1.3}{-3} &  4.6  \\ \hline
                sq80 &     4 & \num{1.6}{-4} &  4.0 & \num{4.2}{-4} &  4.5 & \num{1.3}{-4} &  4.2 & \num{4.2}{-4} &  4.5 & \num{2.9}{-4} &  4.6  \\ \hline
               sq160 &     8 & \num{4.0}{-5} &  3.9 & \num{9.8}{-5} &  4.3 & \num{3.0}{-5} &  4.2 & \num{9.8}{-5} &  4.3 & \num{6.9}{-5} &  4.2  \\ \hline
    rate             &       &  $1.99$       &      &  $2.14$       &      &  $2.04$       &      &  $2.14$       &      &  $2.16$       &       \\ \hline
\end{tabular}
\caption{Inses bulk: fw.mu02.rem1, solid=(swefos,fos), AM=1, amImp=0.5, vpf=0.1, ts=pc2, nc=1, normalMotion=1, bcs=periodic, bcsTop=displacement, cdv=0.25, $t=.3$, $\mu=.02$, ad2=0, $\rho=1$, $\bar{\rho}=1e-1$, scaleTZStrain=0, , cfl=0.5, Wed Jul 17 18:00:27 2013}\label{table:bulk.fw.mu02.rem1.pc2.swefos.fos.vector}
\end{center}
\end{table}
}
\newcommand{\tableMPVAIAmpLight}{%
\begin{tabular}{|c|c|c|c|c|c|c|c|c|c|c|} \hline
\multicolumn{11}{|c|}{MP-VA, traveling wave, light solid} \\ \hline 
\strutt~~$h_j$~~& $E_j^{(p)}$ & $r$ & $E_j^{(\vv)}$ & $r$ & $E_j^{(\usv)}$ & $r$  & $E_j^{(\vsv)}$  & $r$ & $E_j^{(\sigmasv)}$  & $r$ \\ \hline 1/20      & \num{2.5}{-3} &      & \num{8.3}{-3} &      & \num{2.1}{-3} &      & \num{8.3}{-3} &      & \num{6.2}{-3} &       \\ \hline
 1/40      & \num{6.2}{-4} &  4.0 & \num{1.9}{-3} &  4.3 & \num{5.3}{-4} &  3.9 & \num{1.9}{-3} &  4.3 & \num{1.3}{-3} &  4.6  \\ \hline
 1/80      & \num{1.6}{-4} &  4.0 & \num{4.2}{-4} &  4.5 & \num{1.3}{-4} &  4.2 & \num{4.2}{-4} &  4.5 & \num{2.9}{-4} &  4.6  \\ \hline
 1/160	   & \num{4.0}{-5} &  3.9 & \num{9.8}{-5} &  4.3 & \num{3.0}{-5} &  4.2 & \num{9.8}{-5} &  4.3 & \num{6.9}{-5} &  4.2  \\ \hline
\rateLabel &  $1.99$       &      &  $2.14$       &      &  $2.04$       &      &  $2.14$       &      &  $2.16$       &       \\ \hline
\end{tabular}
}

\begin{figure}[hbt]\tableFontSize
\begin{center}
\tableMPVAIAmpHeavy
\vskip.5\baselineskip
\tableMPVAIAmpMedium
\vskip.5\baselineskip
\tableMPVAIAmpLight
\caption{Traveling wave solution for an viscous incompressible fluid and acoustic solid (model problem MP-VA).  Maximum errors and estimated convergence rates at $t=0.3$, computed using the AMP scheme for a heavy solid, $\delta=10^3$, medium solid, $\delta=1$ and light solid, $\delta=10^{-1}$. 
}
\label{tab:MP-VA1_AMP}
\end{center}
\end{figure}

%% file: texFiles/twMPVEErrorTables.tex
\bogus{
\begin{table}[hbt]\tableFont 
\begin{center}
\begin{tabular}{|l|c|c|c|c|c|c|c|c|c|c|c|} \hline 
grid  & N & $p$ & r & $\vert \mathbf{v} \vert$ & r & $\bar{\mathbf{u}}$ & r & $\bar{\mathbf{v}}$ & r & $\bar{\mathbf{s}}$& r \\ \hline 
sq20 &     1 & \num{1.2}{-2} &      & \num{1.9}{-2} &      & \num{2.4}{-3} &      & \num{1.6}{-2} &      & \num{3.5}{1} &       \\ \hline
sq40 &     2 & \num{2.9}{-3} &  4.1 & \num{3.7}{-3} &  5.1 & \num{4.5}{-4} &  5.4 & \num{3.1}{-3} &  5.2 & \num{9.1}{0} &  3.9  \\ \hline
sq80 &     4 & \num{6.5}{-4} &  4.5 & \num{6.0}{-4} &  6.1 & \num{8.3}{-5} &  5.4 & \num{6.0}{-4} &  5.1 & \num{2.5}{0} &  3.6  \\ \hline
sq160 &    8 & \num{1.5}{-4} &  4.3 & \num{1.3}{-4} &  4.8 & \num{1.6}{-5} &  5.0 & \num{1.2}{-4} &  4.8 & \num{6.8}{-1} &  3.7  \\ \hline  
 rate &      &  $2.11$       &      &  $2.42$       &      &  $2.40$       &      &  $2.34$       &      &  $1.89$       &       \\ \hline
\end{tabular}
\caption{Inses bulk: fw.mu02.rep3, solid=(ewe,fos), AM=1, amImp=0.5, vpf=0.1, ts=pc2, nc=1, normalMotion=0, bcs=periodic, bcsTop=displacement, cdv=0.25, $t=.3$, $\mu=.02$, ad2=0, $\rho=1$, $\bar{\rho}=1e3$, , cfl=0.5, Tue Jul 16  4:34:42 2013}\label{table:bulk.fw.mu02.rep3.pc2.ewe.fos.vector}
\end{center}
\end{table}
}
\newcommand{\tableMPVEAmpHeavy}{%
\begin{tabular}{|c|c|c|c|c|c|c|c|c|c|c|} \hline
\multicolumn{11}{|c|}{MP-VE, traveling wave, heavy solid} \\ \hline 
\strutt~~$h_j$~~& $E_j^{(p)}$ & $r$ & $E_j^{(\vv)}$ & $r$ & $E_j^{(\usv)}$ & $r$  & $E_j^{(\vsv)}$  & $r$ & $E_j^{(\sigmasv)}$  & $r$ \\ \hline
 1/20   & \num{1.2}{-2} &      & \num{1.9}{-2} &      & \num{2.4}{-3} &      & \num{1.6}{-2} &      & \num{3.5}{1} &       \\ \hline   
 1/40   & \num{2.9}{-3} &  4.1 & \num{3.7}{-3} &  5.1 & \num{4.5}{-4} &  5.4 & \num{3.1}{-3} &  5.2 & \num{9.1}{0} &  3.9  \\ \hline    
 1/80   & \num{6.5}{-4} &  4.5 & \num{6.0}{-4} &  6.1 & \num{8.3}{-5} &  5.4 & \num{6.0}{-4} &  5.1 & \num{2.5}{0} &  3.6  \\ \hline    
 1/160	& \num{1.5}{-4} &  4.3 & \num{1.3}{-4} &  4.8 & \num{1.6}{-5} &  5.0 & \num{1.2}{-4} &  4.8 & \num{6.8}{-1} &  3.7  \\ \hline   
\rateLabel &  $2.18$       &      &  $2.21$       &      &  $2.03$       &      &  $2.03$       &      &  $1.95$       &       \\ \hline   
\end{tabular}
}

\bogus{
\begin{table}[hbt]\tableFont 
\begin{center}
\begin{tabular}{|l|c|c|c|c|c|c|c|c|c|c|c|} \hline 
grid  & N & $p$ & r & $\vert \mathbf{v} \vert$ & r & $\bar{\mathbf{u}}$ & r & $\bar{\mathbf{v}}$ & r & $\bar{\mathbf{s}}$& r \\ \hline 
                sq20 &     1 & \num{1.3}{-2} &      & \num{1.2}{-2} &      & \num{2.6}{-3} &      & \num{9.2}{-3} &      & \num{4.4}{-2} &       \\ \hline
                sq40 &     2 & \num{2.6}{-3} &  5.0 & \num{2.7}{-3} &  4.3 & \num{5.9}{-4} &  4.4 & \num{2.2}{-3} &  4.2 & \num{8.5}{-3} &  5.2  \\ \hline
                sq80 &     4 & \num{4.8}{-4} &  5.5 & \num{5.6}{-4} &  4.8 & \num{1.3}{-4} &  4.6 & \num{5.0}{-4} &  4.3 & \num{1.7}{-3} &  5.0  \\ \hline
               sq160 &     8 & \num{9.4}{-5} &  5.1 & \num{1.2}{-4} &  4.5 & \num{2.9}{-5} &  4.5 & \num{1.2}{-4} &  4.3 & \num{4.1}{-4} &  4.2  \\ \hline
    rate             &       &  $2.39$       &      &  $2.19$       &      &  $2.17$       &      &  $2.10$       &      &  $2.26$       &       \\ \hline  
\end{tabular}
\caption{Inses bulk: fw.mu02.rep0, solid=(ewe,fos), AM=1, amImp=0.5, vpf=0.1, ts=pc2, nc=1, normalMotion=0, bcs=periodic, bcsTop=displacement, cdv=0.25, $t=.3$, $\mu=.02$, ad2=0, $\rho=1$, $\bar{\rho}=1.$, , cfl=0.5, Tue Jul 16  5:01:57 2013}\label{table:bulk.fw.mu02.rep0.pc2.ewe.fos.vector}
\end{center}
\end{table}
}
\newcommand{\tableMPVEAmpMedium}{%
\begin{tabular}{|c|c|c|c|c|c|c|c|c|c|c|} \hline
\multicolumn{11}{|c|}{MP-VE, traveling wave, medium solid} \\ \hline 
\strutt~~$h_j$~~& $E_j^{(p)}$ & $r$ & $E_j^{(\vv)}$ & $r$ & $E_j^{(\usv)}$ & $r$  & $E_j^{(\vsv)}$  & $r$ & $E_j^{(\sigmasv)}$  & $r$ \\ \hline
 1/20      & \num{1.3}{-2} &      & \num{1.2}{-2} &      & \num{2.6}{-3} &      & \num{9.2}{-3} &      & \num{4.4}{-2} &       \\ \hline     
 1/40      & \num{2.6}{-3} &  5.0 & \num{2.7}{-3} &  4.3 & \num{5.9}{-4} &  4.4 & \num{2.2}{-3} &  4.2 & \num{8.5}{-3} &  5.2  \\ \hline  
 1/80      & \num{4.8}{-4} &  5.5 & \num{5.6}{-4} &  4.8 & \num{1.3}{-4} &  4.6 & \num{5.0}{-4} &  4.3 & \num{1.7}{-3} &  5.0  \\ \hline  
 1/160	   & \num{9.4}{-5} &  5.1 & \num{1.2}{-4} &  4.5 & \num{2.9}{-5} &  4.5 & \num{1.2}{-4} &  4.3 & \num{4.1}{-4} &  4.2  \\ \hline  
\rateLabel &  $2.39$       &      &  $2.19$       &      &  $2.17$       &      &  $2.10$       &      &  $2.26$       &       \\ \hline  
\end{tabular}
}

\bogus{
\begin{table}[hbt]\tableFont 
\begin{center}
\begin{tabular}{|l|c|c|c|c|c|c|c|c|c|c|c|} \hline 
grid  & N & $p$ & r & $\vert \mathbf{v} \vert$ & r & $\bar{\mathbf{u}}$ & r & $\bar{\mathbf{v}}$ & r & $\bar{\mathbf{s}}$& r \\ \hline 
                sq20 &     1 & \num{3.3}{-3} &      & \num{5.4}{-3} &      & \num{1.5}{-3} &      & \num{5.8}{-3} &      & \num{2.0}{-3} &       \\ \hline
                sq40 &     2 & \num{7.0}{-4} &  4.7 & \num{1.2}{-3} &  4.4 & \num{3.9}{-4} &  3.9 & \num{1.4}{-3} &  4.2 & \num{4.0}{-4} &  5.0  \\ \hline
                sq80 &     4 & \num{1.4}{-4} &  4.9 & \num{2.7}{-4} &  4.6 & \num{8.8}{-5} &  4.4 & \num{3.0}{-4} &  4.6 & \num{8.3}{-5} &  4.8  \\ \hline
               sq160 &     8 & \num{3.0}{-5} &  4.7 & \num{6.0}{-5} &  4.4 & \num{2.0}{-5} &  4.4 & \num{6.6}{-5} &  4.4 & \num{2.9}{-5} &  2.9  \\ \hline
    rate             &       &  $2.26$       &      &  $2.16$       &      &  $2.08$       &      &  $2.16$       &      &  $2.06$       &       \\ \hline   
\end{tabular}
\caption{Inses bulk: fw.mu02.rem1, solid=(ewe,fos), AM=1, amImp=0.5, vpf=0.1, ts=pc2, nc=1, normalMotion=0, bcs=periodic, bcsTop=displacement, cdv=0.25, $t=.3$, $\mu=.02$, ad2=0, $\rho=1$, $\bar{\rho}=.1$, , cfl=0.5, Tue Jul 16  4:58:04 2013}\label{table:bulk.fw.mu02.rem1.pc2.ewe.fos.vector}
\end{center}
\end{table}
}
\newcommand{\tableMPVEAmpLight}{%
\begin{tabular}{|c|c|c|c|c|c|c|c|c|c|c|} \hline
\multicolumn{11}{|c|}{MP-VE, traveling wave, light solid} \\ \hline 
\strutt~~$h_j$~~& $E_j^{(p)}$ & $r$ & $E_j^{(\vv)}$ & $r$ & $E_j^{(\usv)}$ & $r$  & $E_j^{(\vsv)}$  & $r$ & $E_j^{(\sigmasv)}$  & $r$ \\ \hline 1/20       & \num{3.3}{-3} &      & \num{5.4}{-3} &      & \num{1.5}{-3} &      & \num{5.8}{-3} &      & \num{2.0}{-3} &       \\ \hline   
 1/40       & \num{7.0}{-4} &  4.7 & \num{1.2}{-3} &  4.4 & \num{3.9}{-4} &  3.9 & \num{1.4}{-3} &  4.2 & \num{4.0}{-4} &  5.0  \\ \hline   
 1/80       & \num{1.4}{-4} &  4.9 & \num{2.7}{-4} &  4.6 & \num{8.8}{-5} &  4.4 & \num{3.0}{-4} &  4.6 & \num{8.3}{-5} &  4.8  \\ \hline   
 1/160	    & \num{3.0}{-5} &  4.7 & \num{6.0}{-5} &  4.4 & \num{2.0}{-5} &  4.4 & \num{6.6}{-5} &  4.4 & \num{2.9}{-5} &  2.9  \\ \hline   
\rateLabel  &  $2.26$       &      &  $2.16$       &      &  $2.08$       &      &  $2.16$       &      &  $2.06$       &       \\ \hline   
\end{tabular}
}

\begin{figure}[hbt]\tableFontSize
\begin{center}
\tableMPVEAmpHeavy
\vskip.5\baselineskip
\tableMPVEAmpMedium
\vskip.5\baselineskip
\tableMPVEAmpLight
\caption{Traveling wave solution for an inviscid incompressible fluid and elastic solid (model problem MP-VE).  Maximum errors and estimated convergence rates at $t=0.3$, computed using the AMP scheme for a heavy solid, $\delta=10^3$, medium solid, $\delta=1$ and light solid, $\delta=10^{-1}$.
}
\label{tab:MP-VE_AMP}
\end{center}
\end{figure}

%% file: texFiles/conclusions.tex
\section{Conclusions} \label{sec:conclusions}

We have described a stable added-mass partitioned (AMP) algorithm for coupling
incompressible flows with compressible elastic solids. 
The AMP algorithm is apparently the first stable partitioned scheme for incompressible
flows and compressible elastic bulk solids that does not require any
sub-iterations. The Robin 
coupling conditions for the AMP approach are derived from a local characteristic 
decomposition in the solid, and these define relationships between
the fluid velocity and traction on the interface in terms
of the outgoing characteristic in the solid. For fractional-step
incompressible flow algorithms, such as the one used in this paper, an alternative Robin interface
condition was derived that defines a mixed boundary condition on
the pressure; this can be used when solving the pressure equation.

A normal mode analysis for a model FSI problem showed that the AMP algorithm is stable
even for very light solids when the added-mass effects are large. An analysis of a 
traditional partitioned algorithm showed that the scheme 
is, in fact, formally unconditionally unstable for any ratio of the
solid to fluid density, no matter how large.  The traditional scheme may be stable
on a coarse grid but becomes unstable for a sufficiently fine grid.

We derived exact traveling wave solutions for a number of FSI model problems
and these were used to verify the stability and accuracy of the AMP algorithm. Numerical
results confirmed that the AMP scheme was stable for light, medium and heavy solids, and 
the scheme was second-order accurate in the maximum-norm.

The AMP algorithm was described for the case of FSI problems with small perturbations in the solution with respected to a fixed interface position.  In future work, we plan to incorporate the AMP approximations into our deforming
composite grid framework~\cite{fsi2012} that solves the full three-dimensional and nonlinear
incompressible Navier-Stokes equations and that supports large structural deformations.

%% file: texFiles/exactSolutions.tex
\section{Traveling wave exact solutions for the  FSI model problems} \label{sec:travelingWave}

\renewcommand{\usvh}{{\hat {\bf u}}}
\renewcommand{\ush}{{\hat u}}

In this section we present exact {\em traveling wave} solutions to the three FSI model
problems MP-IA, MP-VA and MP-VE defined in Section~\ref{sec:modelProblems}.
In each case we look for traveling wave solutions of the form
\begin{equation}
  \vv(x,y,t)=\vvh(y)\,\expkw,  \qquad p(x,y,t)=\ph(y)\,\expkw, \qquad \usv(x,y,t) = \usvh(y)\,\expkw,
\label{eq:travelingwave}
\end{equation}
for the fluid velocity, fluid pressure and solid displacement, respectively.
The solutions are assumed to be periodic in the horizontal coordinate $x$ with wave number $k$, and have frequency $\omega$ (possibly complex) in time $t$.  Since the fundamental period in the horizontal direction is $L$ we are interested in values of $k$ that are integer multiples of $2\pi/L$.
After substituting these forms into the governing equations, the result is a
system of ordinary differential equations for $\vvh(y)$, $\ph(y)$ and $\usvh(y)$ in terms of the vertical coordinate $y$ subject to homogeneous
boundary and interface conditions. This system has non-trivial solutions (eigenfunctions) provided $\omega$ and $k$ are solutions of a certain determinant condition which defines the dispersion relation for each model problem.

For the fluid solutions it is helpful to define
\begin{align*}
&   \Ck=\cosh(k H),\qquad \Sk=\sinh(k H),\qquad \Ca =\cosh(\alpha H),\qquad \Sa =\sinh(\alpha H),
\end{align*}
where
\[
\alpha^2=k^2- \frac{i\rho\omega}{\mu}.
\]
For the solid solutions we define
\begin{align*}
   C_a=\cosh(a\Hs), \qquad S_a=\sinh(a\Hs), \qquad C_b=\cosh(b\Hs), \qquad S_b=\sinh(b\Hs),
\end{align*}
where
\[
a^2=k^2-{\omega^2\over \cp^2}, \qquad b^2=k^2-{\omega^2\over \cs^2}.
\]

\subsection{Traveling wave solution for MP-IA}

Solutions of the system of ODEs and boundary conditions at $y=-H$ and $y=\Hs$ for the MP-IA model problem can be written in the form
\begin{equation}
\begin{array}{ll}
 \vh_1(y) = i \Af \cosh(k (y+H)), \qquad &\vh_2(y) = \Af \sinh(k (y+H)) , \medskip \\
 \displaystyle{
  \ph(y)  =  \frac{i\rho\omega}{k} \Af \cosh(k (y+H))} ,  &\ush_2(y) =  \As \sinh(a (y-\Hs)) ,
\end{array}
\label{eq:coefficients1}
\end{equation}
where $\Af$ and $\As$ are constants.  Application of the interface conditions at $y=0$ gives the matrix problem
\begin{align*}
    M_{{\rm IA}}  \begin{bmatrix} \Af \\ \As  \end{bmatrix} = 
  \begin{bmatrix}
      \Sk            &  -i\omega\Sas \\
     i\rho\omega\Ck & \rhos\cp^2 a k\Cas    &      
 \end{bmatrix} 
 \begin{bmatrix} \Af \\ \As  \end{bmatrix} = 0 .
\end{align*}
For non-trivial solutions we require values of $\omega$ and $k$ for which $\det(M_{{\rm IA}})=0$. This condition leads to the dispersion relation, 
\begin{equation}
  W_{{\rm IA}}(\omega,k) =  \rhos\cp^2 a k \Sk\Cas - \rho\omega^2 \Ck\Sas = 0 .
  \label{eq:dispersion1}
\end{equation}
Given a solution of~\eqref{eq:dispersion1} for $\omega$ and $k$ (usually values for $k\in\Real$ and the other parameters
are specified, and then~\eqref{eq:dispersion1} becomes a nonlinear equation for $\omega$), 
$\Af$ can be determined in terms of the free parameter $\As$ as 
\begin{align*}
   \Af = i\omega\frac{\Sas}{\Sk} \As. 
\end{align*}
We choose $\As$ so that the maximum displacement of the interface is equal to the real-valued amplitude parameter $\ampe$, 
i.e., $\vert \ush_2(0)\vert=\ampe$.  The real and imaginary parts of the solutions in~\eqref{eq:travelingwave} with component coefficient functions given in~\eqref{eq:coefficients1} define real solutions to the model problem MP-IA.

An analysis of the dispersion relation in~\eqref{eq:dispersion1} shows that there is a plus-minus pair of real solutions for $\omega$ satisfying $0<\vert\omega\vert<\cp\vert k\vert$ if $\rhos/\rho<(k\Hs)\coth(kH)$, and there are an infinite number of real plus-minus pairs satisfying $\vert\omega\vert>\cp\vert k\vert$.  In the limit of large $\vert\omega\vert$, it can be shown that these solutions have the asymptotic form $\omega \sim \pm n \pi \cp/\Hs$ for large integers $n$. 
For example, if $k=2\pi$, $\rhos=\lambdas=\mus=0.1$, $\Hs=0.5$, $\rho=1$ and $H=1$, then $\omega\approx3.36460699$ is the one positive solution satisfying $0<\vert\omega\vert<\cp\vert k\vert$, while $\omega\approx 15.5134370$ is the smallest positive solution satisfying $\vert\omega\vert>\cp\vert k\vert$.


\subsection{Traveling wave solution for MP-VA}

Solution of the system of ODEs and boundary conditions at $y=H$ for the viscous fluid domain of the model problem MP-VA can be written in the form
\begin{equation}
\begin{array}{l}
\displaystyle{
\hat v_1(y) =i\Af \Big(kS_\alpha\cosh(ky)+\Df\cosh(k(y+H))-\alpha S_k\cosh(\alpha y)\Big)
}\smallskip\\
\displaystyle{
\qquad\qquad +\frac{i\alpha\Bf}{k} \Big(k S_\alpha\cosh(ky)+\Df\cosh(\alpha(y+H))-\alpha S_k\cosh(\alpha y)\Big),
}\medskip\\
\displaystyle{
\hat v_2(y) =\Af \Big(kS_\alpha\sinh(ky)+\Df\sinh(k(y+H))-kS_k\sinh(\alpha y)\Big),
}\smallskip\\
\displaystyle{
\qquad\qquad  +\Bf \Big(\alpha S_\alpha\sinh(ky)+\Df\sinh(\alpha(y+H))-\alpha S_k\sinh(\alpha y)\Big) ,
}\medskip\\
\displaystyle{
\hat p(y)={i\rho\omega\over k}\Big[\Af \Big(kS_\alpha\cosh(ky)+\Df\cosh(k(y+H))\Big)+\Bf \alpha S_\alpha\cosh(ky)\Big],
}
\end{array}
\label{eq:fluidcoefficients}
\end{equation}
where
\[
  \Df = \alpha C_\alpha S_k-kC_kS_\alpha , 
\]
and $\Af$ and $\Bf$ are constants to be determined.  Solutions of the system of ODEs and the boundary condition at $y=-\Hs$ for the solid domain can be taken as $\ush_2(y)$ in~\eqref{eq:coefficients1} with constant $A_s$ to be determined.

The three constants $\Af$, $\Bf$ and $\As$ are determined by the three interface conditions for MP-VA. These leads to the
matrix equation
\begin{align*}
    M_{{\rm VA}}  \begin{bmatrix} \Af \\ \Bf \\ \As \end{bmatrix} = 0 ,
\end{align*}
where $M_{{\rm VA}}=[m_{ij}]\in\Complex^{3\times3}$ with
\begin{align*}
m_{11} &= i (k \Sa-\alpha \Sk+\Df \Ck)  ,\qquad
m_{12} = \frac{i}{k} \alpha (k \Sa-\alpha \Sk+\Df \Ca)  ,\qquad
m_{13} = 0  , \\
m_{21} &= \Df \Sk  ,\qquad
m_{22} = \Df \Sa  ,\qquad
m_{23} = -i \omega \Sas  , \\
m_{31} &= i \rho \frac{\omega}{k} (k \Sa+\Df \Ck)-2 \mu (k^2 \Sa-k \alpha \Sk+\Df k \Ck)   ,\\
m_{32} &= i \rho \frac{\omega}{k} \alpha \Sa-2 \mu (\alpha k \Sa-\alpha^2 \Sk+\Df \alpha \Ca)   ,\\
m_{33} &= \rhos\cp^2 a \Cas   .
\end{align*}
For nontrivial solutions we require $\det(M_{{\rm VA}})=0$, which leads to the dispersion relation
\begin{equation}
     W_{{\rm VA}}(\omega,k) = \det(M_{{\rm VA}}) = 0.
\label{eq:dispersion2}
\end{equation}
The precise form for $W_{{\rm VA}}(\omega,k)$ is messy and not particularly revealing, and so we suppress the details here.  (The form can be found readily using the given components of $M_{\rm VA}$ and the symbolic program Maple, for example, and then solved numerically.)  Once a solution $(\omega,k)$ is found, the constants $\Af$, $\Bf$ can be determined in terms of $\As$ from
\begin{align*}
 \begin{bmatrix}
    m_{11} & m_{12}  \\
    m_{21} & m_{22} 
 \end{bmatrix}
  \begin{bmatrix}  \Af \\ \Bf \end{bmatrix}
  = 
 -\As \begin{bmatrix} m_{13} \\ m_{23}  \end{bmatrix} .
\end{align*}
The parameter $\As$ is an arbitrary constant and can be chosen, for example, so that $\vert\ush_2(0)\vert=\ampe$.

For the dispersion relation in~\eqref{eq:dispersion2} we have investigated solutions numerically.  For a chosen set of parameters of the problem, it appears that there are an infinite number of complex-valued solutions for $\omega$.  For example, if $k=2\pi$, $\rhos=\lambdas=\mus=0.1$, $\Hs=0.5$, $\rho=1$, $\mu=.02$ and $H=1$, then one solution is given by the complex number $\omega\approx(2.79247701,-0.746859802)$.




\subsection{Traveling wave solution for MP-VE}

For the case of the model problem MP-VE we use the solutions in~\eqref{eq:fluidcoefficients} for the fluid.  As noted earlier, these functions satisfy the system of ODEs in the fluid and the boundary condition at $y=H$.  For the elastic solid, we use
\begin{equation}
\begin{array}{l}
\displaystyle{
\hat u_1 =\As \Big(-k^2 S_b\cosh(ay)+\Ds\cosh(a(y-\bar H))+abS_a\cosh(by)\Big)
} \smallskip\\
\displaystyle{
\qquad +\Bs\Big(-k^2 S_b\cosh(ay)+\Ds\cosh(b(y-\bar H))+abS_a\cosh(by)\Big),
}\medskip\\
\displaystyle{
\hat u_2={a\As \over i k}\Big(-k^2S_b\sinh(ay)+\Ds\sinh(a(y-\bar H))+k^2S_a\sinh(by)\Big)
}\smallskip\\
\displaystyle{
\qquad +{k\Bs\over i b}\Big(-abS_b\sinh(ay)+\Ds\sinh(b(y-\bar H))+abS_a\sinh(by)\Big),
}
\end{array}
\label{eq:solidcoefficients}
\end{equation}
where
\[
  \Ds = k^2C_aS_b-abC_bS_a .
\]
The four constants $\Af$, $\Bf$, $\As$ and $\Bs$ are determined by the four interface conditions. These leads to the
matrix equation
\begin{align*}
    M_{{\rm VE}}  \begin{bmatrix} \Af \\ \Bf \\ \As \\ \Bs \end{bmatrix} = 0 ,
\end{align*}
where $M_{{\rm VE}}=[m_{ij}]\in\Complex^{4\times4}$ has components
\begin{align*}
m_{11} &= i (k \Sa-\alpha \Sk+\Df \Ck)  ,\quad
m_{12} =  i\alpha (k \Sa-\alpha \Sk+\Df \Ca)/k , \\ 
m_{13} &= i \omega (-k^2 \Sbs/\Ds+\Cas+a b \Sas/\Ds)  ,\quad 
m_{14} = i \omega (-k^2 \Sbs/\Ds+a b \Sas/\Ds+\Cbs) ,  \\
m_{21} &= \Df \Sk  ,\quad
m_{22} = \Df \Sa  ,\quad
m_{23} = -\omega a/(k \Sas)  ,\quad 
m_{24} = -\omega k/(b \Sbs) ,  \\
m_{31} &= -2 i \mu k \Df \Sk  ,\quad
m_{32} = -\mu (i\alpha^2 \Df \Sa/k+i k \Df \Sa) ,\quad
m_{33} = -2 \mus a \Sas  ,\quad
m_{34} = \mus (-b \Sbs-k^2/b \Sbs),   \\
m_{41} &= i \rho\omega (k \Sa+\Df \Ck)/k-2 \mu (k^2 \Sa-k \alpha \Sk+\Df k \Ck), \quad
m_{42} = i \rho \omega \alpha \Sa/k-2 \mu (\alpha k \Sa-\alpha^2 \Sk+\Df \alpha \Ca) ,  \\
m_{43} &= -i (\lambdas+2 \mus) (-a^2 k \Sbs/\Ds+a^2 \Cas/k+k a b \Sas/\Ds)+i \lambdas k (-k^2 \Sbs/\Ds+\Cas+a b \Sas/\Ds) ,  \\
m_{44} &= -i (\lambdas+2 \mus) (-a^2 k \Sbs/\Ds+k a b \Sas/\Ds+k \Cbs)+i \lambdas k (-k^2 \Sbs/\Ds+a b \Sas/\Ds+\Cbs)  .
\end{align*}
Non-trivial solutions are obtained if $\omega$ and $k$ satisfy the dispersion relation
\begin{align*}
     W_{{\rm VE}}(\omega,k) = \det(M_{{\rm VE}}) = 0.
\end{align*}
For a given solution of $W_{{\rm VE}}(\omega,k)=0$, the constants $\Af$, $\Bf$, $\As$ can be defined in terms
of $\Bs$ from 
\begin{align*}
 \begin{bmatrix}
    m_{11} & m_{12} & m_{13} \\
    m_{21} & m_{22} & m_{23} \\
    m_{31} & m_{32} & m_{33} 
 \end{bmatrix}
  \begin{bmatrix}  \Af \\ \Bf \\ \As \end{bmatrix}
  = 
  - \Bs \begin{bmatrix} m_{14} \\ m_{24} \\ m_{34} \end{bmatrix} .
\end{align*}
In this construction, the parameter $\Bs$ is an arbitrary constant. We choose $\Bs$ so that the maximum displacement of the interface
is $\ampe$, i.e. $\vert\usvh(0)\vert=\sqrt{\ush_1^2(0)+\ush_2^2(0)}=\ampe$. 
The real and imaginary parts of~\eqref{eq:travelingwave} with~\eqref{eq:fluidcoefficients} and~\eqref{eq:solidcoefficients} define real solutions to the model problem MP-VE.

As in the previous case, we investigate solutions of $W_{{\rm VE}}(\omega,k)=0$ numerically.  For example, if $k=2\pi$, $\rhos=\lambdas=\mus=0.1$, $\Hs=0.5$, $\rho=1$, $\mu=.02$ and $H=1$, then one solution is given by $\omega\approx(1.90532196,-0.652436711)$.


%% file: texFiles/procedures.tex
\floatname{algorithm}{Procedure}
\algnewcommand{\LineComment}[1]{\State \(\triangleright\) #1}
\section{AMP time-stepping procedures} \label{sec:timeSteppingProcedures}

This section provides partial implementations of some of the procedures that appear 
in the AMP algorithm described in Section~\ref{sec:timeSteppingAlgorithm}.
Let $\grad_h$ and $\Delta_h$ denote discrete approximations to the operators $\grad$ and $\Delta$, respectively, and let 
$\Lv_\iv^n$ be the discretized form of the right-hand-side to the fluid momentum equation defined by 
\begin{align}
   \Lv_\iv^n \equiv - \frac{1}{\rho} \grad_h p_\iv^n + \frac{\mu}{\rho} \Delta_h \vv_\iv^n . 
\end{align}
Let $\OmegaFh$ denote the set of grid point indices, $\iv$, corresponding to the interior, boundary and interface points of the fluid domain.
Similarly, let $\OmegaSh$ denote the corresponding indices in the 
solid domain, and let $\GammaIh$ denote the set of indices on the interface.

\begin{algorithm}
\caption{Advance the fluid velocity from time $t^{n-1}$ to $t^{n}$, return a predicted velocity $\vv_\iv^{(p)}$. This routine assigns interior, boundary and interface points.}
\begin{algorithmic}[1]
\Procedure{{\blue\tt advanceFluid}}{~$\vv_\iv^{n-1}$, $p_\iv^{n-1}$, $\vv_\iv^{n-2}$, $p_\iv^{n-2}$~}
  \State $\displaystyle%
    \vv_\iv^{(p)} = \vv_\iv^{n-1} + \dt\Big( \frac{3}{2} \Lv_\iv^{n-1} - \frac{1}{2} \Lv_\iv^{n-2}\Big) ,  \qquad \iv\in\OmegaFh
$   \Comment{(Adams-Bashforth).}
  \State \Return{$\vv_\iv^{(p)}$}
\EndProcedure
\end{algorithmic}
\end{algorithm}

\begin{algorithm}
\caption{Advance the fluid velocity from $t^{n-1}$ to $t^{n}$ given an approximation to the state at $t^{n}$. This routine assigns interior, boundary and interface points.}
\begin{algorithmic}[1]
\Procedure{{\blue\tt correctFluid}}{~$\vv^{(p)}_\iv$, $p_\iv^{(p)}$, $\vv_\iv^{n-1}$, $p_\iv^{n-1}$~}
  \State $\displaystyle
    \vv_\iv^{n} = \vv_\iv^{n-1} +  \frac{\dt}{2}\left( \Lv_\iv^{(p)} +  \Lv_\iv^{n-1} \right),   \qquad \iv\in\OmegaFh
$   \Comment{(trapezoidal rule).}
  \State \Return{$\vv_\iv^{n}$}
\EndProcedure
\end{algorithmic}
\end{algorithm}

\begin{algorithm}
\caption{Solve for the pressure at all grid points.}
\begin{algorithmic}[1]
\Procedure{{\blue\tt solveFluidPressureEquation}}{~$\vv_\iv$, $\vsv_\iv$, $\sigmasv_\iv\nv_\iv$, $\dot{\vvs}_\iv$~}
  \LineComment{Solve the following system of equations for the pressure:}
  \State $\Delta_h p_\iv =0, \qquad \iv\in\OmegaFh $  \label{eq:pDiscrete}
  \State $\displaystyle
      -p_\iv  - \frac{\zp\dt}{\rho} \nv_\iv\cdot\grad_h p_\iv 
   = - \nv_\iv^T\tauv_\iv\nv_\iv + \frac{\mu\zp\dt}{\rho} \nv_\iv^T(\grad_h\times\grad_h\times\vv_\iv) 
     + \nv_\iv^T\sigmasv_\iv\nv_\iv +\zp\dt \, \nv_\iv^T\dot{\vvs}_\iv 
            ,\quad \iv\in\GammaIh
$  \label{eq:pbcDiscrete}
  \State \Return{$p_\iv$}
\EndProcedure
\end{algorithmic}
Notes: Line~\ref{eq:pbcDiscrete} is the Robin boundary condition for the pressure~\eqref{eq:AddedMassPressureBC}.
  In practice it is useful to add a divergence damping term to the right-hand-side of the the pressure equation
  on line~\ref{eq:pDiscrete} following~\cite{ICNS,splitStep2003}.
\end{algorithm}

\begin{algorithm}
\caption{Project the interface velocity (if $\beta=0$ do not include traction terms).}
\begin{algorithmic}[1]
\Procedure{{\blue\tt projectInterfaceVelocity}}{~$\vv_\iv$, $p_\iv$, $\vsv_\iv$, $\sigmasv_\iv\nv_\iv$, $\beta$~}
  \State $ \displaystyle
\nv_\iv^T\vv_\iv^I = \frac{\zf}{\zf+\zp} \nv_\iv^T\vv_\iv + \frac{\zp}{\zf+\zp} \nv_\iv^T\vvs_\iv 
                   + \frac{\beta}{\zp + \zf} \Big( \nv_\iv^T\sigmavs_\iv\nv_\iv + p_\iv -\nv_\iv^T\tauv_\iv\nv_\iv  \Big) ,   \quad \iv\in\GammaIh 
  $  
  \State
$ \displaystyle
\tanv_m^T\vv_\iv^I = \frac{\zf}{\zf+\zs} \tanv_m^T\vv_\iv + \frac{\zs}{\zf+\zs} \tanv_m^T\vvs_\iv 
      + \frac{\beta}{\zs + \zf} \Big( \tanv_m^T\sigmavs_\iv\nv_\iv  -\tanv_m^T\tauv_\iv\nv_\iv  \Big) , \quad m=1,2,  \quad \iv\in\GammaIh 
$ 
  \State \Return{$\vv_\iv^I$}
\EndProcedure
\end{algorithmic}
\end{algorithm}

\begin{algorithm}[H]
\caption{Apply the boundary conditions on the fluid velocity. Return the interface velocity and traction.}
\begin{algorithmic}[1]
\Procedure{{\blue\tt assignFluidVelocityBoundaryConditions}}{~$\vv_\iv$, $p_\iv$, $\vsv_\iv$, $\sigmasv_\iv\nv_\iv$~}
  \LineComment{Solve the following equations to determine $\vv_\iv$ on the boundary and ghost line.}
  \State $\displaystyle
  \rho\tanv_m^T( \vv_\iv-\vv_\iv^{n})/\dt +  (\mu/2)\tanv_m^T\big(\Delta_h\vv_\iv + \Delta_h\vv_\iv^{n}\big)
         = (1/2)\tanv_m^T\big( \grad_h p_\iv + \grad_h p_\iv^{n}\big)  , 
          \quad \iv\in\GammaIh, $ \label{eq:EqnOnBoundary}
  \State $\displaystyle \zs \tanv_m^T\vv_\iv + \mu\tanv_m^T(\grad_h\vv_\iv+ (\grad_h\vv_\iv)^T )\nv_\iv =\tanv_m^T\sigmasv_\iv\nv_\iv+ \zs\tanv_m^T\vsv_\iv ,
       \quad \iv\in\GammaIh  $ \Comment{(Eqn~\eqref{eq:solidVtangentialVsubproblem})} \label{eq:AMPtv}
  \State $  \grad_h\cdot\vv_\iv = 0 ,  \quad\iv\in\GammaIh, $ 
  \LineComment{Determine the new interface values.}
  \State $ \vv_\iv^I=\vv_\iv, \quad (\sigmav\nv)^I_\iv = -p_\iv \nv_\iv + \tauv_\iv\nv_\iv,  \quad \iv\in\GammaIh$ 
  \State \Return{$(\vv_\iv^I, (\sigmav\nv)^I_\iv)$}
\EndProcedure
\end{algorithmic}
Note: When the velocity is advanced treating $\mu\Delta\vv$ explicitly in time,
the AMP boundary condition~\eqref{eq:solidVtangentialVsubproblem} (line~\ref{eq:AMPtv}) is combined
with the interior equation on the boundary (line~\ref{eq:EqnOnBoundary}) to determine values on the boundary
and first ghost line.  This is done since for {\em small} $\zs$,
equation~\eqref{eq:solidVtangentialVsubproblem} can be thought of as primarily
determining the ghost values, while for {\em large} $\zs$
equation~\eqref{eq:solidVtangentialVsubproblem} primarily sets the boundary
value. When $\mu\Delta\vv$ is treated implicitly in time,
equation~\eqref{eq:solidVtangentialVsubproblem} is incorporated directly into
the implicit system, to the same effect.
\end{algorithm}